\documentclass[11pt,a4paper]{article}
\usepackage{bbm}
\usepackage{amsmath}
\usepackage{amsfonts}
\usepackage{amssymb}
\usepackage{graphicx}
\usepackage{booktabs}
\usepackage{tabularx}
\usepackage{lscape,url}

\newtheorem{te}{Theorem}
\newtheorem{lemma}[te]{Lemma}
\newtheorem{remark}[te]{Remark}

\newtheorem{example}[te]{Example}
\newtheorem{question}[te]{Question}
\newtheorem{corollary}[te]{Corollary}
\newtheorem{conjecture}[te]{Conjecture}
\newenvironment{proof} {\par \noindent \textbf{Proof: }}{\QED\par \bigskip \par}
\newcommand{\QED}{\hfill$\square$}
\newcommand{\rz}{\vspace{0.1cm}}

\title { \bigskip
    \bf {Chordality and hyperbolicity of a  graph}
    \thanks{Supported by  Science and Technology Commission of Shanghai Municipality  (No. 08QA14036 and  No. 09XD1402500),   Chinese Ministry
of Education (No. 108056), and  National Natural Science Foundation
of China (No. 10871128). } }

\author
{
{\large \sc Yaokun Wu\footnotemark[3],  Chengpeng Zhang  } \\
{\em \normalsize Department of  Mathematics, Shanghai Jiao Tong University} \\
{\em \normalsize   800 Dongchuan Road, Shanghai, 200240, China }
\\   {\normalsize E-mail:
{  \tt \{ykwu, cpzhang\}@sjtu.edu.cn }   }\\
 }
\date{}

\begin{document}

\maketitle

\vspace{-0.5cm}

\begin{abstract} Let $G$ be a connected graph with the usual
shortest-path metric $d$.    The graph  $G$ is $\delta$-hyperbolic
provided for any vertices $x,y,u,v$ in it, the two larger of the
three sums $d(u,v)+d(x,y),d(u,x)+d(v,y)$ and $d(u,y)+d(v,x)$ differ
by at most $2\delta.$  The graph  $G$ is  $k$-chordal provided it
has no induced cycle of length greater than $k.$   Brinkmann, Koolen
and Moulton find that every $3$-chordal graph is $1$-hyperbolic and
is not  $\frac{1}{2}$-hyperbolic if and only if it contains one of
two special graphs as an isometric subgraph. For every $k\geq 4,$ we
show that a $k$-chordal graph must be $\frac{\lfloor
\frac{k}{2}\rfloor}{2}$-hyperbolic  and there does exist a
$k$-chordal graph which is not   $\frac{\lfloor \frac{k-2}{2}\rfloor
}{2}$-hyperbolic. Moreover, we   prove that
        a $5$-chordal graph is $\frac{1}{2}$-hyperbolic if and only if it
does not contain any of a list
  of six special   graphs    (See Fig.   \ref{fig0})  as an isometric subgraph.
\end{abstract}


        \footnotetext[3]
        {        Fax:
86-21-54743152
       }

\section{Introduction}

\subsection{Tree-likeness}

Trees are graphs with some very distinctive and fundamental properties and it is legitimate to ask to what degree those properties can be transferred to more general structures that are tree-like in some sense \cite[p. 253]{Diestel}.
Roughly speaking, tree-likeness stands for  something related to low
 dimensionality, low complexity, efficient information deduction (from local to global),  information-lossless  decomposition (from global into simple pieces)
  and nice shape for efficient implementation of  divide-and-conquer strategy. For the very basic
 interconnection structures like a graph or a  hypergraph,
 tree-likeness is
  naturally reflected by the strength of interconnection, namely its
  connectivity/homotopy type or cyclicity/acyclicity, or just the degree of derivation
  from
  some characterizing conditions of a tree/hypertree and its various associated structures and generalizations.

In vast applications,  one finds that the borderline between
tractable and intractable cases
 may be the tree-like degree of
 the structure to be
dealt with \cite{CYY}.  A support to this   from the fixed-parameter
complexity  point of view is the observation that on various
tree-structures we can design very good algorithms for many purposes
and these algorithms can somehow be lifted to tree-like structures
\cite{ALS,Downey,DMC,KFL}. It is thus very useful to get information
on approximating general structures by tractable structures, namely
tree-like structures. On the other hand, one not only finds it
natural that tree-like structures appear extensively in many fields,
say biology \cite{Dress96}, structured programs \cite{Thorup}  and
database theory \cite{Fagin}, as graphical representations of
various types of hierarchical relationships, but also
 notice surprisingly that many practical structures we
encounter are just tree-like, say the internet
\cite{ABKMRT,Kleinberg,ST} and chemical compounds  \cite{YAM}. This
prompts in many areas the very active study of tree-like structures.
Especially, lots of ways to define/measure a tree-like structure
have been proposed in the literature from many different
considerations, just to name a few, say tree-width \cite{RS,RS86},
tree-length \cite{DG,UY}, combinatorial dimension
\cite{Dress96,Dress84}, $\epsilon$-three-points condition
\cite{Dinitz},
 $\epsilon$-four-points condition \cite{ABKMRT},  asymptotic connectivity
 \cite{Bahls},
 tree-partition-width \cite{Bod,Wood}, tree-degree \cite{CM},
McKee-Scheinerman chordality \cite{McS}, $s$-elimination dimension
\cite{DKT}, linkage (degeneracy) \cite{DKT,KT,LW}, sparsity order
\cite{Laurent},
 persistence \cite{Downey}, cycle rank \cite{CYY,Lee}, various
degrees of acyclicity/cyclicity \cite{Duke,Fagin},  boxicity
\cite{Roberts}, doubling dimension  \cite{Gupta},  Domino treewidth
\cite{Bod}, hypertree-width \cite{GLS}, coverwidth \cite{Chen},
spread-cut-width \cite{CJG}, Kelly-width \cite{Hunter},
 and many other width parameters
  \cite{DMC,HOSG}.   It is clear that many relationships among
  these concepts should be expected as they are all formulated in different ways to represent
     different aspects of our vague but intuitive idea of tree-likeness.  To clarify these relationship helps to bridge the study
  in different fields focusing on different tree-likeness measures and helps to improve our understanding
  of the  universal  tree-like world.

As a small step in pursuing further understanding of tree-likeness,
we take up in this paper the modest task of comparing two parameters
of tree-likeness,
  namely (Gromov) hyperbolicity and chordality of a graph.
    Our main result is that  $k$-chordal graphs must be $\frac{\lfloor
\frac{k}{2}\rfloor}{2}$-hyperbolic when $k\geq 4$  (Theorem
\ref{main}). Besides that, we determine a complete set of
unavoidable
   isometric subgraphs of $5$-chordal  graphs attaining hyperbolicity $1$  (Theorem \ref{main1}),   as a minor
    attempt to   respond to the general  question,
   ``what is the structure of graphs with relative small hyperbolicity'' \cite[p. 62]{BKM01}, and the even more general
   question,
 ``what is the structure of a very tree-like graph''.

The plan of the paper is as follows.  In Sections \ref{chordal}
 and  \ref{hyper} we introduce the two tree-likeness parameters,
 chordality and hyperbolicity, respectively.
   Section \ref{theorem} is devoted to a general discussion of the
   relationship between chordality and hyperbolicity, including a
   presentation of our main results (Theorems \ref{main}  and
   \ref{main1}).
Some  consequences of our main results will be listed in Section
 \ref{oct}.
 In
Section  \ref{parameter}, we study the relationship between several
other tree-likeness parameters and  the main objects of this paper,
that is to say, chordality and hyperbolicity, and make use of these
relationship to connect chordality and hyperbolicity. The
relationship between chordality and hyperbolicity thus obtained by
now  is not as strong as Theorem \ref{main}. But the discussion may
be of some independent interest.
 We present a complete and
self-contained proof of Theorems
 \ref{main} and \ref{main1}   in Section    \ref{proofs} in two stages: some preliminary facts are prepared in Section
 \ref{lemmas} and the final proof appears in Section \ref{Proof}.     Following \cite{BKM01,KM02}, the key to our work is to
  examine the extremal local configurations as described by
  Assumptions I and II (see Section \ref{lemmas}).
Several key   lemmas in Section \ref{lemmas} are basically copied
from \cite{BKM01,KM02}.   It is  often that these lemmas are to be
found as pieces of a long proof of a big statement in
\cite{BKM01,KM02} and so the validity of these technical lemmas
under some weaker assumptions needs to be carefully checked. We
include the complete proofs of them, more or less as they were
presented in \cite{BKM01,KM02}, not only for the convenience of the
reader but also to convince the reader that they do hold in our
setting.

\subsection{Chordality}\label{chordal}

We only consider  simple, unweighted, connected, but not necessarily
finite  graphs.   Any graph $G$  together with      the usual
shortest-path metric on it, $d_G: \ V(G)\times V(G)\mapsto
\{0,1,2,\ldots\}$,  gives rise to a metric space. We often suppress
the subscript and write  $d(x,y)$ instead of   $d_G(x,y)$ when the
graph is known by context.
  Moreover, we
may   use the shorthand $xy$ for $d(x,y)$ to further simplify  the
notation. Note that a pair of vertices $x$ and $y$ form an edge if
and only if $xy=1.$            For $S,T\subseteq V(G)$, we write
$d(S,T)$ for $\min _{x\in S,y\in T}d(x,y)$.  We often omit the
brackets and adopt the convention that  $x$ stands  for the
singleton set $\{x\}$ when no confusion can be caused.

Let  $G$  be a graph.  A  {\em walk of length $n$} in $G$ is a
sequence of vertices
 $x_0,x_1,x_2,\ldots ,x_n$ such that $x_{i-1}x_{i}=1$ for  $i=1,\ldots ,n$. If these $n+1$ vertices are
  pairwise different, we call the sequence
 a {\em path of length $n$}.  A {\em pseudo-cycle} of length $n$ in $G$ is a
cyclic sequence of $n$   vertices $x_1,\ldots ,x_n\in V(G)$ such
that $x_ix_{j}=1$ whenever $j=i+1 \ (\bmod \ n)$; we will reserve
the notation $[x_1x_2\cdots x_n]$ for this pseudo-cycle.      We
call this pseudo-cycle an {\em $n$-cycle}, or a {\em cycle of length
$n$}, if $x_1,\ldots,x_n$ are $n$ different vertices.
   A {\em chord} of a path or cycle is an
edge joining nonconsecutive vertices on  the path or cycle.   An
{\em odd chord} of a cycle of even length is a chord connecting
different vertices the distance between which in the cycle is odd.
 A cycle without chord is called an
{\em induced cycle}, or a {\em chordless cycle}.  For any $n\geq 3,$
the {\em $n$-cycle graph} is the  graph with $n$ vertices
          which   has a chordless   $n$-cycle and we denote this graph by $C_n$.   A subgraph $H$ of a
graph $G$ is \textit{isometric} if for any $u,v\in V(H)$   it holds
$d_H(u,v)=d_G(u,v)$. A $4$-cycle of a graph $G$ is an
  {\em isometric $4$-cycle}
  provided   the subgraph of   $G$   induced by the vertices of this
cycle  is isometric and the subgraph has only those four edges which
are displayed in the cycle. Indeed, this amounts to saying that this
cycle  is an induced/chordless cycle; c.f.  Lemma  \ref{EASY}.

We say that a graph is {\em $k$-chordal} if it does not contain any
induced $n$-cycle for $n>k.$ Clearly,    trees are nothing but
$2$-chordal graphs.  A $3$-chordal graph is usually termed as a {\em
chordal
 graph} and a $4$-chordal graph is often called a  {\em
hole-free
 graph}.    The class of   $k$-chordal graphs is also
discussed under the name $k$-bounded-hole graphs \cite{Gavril}.

 The  {\em chordality}   of   a   graph  $G$ is   the smallest integer $k\geq 2$  such that $G$  is  $k$-chordal
 \cite{BT1}.           Following   \cite{BT1},  we use the notation  $\mathbbm{l}\mathbbm{c}(G)$ for this parameter as it is merely the length of
 the longest chordless cycle in $G$ when $G$ is not a tree.         Note that our use of the concept of chordality is
  basically the same as that used in \cite{CLS,CR} but is very different
 from  the usage  of this term in  \cite{McS}.

The recognition of  $k$-chordal graphs is coNP-complete for
$k=\Theta (n^{\epsilon})$ for any constant $\epsilon >0$
\cite{Uehara}. Especially, to determine the chordality of the
hypercube is attracting much attention  under the name of the
snake-in-the-box problem   due to its connection with some
error-checking codes problem \cite{Klee}. Just like the famous
snake-in-the-box problem, it looks hard to determine the exact value
of the chordality of general grid graphs -- it is only easy to see
that $\mathbbm{l}\mathbbm{c}(G_{m,n})$   should be roughly
proportional to $nm$ when $\min (n,m)>2.$
 Nevertheless, just like
many other tree-likeness parameters, quite a few  natural graph
classes are known to have small chordality \cite{BLS}. We review
  some $5$-chordal  ($4$-chordal)  graphs  in the remainder of this subsection.

An {\em asteroidal triple} ($AT$) of a graph $G$ is a a set of three
vertices of $G$  such that for any pair of them there is a path
connecting the two vertices whose distance to the remaining vertex
is at least two.  A graph is {\em AT-free} if no three vertices form
an $AT$    \cite[p. 114]{BLS}. Obviously, all $AT$-free graphs are
$5$-chordal.  A graph is an {\em interval graph} exactly when it is
both chordal and $AT$-free \cite[Theorem 7.2.6]{BLS}. $AT$-free
graphs also include {\em cocomparability graphs} \cite[Theorem
7.2.7]{BLS}; moreover, all {\em bounded
 tolerance graphs}
are cocomparability
 graphs    \cite{GMT}  \cite[Theorem 2.8]{MA}  and a graph is a   {\em permutation graph} if and only if
  itself and its complement are cocomparability graphs \cite[Theorem 4.7.1]{BLS}.
An important subclass of cocomparability graphs is the class of
 {\em threshold graphs}, which are those graphs without any induced
subgraph isomorphic to the $4$-cycle, the complement of the
$4$-cycle or the path of length $3$  \cite[p. 23]{MA}.

 A graph is {\em weakly chordal } \cite{GMT,Hayward}  when both itself and
 its complement are $4$-chordal.
  Note that all tolerance graphs  \cite{MA}  are domination graphs \cite{Rusu}   and all domination graphs
 are weakly chordal  \cite{DHMO}.
 A graph is {\em strongly chordal} if it is chordal and if every even
 cycle of length at least $6$ in this graph has an odd chord  \cite[p.  21]{GMT}.
    A
graph is {\em distance-hereditary} if each of its induced paths, and
hence each of its connected induced subgraphs, is isometric
\cite{Howorka}.       We call a graph a {\em cograph} provided it
does not contain any induced  path of length $3$ \cite[Theorem
11.3.3]{BLS}. It is easy to see that each cograph is
distance-hereditary and all distance-hereditary
 graphs form a proper subclass
of  $4$-chordal graphs. It is also known that cocomparability graphs
are all $4$-chordal \cite{BT1, Gallai}.

\subsection{Hyperbolicity}\label{hyper}

\subsubsection{Definition  and background}

For any vertices  $x,y,u,v$ of a graph $G,$ put $\delta
_G(x,y,u,v)$, which we often abbreviate to $\delta (x,y,u,v)$, to be
 the difference between the largest and the
second largest of the following three terms:
$$\frac{uv+xy}{2},\, \frac{ux+vy}{2},\, \text{and} \  \frac{uy+vx}{2}.$$ Clearly, $\delta
(x,y,u,v)=0$ if $x,y,u,v$ are not four different vertices. A graph
$G$, viewed as a metric space as mentioned above, is {\em
$\delta$-hyperbolic}  (or tree-like with defect at most $\delta$)
provided for any vertices $x,y,u,v$ in $G$ it holds $\delta
(x,y,u,v)\leq \delta$ and the  (Gromov) {\em hyperbolicity} of $G$,
denoted $\delta^* (G)$, is the minimum half integer $\delta$ such
that $G$ is $\delta$-hyperbolic
\cite{Bow91,BH,CDEHV,CDEHVX,DD,Gromov}. Note that it may happen
$\delta ^*(G)=\infty$.   But for a finite graph $G$,  $\delta^* (G)$
is clearly finite and polynomial time computable.

Note that in some earlier literature the concept of Gromov
hyperbolicity is used
 a little bit different from what we adopt here;
 what we call
$\delta$-hyperbolic here is called $2\delta$-hyperbolic in
\cite{ABKMRT,BC,BC08,BKM01,CE,DHHKMW,Dress96,GL,KM02,MS} and hence
the hyperbolicity of a graph is always an integer according to their
definition. We also refer to \cite{Alonso,Bow91,BH,  Vai} for some
equivalent and very accessible definitions of Gromov hyperbolicity
which involve some other comparable parameters.

The concept of hyperbolicity comes from the work of Gromov in
geometric group theory  which encapsulates many of the global
features of the geometry of complete, simply connected manifolds of
negative curvature \cite[p. 398]{BH}.    This concept   not only
turns out to be strikingly useful in coarse geometry but also
becomes more and more important in many applied fields like
networking and phylogenetics
\cite{CDEHV1,CDEHV,CDEHVX,CE,DraganX,Dress84,DHHKMW,DHM,Dress96,GL,JLB,JLHB,Kleinberg,ST}.
The hyperbolicity of a graph is a way to measure the additive
distortion with which every four-points sub-metric of the given
graph metric embeds into a tree metric \cite{ABKMRT}. Indeed, it is
not hard to check that the hyperbolicity of a tree  is zero -- the
corresponding condition for this  is known as the
  four-point condition (4PC) and is a characterization of
  general tree-like metric spaces
\cite{Dress84,Dress96,Imrich}. Moreover, the fact that hyperbolicity
is a tree-likeness parameter
      is   reflected in the easy fact that the hyperbolicity of a graph is the maximum  hyperbolicity of its 2-connected components --
      This observation implies the classical result    that
$0$-hyperbolic graphs are exactly block graphs, namely those  graphs
in which every $2$-connected subgraph is complete, which are also
known to be   those diamond-free chordal graphs
\cite{BM,DMS,Howorka}.
  More results on bounding hyperbolicity of graphs and characterizing low hyperbolicity graphs  can be found in
\cite{BC,BC08,BKM01,CDEHV1,CDEHV,DG,KM02}; we will only  report  in
Section \ref{MR} some  work most closely related to ours and refer
the readers to corresponding references for many other interesting
unaddressed work.

For any vertex $u\in V(G)$, the {\em Gromov product}, also known as
  the {\em overlap function}, of any two vertices $x$ and $y$ of $G$ with respect to $u$ is equal
to $\frac{1}{2}(xu + yu - xy)$  and is denoted by $(x\cdot y)_u$
\cite[p. 410]{BH}.  As an important context in phylogenetics
\cite{DHHKMW,DHM,Farris},  for any real number $\rho$, the {\em
Farris transform} based at $u$, denoted    $D_{\rho ,u}$, is the
transformation which sends $d_G$ to the map
$$D_{\rho , u}(d_G): V(G)\times V(G)\rightarrow \mathbb{R}:  \ (x,y)\mapsto \rho -(x\cdot
y)_u.$$
   We
say that $G$ is {\em $\delta$-hyperbolic with respect to $u\in
V(G)$} if the following inequality \begin{equation}(x\cdot y)_u\geq
\min ((x\cdot v)_u,(y\cdot v)_u)-\delta \label{EQ}\end{equation}
holds for any vertices $x,y,v$ of $G.$ It is easy to check that the
inequality \eqref{EQ} can be rewritten as
$$xy+uv\leq \max (xu+yv, xv+yu)+2\delta$$ and so we see that $G$ is
$\delta$-hyperbolic if and only if  $G$ is $\delta$-hyperbolic with
respect to every vertex of $G.$
 By a simple but nice argument, Gromov shows that  $G$ is   $2\delta$-hyperbolic provided
 it is  $\delta$-hyperbolic with respect to
 any given vertex \cite[Proposition 2.2]{Alonso} \cite[1.1B]{Gromov}.

The {\em tree-length}   \cite{Dour, DG, Lo, UY} of a graph $G$,
denoted  $\mathbbm{t}\mathbbm{l}(G)$, is the minimum integer $k$
such that there is a chordal graph $G'$ satisfying $V(G)=V(G')$,
$E(G)\subseteq E(G')$ and  $\max (d_G(u,v):\ d_{G'}(u,v)=1)= k.$  We
use the convention that the tree-length of a graph without any edge
is 1.
 It is straightforward
from the definition that   chordal graphs are exactly the graphs of
tree-length $1$.   It is also known that  $AT$-free graphs  and
distance-hereditary graphs have tree-length at most $2$ \cite[p.
367]{Dour}; a way to see this is to use the forthcoming result
relating chordality and tree-length as well as    the fact that
$AT$-free graphs are $5$-chordal and distance-hereditary graphs are
$4$-chordal.

\begin{te}\cite[Lemma 6]{GKKPP}   \cite[Theorem 3.3]{GKKPP1} If $G$ is a $k$-chordal graph, then $\mathbbm{t}\mathbbm{l}(G)\leq \lfloor \frac{k}{2}\rfloor.$
\label{thm9}
\end{te}

\begin{proof}[Outline] To obtain a minimal triangulation of $G,$ it suffices to select a maximal set of pairwise parallel
minimal separators of $G$ and add edges to make each of them a
clique \cite[Theorem 4.6]{Parra}. It is easy to check that each such
new edge connects two points of distance at most  $\lfloor
\frac{k}{2}\rfloor$ apart in $G.$
\end{proof}

The following is an interesting extension  of the classical result
that trees are $0$-hyperbolic and its proof can be given in a way
generalizing     the well-known proof of the latter fact.

\begin{te}
\cite[Proposition 13]{CDEHV}  A graph  $G$ is $k$-hyperbolic
provided its tree-length is no greater than  $k.$
          \label{thm10}
\end{te}

It is noteworthy that a converse of Theorem \ref{thm10} has also
been established, which means that hyperbolicity and tree-length are
comparable parameters of tree-likeness.

\begin{te}\cite[Proposition 14]{CDEHV}   The inequality $\mathbbm{t}\mathbbm{l}(G)\leq 12k+8k\log_2n + 17$ holds for any
 $k$-hyperbolic graph $G$    with $n$ vertices.
\end{te}

\subsubsection{Three examples}

Let us try our hand at three examples to get a feeling of the
concept of hyperbolicity.   The first example says that  graphs with
small diameter, hence those so-called small-world networks, must
have low hyperbolicity. Note that additionally similar simple
results will be reported as    Lemmas \ref{first} and \ref{lem14}.

\begin{example}  \cite[p. 683]{KM02} The hyperbolicity of a graph $G$  with diameter $D$ is at most
$\lfloor \frac{D}{2}\rfloor$. \label{diam}
\end{example}
\begin{proof} Take $x,y,u,v\in V(G)$. Our goal  is to show that   $\delta (x,y,u,v)\leq
\frac{D}{2}$. Without loss of generality, assume that
\begin{equation}  \label{Leipzig}  xy+uv\geq xu+yv\geq xv+yu
\end{equation} and hence \begin{equation}\delta
(x,y,u,v)=\frac{1}{2}((xy+uv)-(xu+yv)). \label{coffee}
\end{equation}
 In the first place, we have
$$
xu+yu\geq xy, ux+vx\geq uv,   xv+yv\geq xy, vy+uy\geq uv.
$$
Summing up these inequalities yields $(xu+yv)+(xv+yu)\geq  xy+uv$,
which, according to Eq. \eqref{Leipzig}, implies that
$$xu+yv\geq \frac{1}{2}(xy+uv).$$
This along with Eq. \eqref{coffee}  gives  $\delta (x,y,u,v)\leq
\frac{1}{4}(xy+uv)\leq \frac{D}{2}.$
 Moreover, if  $\delta (x,y,u,v)=
\frac{D}{2}$, then we have
 \begin{equation}xu+yv=D,
 xv+yv=xy=D, xu+xv=uv=D.
 \label{JS}
 \end{equation}
By adding the equalities in Eq. \eqref{JS} together, we see that
$3D=2(xu+xv+yv)$ and so $D$ must be even.
\end{proof}

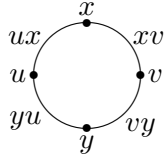
\begin{figure}
\hspace{-33mm}
\unitlength 1mm 
\linethickness{0.4pt}
\ifx\plotpoint\undefined\newsavebox{\plotpoint}\fi \hspace{20mm}
\begin{center}
\begin{picture}(10,10)
\put(10,10){\circle{20}}\put(10,17){\circle*{1}}\put(10,18){\makebox(0,0)[b]{$x$}}
\put(17,10){\circle*{1}}\put(18,10){\makebox(0,0)[l]{$v$}}
\put(10,3){\circle*{1}}\put(10,2){\makebox(0,0)[t]{$y$}}
\put(3,10){\circle*{1}}\put(2,10){\makebox(0,0)[r]{$u$}}
\put(16,15){\makebox(0,0)[l]{$xv$}}\put(15,3){\makebox(0,0)[l]{$vy$}}
\put(4,4){\makebox(0,0)[r]{$yu$}}\put(4,15){\makebox(0,0)[r]{$ux$}}
 \end{picture}
\end{center}
\caption{Four points in an $n$-cycle.}\label{n-cycle}
\end{figure}

The bound asserted by  Example \ref{diam} is clearly not tight when
$D=1.$ But, as can be seen from the next example,   the bound given
in  Example  \ref{diam} in terms of the diameter $D$ is best
possible for every $D\geq 2$. Note that this forthcoming example can
also be seen directly via Example \ref{diam}, as indicated in
\cite[p. 683]{KM02}.

\begin{example} \cite[p. 683]{KM02}   For any $n\geq 3,$ the  chordality  of the $n$-cycle  is $n$ while
the hyperbolicity of the $n$-cycle is \begin{equation}
\delta^*(C_n)=\begin{cases} \lfloor
\frac{n}{4}\rfloor-\frac{1}{2},&\text{if $n\equiv 1\ (\!\!\!\!\!\!\mod  4)$;}\\
\lfloor \frac{n}{4}\rfloor,&\text{else.}
\end{cases}
\label{USTC}
\end{equation}
Note that the diameter of  $C_n$ is  $\lfloor \frac{n}{2}\rfloor$
and
$$ \delta^*(C_n)=\begin{cases}
\frac{\lfloor \frac{n}{2}\rfloor}{2},&\text{if $n\equiv 0\ (\!\!\!\!\!\!\mod 4)$;}\\
\frac{\lfloor \frac{n}{2}\rfloor}{2}-\frac{1}{2},&\text{else.}
\end{cases}$$
\label{exam7}
\end{example}
\begin{proof} To prove Eq. \eqref{USTC}, we need to estimate $\delta
(x,y,u,v)$ for any four vertices $x,y,u,v$ of the $n$-cycle graph
$C_n.$  If there is a geodesic connecting two vertices and passing
through all the four vertices $x,y,u,v$, we surely have $\delta
(x,y,u,v) =0$ just because we know that the hyperbolicity of a path
is $0.$ So, we can assume that the cycle $C_n$ is $[c_0c_1\ldots
c_{n-1}]$ where $c_0=x,c_{xv}=v,c_{xv+vy}=y,c_{xv+vy+yu}=u$ and
\begin{equation}xv+vy+yu+ux=n;
\label{sum}
\end{equation} see Fig.     \ref{n-cycle}.  With no loss of generality, we assume that
\begin{equation}xu-vy\geq |xv-uy|.
\label{morning}
\end{equation} This implies  $ xu+uy \geq   xv+vy
$ and $ vx+xu\geq  vy+yu$. According to the geometric distribution
of the four points, we then come to
$$xy=xv+vy\ \ \text{and} \ \ vu=vy+yu.
$$ It follows that
\begin{equation}xy+vu=(xv+yu)+2vy\label{expo}\end{equation} and
$$xy+vu=xv+vy+vy+yu=(xv+vy+yu)+vy\geq xu+vy.$$  At the moment, we
see that there are only two possibilities, either $xy+vu\geq
xu+vy>xv+yu$ or $xy+vu\geq xv+yu \geq xu+vy.$

 If the first case
happens, we have
\begin{equation}
\begin{array}{cll}\delta(x,y,u,v)&=&\frac{1}{2}(xy+vu- xu-vy)
\\
&= & \frac{n}{2}-xu.\ \ \text{(By Eqs. \eqref{sum} and
\eqref{expo})}
 \end{array}
 \label{Kool}
 \end{equation}
By
 Eqs.  \eqref{sum}  and \eqref{morning} and  $xu+vy>xv+yu$, we see
 that $xu\geq
 \begin{cases} \lfloor \frac{n}{4}\rfloor +1,&\text{if $n\equiv 0,1,2\ (\!\!\!\!\mod  4)$,}\\
 \lfloor \frac{n}{4}\rfloor +2,&\text{if $n\equiv 3\ (\!\!\!\!\mod  4)$,}
\end{cases}
$
 and hence  Eq.  \eqref{Kool}  forces
 \begin{equation}\label{cycle1}\delta(x,y,u,v)= \frac{n}{2} - xu  \leq  \begin{cases}
 \lfloor \frac{n}{4}\rfloor -1,&\text{if $n\equiv 0\ (\!\!\!\!\!\!\mod 4)$;}\\
 \lfloor \frac{n}{4}\rfloor ,&\text{if $n\equiv 2\ (\!\!\!\!\!\!\mod 4)$;}\\
 \lfloor \frac{n}{4}\rfloor -\frac{1}{2},&\text{if $n\equiv 1,3\ (\!\!\!\!\!\!\mod 4)$.}
\end{cases}
\end{equation}

 For the second case, we have
$$
\begin{array}{cll}\delta(x,y,u,v)&=&\frac{1}{2}(xy+vu-
xv-yu)
\\
&= & vy,\ \ \text{(By Eq. \eqref{expo})}
\end{array}$$
and hence by
 Eqs.  \eqref{sum}  and \eqref{morning} and  $xv+yu\geq vy +xu$, we further  obtain
 \begin{equation}\label{cycle2}  \delta(x,y,u,v)=vy\leq \begin{cases}
 \lfloor \frac{n}{4}\rfloor ,&\text{if $n\equiv 0,2,3\ (\!\!\!\!\!\mod 4)$;}\\
 \lfloor \frac{n}{4}\rfloor -1,&\text{if $n\equiv 1\ (\!\!\!\!\!\mod 4)$.}
\end{cases}
\end{equation}

Combining Eqs. \eqref{cycle1} and \eqref{cycle2}  yields
\begin{equation}\label{quartet}
\delta (x,y,u,v)\leq \begin{cases} \lfloor
\frac{n}{4}\rfloor,&\text{if $n\equiv 0,2,3\ (\!\!\!\!\!\mod 4)$;}\\
\lfloor
\frac{n}{4}\rfloor-\frac{1}{2},&\text{if $n\equiv 1\ (\!\!\!\!\!\!\mod  4)$.}\\
\end{cases}
\end{equation}
Taking
\begin{equation*}
 (x,v,y,u)= \begin{cases} (c_0,c_k,c_{2k},c_{3k}),&\text{if $n\equiv 0\ (\!\!\!\!\!\!\mod  4)$,}\\
(c_0,c_k,c_{2k},c_{3k}),&\text{if $n\equiv 1\ (\!\!\!\!\!\!\mod  4)$,}\\
(c_0,c_k,c_{2k},c_{3k+1}),&\text{if $n\equiv 2\ (\!\!\!\!\!\!\mod  4)$,}\\
(c_0,c_{k+1},c_{2k+1},c_{3k+2}),&\text{if $n\equiv 3\
(\!\!\!\!\!\!\mod 4)$,}
\end{cases}
\end{equation*}
we see that Eq. \eqref{quartet} is tight and hence Eq.  \eqref{USTC}
is established.
\end{proof}

For any two graphs  $G_1$ and  $G_2,$ we define its {\em Cartesian
product } $G_1\Box G_2$  to be the graph satisfying $V(G_1\Box
G_2)=V(G_1)\times V(G_2)$ and $d_{G_1\Box
G_2}((u_1,u_2),(v_1,v_2))=d_{G_1}(u_1,v_1)+d_{G_2}(u_2,v_2)$
\cite[\S 1.4]{IK}.

\begin{example}  \label{exam3} Let $G_1$ and  $G_2$ be two graphs satisfying  $\delta ^*(G_1)=\delta^*
(G_2)=0.$ Then $\delta^*(G_1\Box G_2)=\min (D_1,D_2)$ where $D_1$
and $D_2$ are the diameters of  $G_1$ and  $G_2$, respectively.
\end{example}

\begin{proof}  For any $v\in V(G_1\Box G_2)$, we often use
the convention that $v=(v_1,v_2)$ for $v_1\in V(G_1)$ and $v_2\in
V(G_2)$.  For any $u,v\in V(G_1\Box G_2)$, we write  $uv$ for
$d_{G_1\Box G_2}(u,v)$, $(uv)_1$  for $d_{G_1}(u_1,v_1)$, $(uv)_2$
for $d_{G_2}(u_2,v_2)$  and we  use $\delta$ for $\delta_{G_1\Box
G_2}$.

Take  $a,b\in V(G_1)$ such that $d_{G_1}(a,b)=D_1$ and take $c,d\in
V(G_2)$ such that $d_{G_2}(c,d)=D_2$. Set
$x=(a,c),y=(a,d),u=(b,c),v=(b,d)$. It is straightforward that
$\delta (x,y,u,v)=\min (D_1, D_2)$.

To complete the proof, we pick any four vertices  $x,y,u,v$ of
$G_1\Box G_2$ and aim to show that \begin{equation}\delta
(x,y,u,v)\leq \min (D_1, D_2). \label{sunny}\end{equation}  Let
$A=xy+uv,$ $A_1=(xy)_1+(uv)_1$, $A_2=(xy)_2+(uv)_2$, $B=xu+yv,$
$B_1=(xu)_1+(yv)_1$, $B_2=(xu)_2+(yv)_2$,   $C=xv+yu,$
  $C_1=(xv)_1+(yu)_1$, $C_2=(xv)_2+(yu)_2.$ Because
  $\delta^*(G_1)=\delta^*(G_2)=0$, we can suppose $A_1=\max (B_1, C_1)$ and
  $A_2=\max (B_2,C_2)$.

  If it happens either  $(A_1,A_2)=(B_1,B_2)$ or
  $(A_1,A_2)=(C_1,C_2)$, we can immediately conclude that $\delta
  (x,y,u,v)=0.$ By symmetry between    $B$  and  $C$ and between $G_1$  and  $G_2,$   it thus remains to
  deduce Eq.  \eqref{sunny}  under the condition that
\begin{equation*}B\geq C, A_1=B_1>C_1, \ \text{and}\ A_2=C_2>B_2.
\label{KAIST}
\end{equation*}

Since $A_1=B_1,$  we have $\delta
(x,y,u,v)=\frac{A-B}{2}=\frac{A_2-B_2}{2}$. We proceed with a direct
computation and find \begin{equation}\label{eq14}\delta (x,y,u,v)
=\frac{((xy)_2-(xu)_2)+((uv)_2-(yv)_2 )}{2}\leq (yu)_2\leq D_2.
\end{equation}
Making use of $A_2=C_2$ and  $B\geq C,$  we can obtain  instead
\begin{equation}\label{eq15}
\begin{array}{cll}\delta(x,y,u,v)&\leq&\frac{A_2-B_2+B-C}{2}=\frac{A_2-C_2+B_1-C_1}{2}=\frac{B_1-C_1}{2}
\\
&= & \frac{((xu)_1-(xv)_1)+((yv)_1-(yu)_1)}{2}\leq (uv)_1\leq D_1
\end{array}
 \end{equation} Combining  Eqs. \eqref{eq14} and \eqref{eq15} we now get Eq. \eqref{sunny}, as desired.
\end{proof}

\begin{remark}\label{rem8}
For any $t$ natural numbers $m_1,\ldots , m_t$, the  {\em
$t$-dimensional
 grid graph} $G_{m_1,\ldots , m_t}$  is the graph with vertex set $\{
1,2,\ldots,m_1\}\times \cdots \times \{ 1,2,\ldots,m_t\}$ and
$(i_1,\ldots ,i_t)$ and $(j_1,\ldots, j_t)$ are adjacent in
$G_{m_1\ldots,m_t}$ if any only if $\sum_{p=1}^t(i_p-j_p)^2=1.$
Example \ref{exam3}  implies that $\delta^*(G_{m_1,m_2})=\min
(m_1,m_2)-1$ and hence $G_{m,m}$ provides another example that the
bound reported in Example \ref{diam} is tight. It might be
interesting to determine the hyperbolicity of  $t$-dimensional
 grid graphs for $t\geq 3.$
\end{remark}

\begin{remark}\label{grid}
 Dourisboure and Gavoille   show  that the tree-length of $G_{n,m}$   is $ \min (n,m)$ if $n\not= m$ or $n=m$ is even and is $n-1$ if
$n=m$ is odd \cite[Theorem 3]{DG}. Remark \ref{rem8}  tells us that
$\delta^*(G_{n,m})=\min(m,n)-1$.  This says that Theorem \ref{thm10}
is   tight.
\end{remark}


\section{Chordality vs. hyperbolicity}\label{MR}

\subsection{Main results}\label{theorem}

Firstly, we point out that
 a graph with low hyperbolicity may have large
chordality. Indeed, take any graph $G$ and form the new graph $G'$
by adding an additional vertex and connecting this new vertex with
every vertex of $G$. It is obvious that  $\delta^*(G')\leq 1$ while
$\mathbbm{l}\mathbbm{c}(G')= \mathbbm{l}\mathbbm{c}(G)$ if  $G$ is
not a tree. Moreover, it is equally easy to see that $G'$ is even
$\frac{1}{2}$-hyperbolic if $G$ does not have any induced $4$-cycle
\cite[p. 695]{KM02}. Surely, this example does not preclude the
possibility that for many important graph classes we can bound their
chordality in terms of their hyperbolicity.

One of our  main results says that   hyperbolicity can be bounded
from above in terms of chordality.

\begin{te} \label{main} For each $k\geq 4,$ all  $k$-chordal graphs are  $\frac{\lfloor
\frac{k}{2}\rfloor}{2}$-hyperbolic.
\end{te}

\begin{remark}
  A graph is {\em bridged} \cite{AF,LS}   if it
does not contain any finite isometric cycles of length at least
four, or equivalently, if  it   is cop-win and has no chordless
cycle of length $4$ or $5$. In contrast to    Theorem \ref{main}, it
is interesting to note that the hyperbolicity of bridged graphs can
be arbitrarily high \cite[p. 684]{KM02}.
\end{remark}

\begin{remark}  Bandelt and Chepoi \cite[\S 5.2]{BC08}  make the remark that
 ``a characterization of
 all   $1$-hyperbolic   graphs   by   forbidden   isometric   subgraphs
is   not
 in   sight,  in  as  much  as  isometric  cycles  of  lengths   up  to  $7$   may   occur,   thus
complicating   the   picture''.  Note that our Theorem \ref{main}
says that all $5$-chordal graphs are $1$-hyperbolic and hence the
appearance of those chordless $6$-cycles and chordless $7$-cycles
may be a real headache to deal with in pursuing a characterization
of all $1$-hyperbolic graphs.
\end{remark}

\begin{figure}
\hspace{-33mm}
\unitlength 1mm 
\linethickness{0.4pt}
\ifx\plotpoint\undefined\newsavebox{\plotpoint}\fi \hspace{30mm}
\begin{center}
\begin{picture}(150,30)

\put(60,25){\circle*{1}}\put(60,28){\makebox(0,0)[t]{$v_8$}}

\put(55,15){\circle*{1}}\put(51,15){\makebox(0,0)[l]{$v_7$}}

\put(50,5){\circle*{1}}\put(50,1){\makebox(0,0)[b]{$v_6$}}
 \put(60,5){\circle*{1}}\put(60,1){\makebox(0,0)[b]{$v_5$}}
\put(70,25){\circle*{1}}\put(70,29){\makebox(0,0)[t]{$v_1$}}
\put(80,25){\circle*{1}}\put(80,28){\makebox(0,0)[t]{$v_2$}}
\put(75,15){\circle*{1}}\put(79,15){\makebox(0,0)[r]{$v_3$}}

\put(70,5){\circle*{1}}\put(70,1){\makebox(0,0)[b]{$v_4$}}

\qbezier(60,25)(55,15)(50,5)\qbezier(80,25)(75,15)(70,5)
\qbezier(55,15)(58,9)(60,5)\qbezier(70,25)(73,19)(75,15)\qbezier(50,5)(60,5)(70,5)
\qbezier(60,25)(70,25)(80,25)
\end{picture}
\end{center}

\caption{The outerplanar graph $F_2$ has  chordality $6$,
hyperbolicity $\frac{3}{2}$, and  tree-length $2$.}\label{counter}
\end{figure}
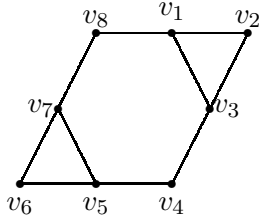

\begin{example}\label{exam8} For  any $t\geq 2$ we set
  $F_t$ to  be the  graph obtained from  the $4t$-cycle $[v_1v_2\cdots v_{4t}]$ by adding the two
 edges $\{v_1,v_3\}$ and $ \{v_{2t+1},v_{2t+3}\}$;  see  Fig. \ref{counter} for an illustartion of  $F_2$. Clearly,
 $\delta(v_2, v_{t+2},v_{2t+2}, v_{3t+2})=t-\frac{1}{2}$. Furthermore, we
 can check that
  $\mathbbm{l}\mathbbm{c}(F_t)=4t-2$ and $\delta^*(F_t)=t-\frac{1}{2}=\delta(v_2, v_{t+2},v_{2t+2},
  v_{3t+2})=\frac{\mathbbm{l}\mathbbm{c}(F_t)}{4}.$  $F_t$ is clearly an outerplanar graph. Thus, applying the result that   $\mathbbm{t}\mathbbm{l}(G)=\lceil \frac{\mathbbm{l}\mathbbm{c}(G)}{3}\rceil$   for
every outerplanar graph $G$   \cite[Theorem 1]{DG},  we even know
that  $\mathbbm{t}\mathbbm{l}(F_t)=\lceil \frac{4t-2}{3}\rceil$.
\end{example}

It is clear that if the bound claimed by  Theorem \ref{main} is
tight for $k=4t$ ($k=4t-2$) then it is tight for $k=4t+1$
($k=4t-1$). Consequently, Examples  \ref{exam7} and \ref{exam8}
indeed mean that the bound reported in Theorem \ref{main} is tight
for every $k\geq 4.$
      Surely,   the logical next   step would be to characterize     all those extremal
      graphs $G$
satisfying
\begin{equation}
\delta^*(G)=\frac{\lfloor
\frac{\mathbbm{l}\mathbbm{c}(G)}{2}\rfloor}{2}. \label{extremal}
\end{equation}
       However,   there seems to be still   a
long haul ahead in this direction.

\begin{remark} \label{rem13}
For any graph $G$ and any positive number $t$, we put $S^t(G)$ to be
a {\em subdivision graph} of $G$, which is obtained from $G$ by
replacing each edge $\{u,v\}$ of $G$ by a path $u,n_{u,v}^1,\ldots,
  n_{u,v}^{t-1}, v$ of length $t$ connecting $u$ and $v$ through a sequence of new vertices $n_{u,v}^1,\ldots,
   n_{u,v}^{t-1}$  (we surly require that  $n_{v,u}^q=n_{u,v}^{t-q}$).  For any four vertices $x,y,u,v\in V(G)$,
we obviously have $\delta_{S^t(G)}(x,y,u,v)=t\delta_G(x,y,u,v)$ and
so $\delta^*(S^t(G))\geq   t\delta^*(G)$.
 Instead of the trivial fact    $\mathbbm{l}\mathbbm{c} (S^t(G))\geq
   t \mathbbm{l}\mathbbm{c}(G)$, if the good shape of  $G$ permits us to deduce  a good upper bound
   of  $\mathbbm{l}\mathbbm{c} (S^t(G))$ in terms of
   $\mathbbm{l}\mathbbm{c}(G)$, we will see that $\delta^*(S^t(G))$
   is   high relative to  $\mathbbm{l}\mathbbm{c} (S^t(G))$
   provided so is $G.$  Recall that  the cycles whose lengths are divisible by $4$ as discussed  in Example
   \ref{exam7} are used to demonstrate  the tightness of the bound given in Theorem
   \ref{main}; also observe  that the graphs suggested by Example \ref{exam8} is nothing but a slight
   ``perturbation''
of cycles of length divisible by $4.$ Since $C_{4t}=S^t(C_{4})$,
these examples can be said to be generated by the ``seed''  $C_4.$
It might
   deserve to look for some other good ``seeds" from which we can use
   the above subdivision operation or its variant to produce graphs
   satisfying Eq. \eqref{extremal}.
\end{remark}

\begin{figure}
\hspace{-33mm}
\unitlength 1mm 
\linethickness{0.4pt}
\ifx\plotpoint\undefined\newsavebox{\plotpoint}\fi \hspace{30mm}
\begin{center}
\begin{picture}(150,100)

\put(30,60){\circle*{1}}\put(33,57){\makebox(0,0)[br]{$y$}}
\put(0,60){\circle*{1}}\put(-3,57){\makebox(0,0)[bl]{$u$}}
\put(30,90){\circle*{1}}\put(33,93){\makebox(0,0)[tr]{$v$}}
\put(0,90){\circle*{1}}\put(-3,93){\makebox(0,0)[tl]{$x$}}
\put(15,50){\makebox(0,0)[l]{$C_{4}$}}

\qbezier(30,90)(20,90)(0,90)\qbezier(0,90)(0,80)(0,60)\qbezier(0,60)(20,60)(30,60)\qbezier(30,60)(30,80)(30,90)

\put(70,60){\circle*{1}}\put(70,56){\makebox(0,0)[b]{$y$}}

\put(50,80){\circle*{1}}\put(47,80){\makebox(0,0)[l]{$u$}}
\put(90,80){\circle*{1}}\put(93,80){\makebox(0,0)[r]{$v$}}

\put(70,100){\circle*{1}}\put(70,103){\makebox(0,0)[t]{$x$}}

\put(60,90){\circle*{1}}\put(57,90){\makebox(0,0)[l]{$a$}}
\put(80,90){\circle*{1}}\put(83,90){\makebox(0,0)[r]{$b$}}
\put(60,70){\circle*{1}}\put(57,70){\makebox(0,0)[l]{$c$}}
\put(80,70){\circle*{1}}\put(83,70){\makebox(0,0)[r]{$d$}}
\put(70,50){\makebox(0,0)[l]{$H_{1}$}}

\qbezier(50,80)(60,70)(70,60)\qbezier(70,60)(80,70)(90,80)\qbezier(90,80)(80,90)(70,100)\qbezier(70,100)(60,90)(50,80)
\qbezier(60,90)(70,80)(80,70)\qbezier(60,90)(60,80)(60,70)\qbezier(60,70)(70,70)(80,70)\qbezier(80,70)(80,80)(80,90)
\qbezier(60,90)(70,90)(80,90)

\put(120,60){\circle*{1}}\put(120,56){\makebox(0,0)[b]{$y$}}

\put(100,80){\circle*{1}}\put(97,80){\makebox(0,0)[l]{$u$}}
\put(140,80){\circle*{1}}\put(143,80){\makebox(0,0)[r]{$v$}}

\put(120,100){\circle*{1}}\put(120,103){\makebox(0,0)[t]{$x$}}

\put(110,90){\circle*{1}}\put(107,90){\makebox(0,0)[l]{$a$}}
\put(130,90){\circle*{1}}\put(133,90){\makebox(0,0)[r]{$b$}}
\put(110,70){\circle*{1}}\put(107,70){\makebox(0,0)[l]{$c$}}
\put(130,70){\circle*{1}}\put(133,70){\makebox(0,0)[r]{$d$}}
\put(120,50){\makebox(0,0)[l]{$H_{2}$}}

\qbezier(100,80)(110,70)(120,60)\qbezier(120,60)(130,70)(140,80)\qbezier(140,80)(130,90)(120,100)\qbezier(120,100)(110,90)(100,80)
\qbezier(110,90)(120,80)(130,70)\qbezier(110,70)(120,80)(130,90)\qbezier(130,90)(130,80)(130,70)\qbezier(110,90)(120,90)(130,90)
\qbezier(110,90)(110,80)(110,70)\qbezier(110,70)(120,70)(130,70)

\put(20,0){\circle*{1}}\put(20,-4){\makebox(0,0)[b]{$y$}}

\put(0,20){\circle*{1}}\put(-3,20){\makebox(0,0)[l]{$u$}}
\put(40,20){\circle*{1}}\put(43,20){\makebox(0,0)[r]{$v$}}

\put(20,30){\circle*{1}}\put(20,33){\makebox(0,0)[t]{$x$}}

\put(10,10){\circle*{1}}\put(7,10){\makebox(0,0)[l]{$c$}}
\put(30,10){\circle*{1}}\put(33,10){\makebox(0,0)[r]{$d$}}
\put(22,-6){\makebox(0,0)[l]{$H_{3}$}}

\qbezier(0,20)(10,10)(20,0)\qbezier(20,0)(30,10)(40,20)\qbezier(40,20)(30,25)(20,30)\qbezier(20,30)(10,25)(0,20)
 \qbezier(10,10)(20,10)(30,10)

\put(70,0){\circle*{1}}\put(70,-4){\makebox(0,0)[b]{$y$}}

\put(50,20){\circle*{1}}\put(47,20){\makebox(0,0)[l]{$u$}}
\put(90,20){\circle*{1}}\put(93,20){\makebox(0,0)[r]{$v$}}

\put(70,40){\circle*{1}}\put(70,43){\makebox(0,0)[t]{$x$}}

\put(60,30){\circle*{1}}\put(57,30){\makebox(0,0)[l]{$a$}}
\put(80,30){\circle*{1}}\put(83,30){\makebox(0,0)[r]{$b$}}
\put(60,10){\circle*{1}}\put(57,10){\makebox(0,0)[l]{$c$}}
\put(80,10){\circle*{1}}\put(83,10){\makebox(0,0)[r]{$d$}}
\put(72,-6){\makebox(0,0)[l]{$H_{4}$}}

\qbezier(50,20)(60,10)(70,0)\qbezier(70,0)(80,10)(90,20)\qbezier(90,20)(80,30)(70,40)\qbezier(70,40)(60,30)(50,20)
\qbezier(60,30)(70,20)(80,10) \qbezier(80,30)(70,20)(60,10)

\put(120,0){\circle*{1}}\put(120,-4){\makebox(0,0)[b]{$y$}}
\put(110,10){\circle*{1}}\put(107,10){\makebox(0,0)[l]{$c$}}
\put(110,30){\circle*{1}}\put(107,30){\makebox(0,0)[l]{$a$}}
\put(130,30){\circle*{1}}\put(133,30){\makebox(0,0)[r]{$b$}}
\put(130,10){\circle*{1}}\put(133,10){\makebox(0,0)[r]{$d$}}
\put(100,20){\circle*{1}}\put(96,20){\makebox(0,0)[l]{$u$}}
\put(140,20){\circle*{1}}\put(143,20){\makebox(0,0)[r]{$v$}}
\put(120,40){\circle*{1}}\put(120,43){\makebox(0,0)[t]{$x$}}
\qbezier(100,20)(110,10)(120,0)\qbezier(140,20)(130,10)(120,0)\qbezier(110,30)(120,20)(130,10)
\qbezier(120,40)(110,30)(100,20)\qbezier(120,40)(130,30)(140,20)
\put(122,-6){$H_5$}
 \end{picture}

\end{center}
\caption{Six  $5$-chordal graphs with hyperbolicity
$1$.}\label{fig0}
\end{figure}
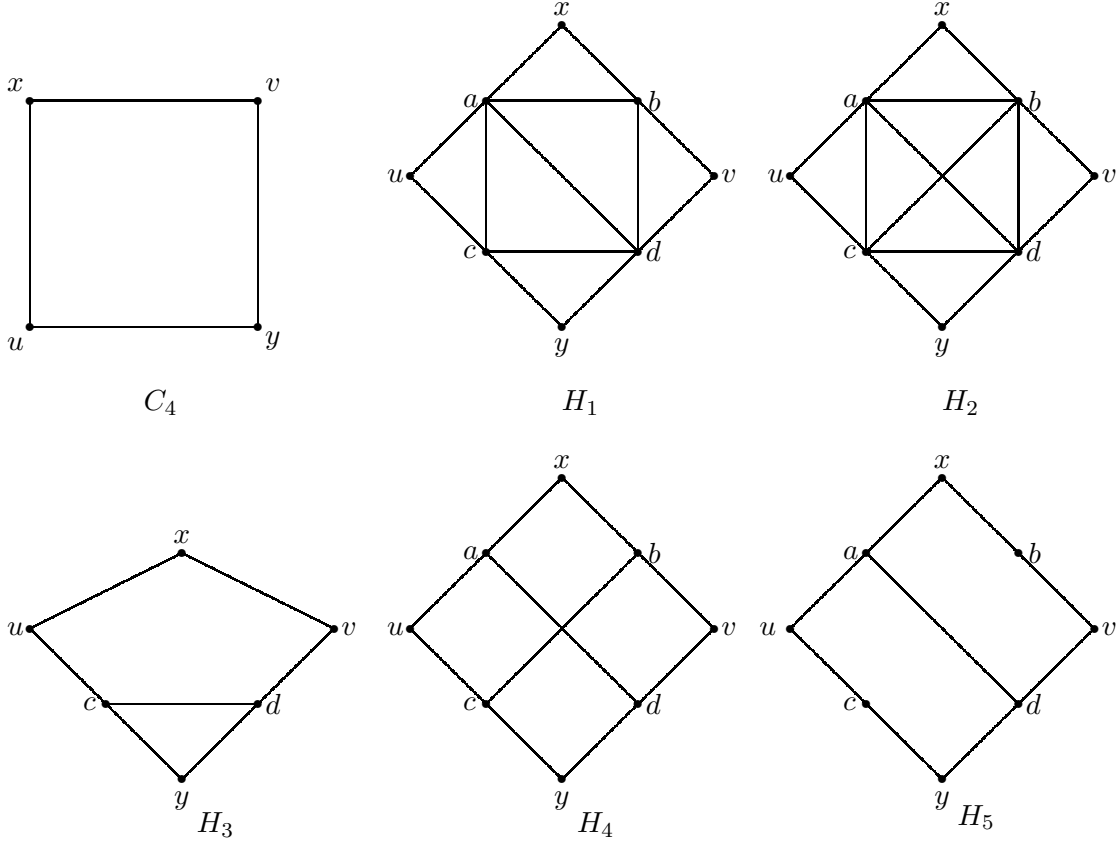

Let  $C_4$,  $H_1$,  $H_2$, $H_3$, $H_4$  and  $H_5$ be the  graphs
displayed in Fig. \ref{fig0}.  It is simple to check that each of
them has hyperbolicity $1$ and is $5$-chordal.   Besides Theorem
\ref{main}, another main contribution of this paper is the
following, which says that  $5$-chordal graphs will be
$\frac{1}{2}$-hyperbolic as soon as these six obvious obstructions
do not occur.


\begin{te}\label{main1}
 A $5$-chordal graph has hyperbolicity one  if and only if  one
of $C_4,H_1,H_2,H_3,H_4, H_5$ appears  as an isometric subgraph of
it.
\end{te}

Returning to Remark \ref{rem13}, it is natural to investigate if
some graphs mentioned in Theorem \ref{main1} besides  $C_4$ can be
used as ``good seeds''. The next example comes from Gavoille
\cite{CGa}.

\begin{figure}
\hspace{-33mm}
\unitlength 1mm 
\linethickness{0.6pt}
\ifx\plotpoint\undefined\newsavebox{\plotpoint}\fi \hspace{30mm}
\begin{center}
\begin{picture}(50,40)

\put(20,0){\circle*{1}}\put(20,-4){\makebox(0,0)[b]{$y$}}
\put(0,20){\circle*{1}}\put(-1,20){\makebox(0,0)[r]{$u$}}
\put(40,20){\circle*{1}}\put(41,20){\makebox(0,0)[l]{$v$}}
\put(20,40){\circle*{1}}\put(20,44){\makebox(0,0)[t]{$x$}}
\put(10,30){\circle*{1}}\put(7,30){\makebox(0,0)[l]{$a$}}
\put(10,10){\circle*{1}}\put(9,10){\makebox(0,0)[r]{$c$}}
\put(30,10){\circle*{1}}\put(31,10){\makebox(0,0)[l]{$d$}}
\put(30,30){\circle*{1}}\put(31,30){\makebox(0,0)[l]{$b$}}
\put(5,25){\circle*{1}}\put(4,25){\makebox(0,0)[r]{$u_a$}}
\put(5,15){\circle*{1}}\put(4,15){\makebox(0,0)[r]{$u_c$}}
\put(15,5){\circle*{1}}\put(14,5){\makebox(0,0)[r]{$y_c$}}
\put(25,5){\circle*{1}}\put(26,5){\makebox(0,0)[l]{$y_d$}}
\put(35,15){\circle*{1}}\put(36,15){\makebox(0,0)[l]{$v_d$}}
\put(35,25){\circle*{1}}\put(36,25){\makebox(0,0)[l]{$v_b$}}
\put(15,35){\circle*{1}}\put(14,35){\makebox(0,0)[r]{$x_a$}}
\put(25,35){\circle*{1}}\put(26,35){\makebox(0,0)[l]{$x_b$}}
\qbezier[30](0,20)(10,30)(20,40)\qbezier[30](20,40)(30,30)(40,20)\qbezier[30](20,0)(30,10)(40,20)\qbezier[30](20,0)(10,10)(0,20)
\qbezier[20](10,30)(10,30)(10,10)
\qbezier[20](10,30)(20,30)(30,30)\qbezier[20](30,10)(30,10)(30,30)
\qbezier[20](10,10)(30,10)(30,10)
\qbezier(5,15)(5,20)(5,25)\qbezier(15,5)(20,5)(25,5)\qbezier(35,15)(35,20)(35,25)\qbezier(15,35)(20,35)(25,35)
\qbezier[20](10,30)(20,20)(30,10)\qbezier[20](10,10)(20,20)(30,30)
 \end{picture}
\end{center}
\caption{$\mathbb{G}_{4t}^q$.}\label{example13}
\end{figure}

\begin{example}\label{ep} \cite{CGa} Let $t,q$ be two positive integer with
$q<t$ and let $H_2$ be the graph shown in the  upper-right corner of
Fig. \ref{fig0}. We construct a planar graph $\mathbb{G}_{4t}^q$
from $S^t(H_2)$ as follows: let $u_a=n_{u,a}^{q}$,
$u_c=n_{c,u}^{q-1}$, $y_c=n_{y,c}^{q}$, $y_d=n_{d,y}^{q-1}$,
$v_d=n_{v,d}^{q}$, $v_b=n_{b,v}^{q-1}$, $x_b=n_{x,b}^{q}$,
$x_a=n_{a,x}^{q-1}$, and then add the new edges
  $\{u_a, u_c\}$, $\{x_a, x_b\}$, $\{y_c,
 y_d\}$, $\{v_b, v_d\}$ to  $S^t(H_2)$; see Fig. \ref{example13}. It
can be checked  that  $C=[u_a \cdots a \cdots x_ax_b\cdots b \cdots
v_bv_d\cdots d\cdots y_dy_c\cdots c \cdots u_c]$ is an isometric
  $4t$-cycle of $\mathbb{G}_{4t}^q$ and that $\mathbbm{l}\mathbbm{c}(\mathbb{G}_{4t}^q) =4t.$
It is also easy to see that $\delta_{\mathbb{G}_{4t}^q}(u,y,v,x)=t $
and thus  Theorem \ref{main} tells us that
$\delta^*(\mathbb{G}_{4t}^q)=t$.
\end{example}

Motivated by the above construction of Gavoille, we construct the
next graph family whose chordality parameters are $1$ modulo $4$.

\begin{example}
By deleting the edge $\{ y_c, n_{d,y}^{q-1}\}$ and adding a new edge
$\{ y_c,  n_{d,y}^{q}\}$, we obtain   from $\mathbb{G}_{4t}^q$ a
graph  $\mathbb{G}_{4t+1}^q$. Using similar analysis like Example
\ref{ep}, we find that   $\mathbbm{l}\mathbbm{c}(\mathbb{G}_{4t}^q)
=4t+1$ and $\delta^*(\mathbb{G}_{4t+1}^q)=t=\frac{\lfloor
\frac{4t+1}{2}\rfloor }{2}$.
\end{example}

Similar constructions using $F_2$ (see Fig. \ref{counter})    as the
``seed" will lead to corresponding extremal graphs whose chordality
parameters are $2$ or $3$ modulo $4$.

\begin{figure}
\hspace{-33mm}
\unitlength 1mm 
\linethickness{0.4pt}
\ifx\plotpoint\undefined\newsavebox{\plotpoint}\fi \hspace{30mm}
\begin{center}
\begin{picture}(150,30)

\put(60,25){\circle*{1}}\put(60,28){\makebox(0,0)[t]{$v_8$}}

\put(55,15){\circle*{1}}\put(51,15){\makebox(0,0)[l]{$v_7$}}

\put(53,11){\circle*{1}}\put(52,11){\makebox(0,0)[r]{$v_{67}$}}

\put(55,5){\circle*{1}}\put(55,4){\makebox(0,0)[t]{$v_{65}$}}

\put(50,5){\circle*{1}}\put(50,1){\makebox(0,0)[b]{$v_6$}}
 \put(60,5){\circle*{1}}\put(61,1){\makebox(0,0)[b]{$v_5$}}
\put(70,25){\circle*{1}}\put(69,26){\makebox(0,0)[b]{$v_1$}}
\put(78,21){\circle*{1}}\put(79,21){\makebox(0,0)[l]{$v_{23}$}}
\put(75,25){\circle*{1}}\put(75,26){\makebox(0,0)[b]{$v_{21}$}}
\put(80,25){\circle*{1}}\put(81,28){\makebox(0,0)[t]{$v_2$}}
\put(75,15){\circle*{1}}\put(79,15){\makebox(0,0)[r]{$v_3$}}

\put(70,5){\circle*{1}}\put(70,1){\makebox(0,0)[b]{$v_4$}}

\qbezier[20](60,25)(55,15)(50,5)\qbezier[20](80,25)(75,15)(70,5)
\qbezier[20](55,15)(58,9)(60,5)\qbezier[20](70,25)(73,19)(75,15)\qbezier[20](50,5)(60,5)(70,5)
\qbezier[20](60,25)(70,25)(80,25)\qbezier(53,11)(53,11)(55,5)\qbezier(78,21)(78,21)(75,25)
\end{picture}
\end{center}

\caption{$\mathbb{G}_{6(2t+1)}^q$.}\label{expo}
\end{figure}

\begin{example} \label{nuli} Let $t>q$ be two positive integers. We
construct an outerplanar graph $\mathbb{G}_{6(2t+1)}^q$ by adding
two new edges   $\{v_{21}, v_{23}\}$  and  $\{v_{65},
 v_{67}\}$
 to the graph $S^{2t+1}(F_2)$
 where
  $v_{21}=n_{v_2,v_1}^{q}$, $v_{23}=n_{v_3,v_2}^{q-1}$,
$v_{65}=n_{v_6,v_5}^{q}$, $v_{67}=n_{v_7,v_6}^{q-1}$; see Fig.
 \ref{expo} for an illustration. It is not hard to check that $\mathbbm{l}\mathbbm{c}(\mathbb{G}_{6(2t+1)}^q)
=6(2t+1)$ and  $\delta^*(\mathbb{G}_{6(2t+1)}^q)=3t+\frac{3}{2}$.
Moreover, if we replace the edge $\{v_{21}, v_{23}\}$
 by the edge $\{v_{21}, n_{v_3,v_2}^{q}\}      $, then we obtain from $\mathbb{G}_{6(2t+1)}^q$
 another outerplanar graph $\mathbb{G}_{6(2t+1)+1}^q$ for which we
 have
$\mathbbm{l}\mathbbm{c}(\mathbb{G}_{6(2t+1)+1}^q) =6(2t+1)+1$ and
$\delta^*(\mathbb{G}_{6(2t+1)+1}^q)=3t+\frac{3}{2}$.
\end{example}

\begin{figure}
\hspace{-33mm}
\unitlength 1mm 
\linethickness{0.4pt}
\ifx\plotpoint\undefined\newsavebox{\plotpoint}\fi \hspace{30mm}
\begin{center}
\begin{picture}(120,60)

\put(0,25){\circle*{1}}

\put(10,15){\circle*{1}} \put(20,25){\circle*{1}}

\put(20,5){\circle*{1}}

\put(30,5){\circle*{1}} \put(30,25){\circle*{1}}
\put(50,25){\circle*{1}}\put(40,15){\circle*{1}}

\put(25,15){\circle*{1}}

\qbezier(0,25)(10,15)(20,5)\qbezier(20,5)(25,5)(30,5)\qbezier(30,5)(40,15)(50,25)\qbezier(50,25)(20,25)(0,25)
\qbezier(20,5)(25,15)(30,25)
\qbezier(20,25)(25,15)(30,5)\qbezier(10,15)(25,15)(40,15)\qbezier(20,25)(15,20)(10,15)\qbezier(40,15)(35,20)(30,25)

\put(25,-5){\makebox(0,0){$G_1$}}

\put(60,25){\circle*{1}}

\put(55,15){\circle*{1}} \put(65,15){\circle*{1}}

\put(50,5){\circle*{1}}

\put(60,5){\circle*{1}} \put(70,25){\circle*{1}}
\put(80,25){\circle*{1}}\put(75,15){\circle*{1}}

\put(70,5){\circle*{1}}

\qbezier(60,25)(55,15)(50,5)\qbezier(60,5)(65,15)(70,25)\qbezier(80,25)(75,15)(70,5)\qbezier(55,15)(65,15)(75,15)
\qbezier(60,25)(65,15)(70,5)
\qbezier(70,25)(65,15)(60,5)\qbezier(55,15)(58,9)(60,5)\qbezier(70,25)(73,19)(75,15)\qbezier(50,5)(60,5)(70,5)
\qbezier(60,25)(70,25)(80,25)

\put(60,-5){\makebox(0,0){$G_2$}}

\put(55,45){\circle*{1}}\put(75,55){\circle*{1}}\put(75,35){\circle*{1}}\put(85,45){\circle*{1}}

\put(65,55){\circle*{1}}\put(70,45){\circle*{1}}

\put(65,35){\circle*{1}}
\qbezier(55,45)(60,50)(65,55)\qbezier(55,45)(60,40)(65,35)\qbezier(65,35)(70,35)(75,35)
\qbezier(75,35)(80,40)(85,45)\qbezier(75,55)(80,50)(85,45)\qbezier(65,55)(70,55)(75,55)
\qbezier(65,55)(70,45)(75,35)

\put(70,30){\makebox(0,0){$G_3$}}

\put(5,45){\circle*{1}}\put(15,35){\circle*{1}}\put(25,35){\circle*{1}}
\put(35,45){\circle*{1}}\put(15,55){\circle*{1}}\put(25,55){\circle*{1}}
\qbezier(15,35)(10,40)(5,45)\qbezier(15,35)(20,35)(25,35)\qbezier(25,35)(30,40)(35,45)
\qbezier(25,55)(30,50)(35,45)\qbezier(15,55)(20,55)(25,55)\qbezier(5,45)(10,50)(15,55)
\put(20,30){\makebox(0,0){$C_6$}}
 \end{picture}
\end{center}

\caption{Four  graphs with hyperbolicity $1$ and chordality
$6$.}\label{figconjecture}
\end{figure}
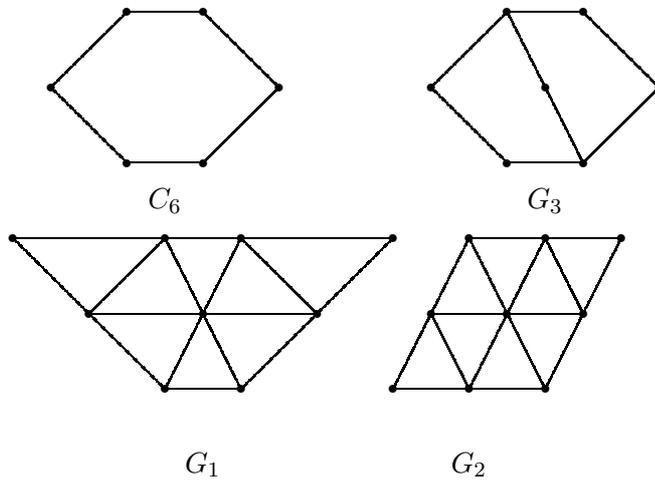

Let  $C_6,G_1,G_2,G_3$ be the graphs depicted in Fig.
 \ref{figconjecture}. It is clear that
 $G_1,G_2,G_3,C_4,C_6,H_i,i=1,\ldots,5,$ are $6$-chordal graphs with
 hyperbolicity $1$.

\begin{conjecture} \label{conj14}  A  $6$-chordal graph   is
$\frac{1}{2}$-hyperbolic if and only if it does not contain any of a
list of ten special graphs $G_1,G_2,G_3,C_4,C_6,H_i,i=1,\ldots,5,$
as an isometric subgraph.
\end{conjecture}

Let  $E_1$ and $E_2$ be the graphs  depicted  in Fig. \ref{fig
bridged}.     In comparison with Conjecture  \ref{conj14}, when we
remove the $6$-chordal restriction, we can present the following
characterization   of all $\frac{1}{2}$-hyperbolic graphs obtained
by Bandelt and Chepoi \cite{BC}.  We refer to \cite[Fact 1]{BC} for
two other characterizations; also see \cite{FJ,SC}.

 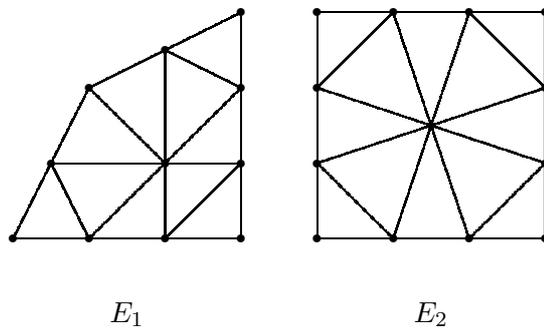
\begin{figure}
\hspace{-33mm}
\unitlength 1mm 
\linethickness{0.4pt}
\ifx\plotpoint\undefined\newsavebox{\plotpoint}\fi \hspace{30mm}
\begin{center}
\begin{picture}(80,40)
\put(10,25){\circle*{1}}

\put(5,15){\circle*{1}} \put(20,30){\circle*{1}}

\put(0,5){\circle*{1}}

\put(10,5){\circle*{1}} \put(20,15){\circle*{1}}
\put(30,25){\circle*{1}}
\put(20,5){\circle*{1}}\put(30,5){\circle*{1}}
\put(30,35){\circle*{1}}\put(30,15){\circle*{1}}

\qbezier(10,25)(5,15)(0,5)\qbezier(10,5)(20,15)(30,25)\qbezier(30,35)(30,15)(30,5)\qbezier(5,15)(20,15)(30,15)
\qbezier(10,25)(15,20)(20,15)
\qbezier(20,30)(20,15)(20,5)\qbezier(5,15)(8,9)(10,5)\qbezier(10,25)(20,30)(30,35)\qbezier(0,5)(10,5)(30,5)
\qbezier(20,5)(25,10)(30,15)\qbezier(20,30)(26,27)(30,25)

\put(15,-5){\makebox(0,0){$E_1$}}

\put(40,5){\circle*{1}}\put(50,5){\circle*{1}}\put(60,5){\circle*{1}}\put(70,5){\circle*{1}}

\put(40,15){\circle*{1}}\put(70,15){\circle*{1}}\put(40,25){\circle*{1}}

\put(70,25){\circle*{1}}\put(40,35){\circle*{1}}\put(50,35){\circle*{1}}\put(60,35){\circle*{1}}
\put(70,35){\circle*{1}}\put(55,20){\circle*{1}}

\qbezier(40,5)(60,5)(70,5)\qbezier(40,5)(40,10)(40,35)\qbezier(40,35)(50,35)(70,35)
\qbezier(70,35)(70,15)(70,5)\qbezier(40,15)(55,20)(70,25)\qbezier(40,25)(55,20)(70,15)
\qbezier(60,35)(55,20)(50,5)\qbezier(50,35)(55,20)(60,5)\qbezier(50,35)(45,30)(40,25)
\qbezier(40,15)(45,10)(50,5)\qbezier(60,35)(65,30)(70,25)

\qbezier(60,5)(65,10)(70,15)

\put(55,-5){\makebox(0,0){$E_2$}}
 \end{picture}
\end{center}

\caption{Two bridged graphs with hyperbolicity $1$.}\label{fig
bridged}
\end{figure}

\begin{te}  \cite[p. 325]{BC}  A graph  $G$ is $\frac{1}{2}$-hyperbolic if and only if
$G$ does not contain isometric $n$-cycles for any $n>5$,  for any
two vertices $x$ and $y$ of  $G$ one cannot find two non-adjacent
neighbors of  $x$  which are both closer to $y$ in  $G$ than $x$,
and none of the six graphs $H_1,H_2,G_1,G_2,E_1,E_2$ occurs as an
isometric subgraph of $G.$ \label{thmBC}
\end{te}

\begin{remark} Instead of Theorem \ref{thmBC}, it would be interesting to determine, if possible, a finite
 list of graphs such  that that a graph is $\frac{1}{2}$-hyperbolic if and only if
it does not include any graph from that list as an isometric
subgraph.  Koolen and Moulton point out a possible approach to
deduce such kind of a characterization in \cite[p. 696]{KM02}.
\end{remark}

Note that a $5$-chordal graph cannot contain any  isometric
$n$-cycle for $n>5$.  It is also easy to see that
$\mathbbm{l}\mathbbm{c}(G_1)=\mathbbm{l}\mathbbm{c}(G_2)=6,\mathbbm{l}\mathbbm{c}(E_1)=7,
\mathbbm{l}\mathbbm{c}(E_2)=8$. Therefore, we obtain the following
easy consequence of Theorem \ref{thmBC}. It is interesting to
compare it with Theorems \ref{main} and \ref{main1}.

\begin{corollary}   A  $5$-chordal graph $G$ is $\frac{1}{2}$-hyperbolic if and only if it does not
contain the graph   $H_1$ and  $H_2$ as isometric subgraphs and for
any two vertices $x$ and $y$ of  $G$ the neighbors of  $x$ which are
closer to $y$ than  $x$ are pairwise adjacent.
\end{corollary}

\subsection{Some consequences}\label{oct}

Note that  $\mathbbm{l}\mathbbm{c}(C_4)=4,
\mathbbm{l}\mathbbm{c}(H_1)=\mathbbm{l}\mathbbm{c}(H_2)=3,
\mathbbm{l}\mathbbm{c}(H_3)=\mathbbm{l}\mathbbm{c}(H_4)=\mathbbm{l}\mathbbm{c}(H_5)=5.$
The next two results follow immediately from Theorem \ref{main1}.

\begin{corollary}   Every  $4$-chordal graph must be $1$-hyperbolic and it has
hyperbolicity one  if and only if it contains one of $C_4$, $H_1$
and $H_2$ as an isometric subgraph.  \label{cor7}
\end{corollary}

\begin{corollary}   \cite[Theorem 1.1]{BKM01}  Every  chordal
graph is $1$-hyperbolic and it has hyperbolicity one if and only if
it contains either $H_1$ or $H_2$ as an isometric
subgraph.\label{BKM}
\end{corollary}

We remark that as long as every $4$-chordal graph is 1-hyperbolic is
known, Corollary \ref{cor7} also immediately follows from Corollary
\ref{BKM}.  In addition, it is noteworthy  that the first part of
Corollary \ref{BKM}, namely every chordal graph is $1$-hyperbolic is
immediate from Theorem \ref{thm10} as chordal graphs have
tree-length $1$.

\begin{corollary}  Each  weakly chordal graph is $1$-hyperbolic and has
hyperbolicty one if and only if it contains one of $C_4,H_1,H_2$ as
an isometric subgraph.
\end{corollary}

\begin{proof} By definition, each weakly chordal graph is $4$-chordal. It is also easy to check that   that $C_4,H_1$ and $H_2$
are all weakly chordal.  Hence, the result follows from  Corollary
\ref{cor7}.
\end{proof}

\begin{corollary} All strongly chordal graphs are
$\frac{1}{2}$-hyperbolic.
\end{corollary}
\begin{proof} Note that the cycle $C=[x,a,u,c,y,d,v,b]$ in $H_1$
and $H_2$  does not have any odd chord and hence neither  $H_1$ nor
$H_2$ can appear as an isometric subgraph of a strongly chordal
graph.  Since strongly chordal graphs must be chordal graphs, this
result holds by Corollary \ref{BKM}.
\end{proof}

\begin{corollary} All threshold graphs are  $\frac{1}{2}$-hyperbolic.
\end{corollary}
\begin{proof} It is obvious that    threshold graphs are chordal as they   contain neither   $4$-cycle
nor path of length $3$ as induced subgraph.  Since the subgraph
induced by  $x,u,b,c$ in either $H_1$  or  $H_2$ is just the
complement of $C_4$,   the result follows from Corollary \ref{BKM}
and the definition of a threshold graph.
\end{proof}

\begin{corollary}  Every  $AT$-free graph  is  $1$-hyperbolic and it has hyperbolicity one if and only if it contains $C_4$
as an isometric subgraph. \label{AT-free}
\end{corollary}
\begin{proof}   First observe that an
 $AT$-free graph  must be $5$-chordal. Further notice that
the triple $u,y,v$ is an $AT$ in any of the graphs $H_1,\ldots,H_5.$
Now, an application of Theorem  \ref{main1}  concludes the proof.
\end{proof}

\begin{corollary} \label{com}   A cocomparability graph  is  $1$-hyperbolic and has hyperbolicity one if and only if it contains $C_4$
as an isometric subgraph.
\end{corollary}
\begin{proof} We know that cocomparability graphs are $AT$-free and
 $C_4$ is a cocomparability graph. Thus the result
comes directly  from Corollary  \ref{AT-free}.   The deduction of
this result can also be made via Corollary \ref{cor7} and the fact
that cocomparability graphs are $4$-chordal \cite{BT1, Gallai}.
\end{proof}

\begin{corollary}  A permutation graph is $1$-hyperbolic and has hyperbolicity one if and only
if it contains $C_4$ as an isometric subgraph.
\end{corollary}

\begin{proof} Every permutation graph is a cocomparability graph and   $C_4$ is a permutation
graph. So,  the result follows from Corollary \ref{com}.
\end{proof}

\begin{corollary} \cite[p. 16]{BC08}
A distance-hereditary graph
 is always   $1$-hyperbolic and is    $\frac{1}{2}$-hyperbolic exactly  when  it  is
 chordal, or equivalently, when it contains no induced $4$-cycle.
 \label{distance}
 \end{corollary}

 \begin{proof} It is easy to see that distance-hereditary graphs
 must be $4$-chordal and can contain  neither  $H_1$ nor $H_2$ as an
isometric subgraph. The result now follows from Corollary
\ref{cor7}.
\end{proof}

\begin{corollary} A cograph is $1$-hyperbolic and has hyperbolicity one if  and only if it contains
$C_4$ as an isometric subgraph.
\end{corollary}

\begin{proof} We know that $C_4$ is a cograph and every  cograph is   ditance-hereditary. Applying  Corollary
\ref{distance}  yields  the  required result.
\end{proof}

\section{Relevant tree-likeness parameters}\label{parameter}

\subsection{Tree-length}

It turns out that tree-length is a very useful concept for
connecting chordality and hyperbolicity.
 Indeed, the following theorem, which can be read from  Theorem \ref{main} (Corollary \ref{BKM}), comes directly from
 Theorems
\ref{thm9} and   \ref{thm10}. This result is firstly notified to us
by Dragan \cite{Dragan} and is presumably in the folklore.

\begin{te} For any $k\geq 3,$ every $k$-chordal graph is  $\lfloor \frac{k}{2}\rfloor $-hyperbolic.\label{June}
\end{te}

In view of Remark \ref{grid}, to get better estimate than  Theorem
\ref{June} along the same approach
  one may try to beef
up   Theorem \ref{thm9}. We point out that Dourisboure  and Gavoille
\cite[Question 1]{DG} posed as an open problem that whether or not
\begin{equation}  \label{eq22} \mathbbm{t}\mathbbm{l}(G)\leq \lceil \frac{\mathbbm{l}\mathbbm{c}(G)}{3}\rceil
\end{equation} is true.
The {\em $k$th-power} of a graph $G$, denoted  $G^k,$ is the graph
with $V(G)$ as vertex set and there is an edge connecting two
vertices $u$ and $v$ if and only if $d_G(u,v)\leq k.$  Let us
interpret the problem of Dourisboure  and Gavoille as a Chordal
Graph Sandwich Problem:

\begin{question}
For any  graph  $G$,   is there always a chordal  graph $H$  such
that  $V(H)=V(G)=V(G^{\lceil
\frac{\mathbbm{l}\mathbbm{c}(G)}{3}\rceil})$ and $E(G)\subseteq E(H)
\subseteq G^{\lceil \frac{\mathbbm{l}\mathbbm{c}(G)}{3}\rceil}$?
\end{question}

If  \eqref{eq22} can be established, it will be the best we can
expect in the sense that $\mathbbm{t}\mathbbm{l}(G)=\lceil
\frac{\mathbbm{l}\mathbbm{c}(G)}{3}\rceil$   for every outerplanar
graph $G$   \cite[Theorem 1]{DG}.



\subsection{Approximating trees, slimness and thinness}

We introduce in this subsection two general approaches to connect
chordality with hyperbolicity. A result is given together with a
proof only when that proof is short and  when we do not find it
appear very explicitly elsewhere. This section also aims to provide
the reader a warm-up before entering the longer proof in the main
part of this paper.

A
 result   weaker than Theorem \ref{main} (Theorem \ref{June}) and reported in
\cite[p. 64]{CE} as well as \cite[p. 3]{CDEHVX} is that  each
$k$-chordal graph is $k$-hyperbolic. The two approaches to be
reported below by far basically only lead to
 this weaker result.
 Despite of this, it might be interesting  to see  different ways of bounding
  hyperbolicity  in terms of chordality via the use of some other intermediate tree-likeness parameters.

The first approach is to look at distance approximating trees. A
tree $T$ is a {\em distance $t$-approximating tree} of a graph $G$
provided $V(T)=V(G)$ and $|d_G(u,v)-d_T(u,v)|\leq t$ for any $u,v\in
V(G)$ \cite{Balint,BCD,CD,DY}. It is well-known that a graph with a
good distance approximating tree    will have low hyperbolicity,
which is briefly mentioned in  \cite[p. 3]{CDEHVX}   and  \cite[p.
64]{CE}   and is in the same spirit of a general result on
hyperbolic geodesic metric spaces \cite[p. 402, Theorem 1.9]{BH}. We
make this point clear in the following simple lemma.

\begin{lemma}\label{lem1}  Let $G$   be a graph and $t$ be a nonnegative integer.    If $G$  has a distance $t$-approximating tree
$T$, then $G$ is $2t$-hyperbolic.
\end{lemma}
\begin{proof}
For any $x,y,u,v\in V(G)$, our aim is to show that
$\delta_G(x,y,u,v)\leq 2t.$    Assume, as we may, that
$d_G(x,y)+d_G(u,v)\geq d_G(x,u)+d_G(y,v)\geq d_G(x,v)+d_G(y,u).$
Since the tree metric  $d_T$ is a four-point  inequality   metric
(or additive metric)   \cite{DD}, we know that $\delta^*(T)=0$ and
so the following three cases are exhaustive.

\paragraph {\sc Case 1:}
$d_T(x,y)+d_T(u,v)= d_T(x,u)+d_T(y,v)\geq d_T(x,v)+d_T(y,u).$

$\delta_G(x,y,u,v)=\frac{1}{2}(d_G(x,y)+d_G(u,v))
-\frac{1}{2}(d_G(x,u)+d_G(y,v))\leq \frac{1}{2}(d_T(x,y)+d_T(u,v)+
2t ) -\frac{1}{2}(d_T(x,u)+d_T(y,v)-   2t)=2t. $

\paragraph {\sc Case 2:}
$d_T(x,y)+d_T(u,v)= d_T(x,v)+d_T(y,u)\geq d_T(x,u)+d_T(y,v) $

$\delta_G(x,y,u,v)=\frac{1}{2}(d_G(x,y)+d_G(u,v))
-\frac{1}{2}(d_G(x,u)+d_G(y,v))\leq \frac{1}{2}(d_G(x,y)+d_G(u,v))
-\frac{1}{2}(d_G(x,v)+d_G(y,u))\leq  \frac{1}{2} (d_T(x,y)+d_T(u,v)+
2t ) -\frac{1}{2}(d_T(x,v)+d_T(y,u)-   2t) = 2t. $

\paragraph {\sc Case 3:}  $ d_T(x,v)+d_T(y,u)= d_T(x,u)+d_T(y,v) \geq   d_T(x,y)+d_T(u,v).$

$\delta_G(x,y,u,v)=\frac{1}{2}(d_G(x,y)+d_G(u,v))
-\frac{1}{2}(d_G(x,u)+d_G(y,v))\leq \frac{1}{2}(d_T(x,y)+d_T(u,v)+
2t) -\frac{1}{2}(d_T(x,u)+d_T(y,v)-  2t)\leq
\frac{1}{2}(d_T(x,u)+d_T(y,v)+ 2t ) -\frac{1}{2}(d_T(x,u)+d_T(y,v)-
2t)= 2t. $
\end{proof}

After showing that the existence of good distance approximating tree
  guarantees low hyperbolicity, in order  to   connect chordality
with hyperbolicity,  we need to make sure  that low chordality
graphs have good distance approximating trees \cite{BCD,CD}  . Here
is an exact result.

\begin{te}  \cite{CD}   Let  $G$ be a $k$-chordal graph. Then, there is a tree
$T$ with $V(T)=V(G)$  such that for any $u,v\in V(G)$ it holds
\begin{equation*} \label{eq:1} |d_G(u,v)-d_T(u,v)|\leq \left\{
\begin{aligned}
          \lfloor
\frac{k}{2 }\rfloor +2, & \ \ \ \ \ \text{if}\ \   k=4,5, \\
               \lfloor  \frac{k}{2 }\rfloor +1, &  \ \ \ \ \ \text{else.}
                          \end{aligned} \right.
                          \end{equation*}
                           \label{Chepoi}
\end{te}

The other  possible approach to connect hyperbolicity and chordality
is via  the concept of the thinness/slimness of geodesic triangles.
This approach also consists of two parts, one is to show that a
graph with low thinness/slimness has low hyperbolicity, as
summarized in \cite[Proposition 1]{CDEHV}, and the other part is to
show that low chordality implies low thinness/slimness.

\begin{figure}
\hspace{-33mm}
\unitlength 1mm 
\linethickness{0.4pt}
\ifx\plotpoint\undefined\newsavebox{\plotpoint}\fi \hspace{30mm}
\begin{center}
\begin{picture}(30,20)
\put(10,5.5){\circle{11}}
\put(0,0){\circle*{1}}\put(0,-5){\makebox(0,0)[b]{$x$}}
\put(20,0){\circle*{1}}\put(20,-5){\makebox(0,0)[b]{$y$}}
\put(0,0){\line(1,0){20}}
\put(10,17.3){\circle*{1}}\put(10,22){\makebox(0,0)[t]{$z$}}
\put(5,8.5){\circle*{1}}\put(-5,8.5){\makebox(0,0)[l]{$m_y^{z,x}$}}
\put(15,8.5){\circle*{1}}\put(25,8.5){\makebox(0,0)[r]{$m_x^{y,z}$}}
\put(10,0){\circle*{1}}\put(10,-5){\makebox(0,0)[b]{$m_z^{x,y}$}}
\qbezier(0,0)(5,8.5)(10,17.3)\qbezier(20,0)(15,8.5)(10,17.3)
 \end{picture}
\end{center}
\caption{An illustration of Gromov product.}\label{fig-1}
\end{figure}
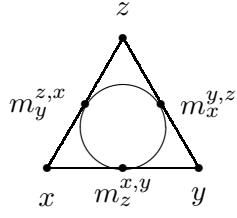

Given a   graph $G,$ we can put an orientation on it by choosing two
maps $\partial_0$ and $\partial_1$ from  $E(G)$ to $V(G)$ such that
each edge $e$ just have $\partial_0(e)$ and $\partial_1(e)$ as its
two endpoints.  The discrete metric space $(V(G),d_G)$ can then be
naturally embedded into the {\em metric graph} \cite[p. 7]{BH}
$(X_G, \widetilde{d_G})$, where $X_G$ is the quotient space of
$E(G)\times [0,1]$ under the identification of $(e,i)$ and $(e',i')$
whenever $\partial _i(e)=\partial_{i'}(e')$ for any $e,e'\in E(G)$
and $i,i'\in \{0,1\},$ and $\widetilde{d_G}$ is the   metric on
$X_G$ satisfying $\widetilde{d_G}((e,t),(e,t'))=
    |t-t'|$ if $ e=e'$ and  $\widetilde{d_G}((e,t),(e,t'))=
\min (d_G(\partial_0(e),\partial_0(e'))+t+t',
d_G(\partial_0(e),\partial_1(e'))+t+1-t',
d_G(\partial_1(e),\partial_0(e'))+1-t+t',d_G(\partial_1(e),\partial_1(e'))+2-t-t'
)$ else. It is easy to see that the definition of  $(X_G,
    \widetilde{d_G})$ is indeed independent of the orientation of
    $G.$  Also, any cycle of $G$ naturally corresponds to a
    circle, namely one-dimensional sphere, embedded in $X_G.$
For any two points   $x,y\in X_G$, there is a not necessarily unique
geodesic connecting them in $(X_G, \widetilde{d_G})$, which we will
use the notation $[x,y]$ if no ambiguity arises.
   We say that $(x,y,z)$ is {\em $(\delta_1,\delta_2)$-thin} provided for any choice
   of the geodesics $[x,y], [y,z], [z,x]$ and  $m_x^{y,z}\in [y,z],m_y^{z,x}\in [z,x],
   m_z^{x,y}\in [x,y]$ satisfying \begin{equation}\left\{ \begin{array}{l }
    \widetilde{d_G}(m_z^{x,y},x)=\widetilde{d_G}(m_y^{z,x},x)=(y\cdot z)_x,\\
\widetilde{d_G}(m_x^{y,z},y)=\widetilde{d_G}(m_z^{x,y},y)=(z\cdot x)_y,\\
   \widetilde{d_G}(m_y^{z,x},z)=\widetilde{d_G}(m_x^{y,z},z)=(x\cdot
   y)_z,
\end{array}\right.\label{E3}
\end{equation}
(Figure \ref{fig-1} is an illustration of \eqref{E3} as well as a
widely-used geometric interpretation of the Gromov product.)  the
following two conditions hold:
   \begin{itemize} \item[(A)]
$\delta_1 \geq \min (\widetilde{d_G}(m_x^{y,z},m_y^{z,x}),
\widetilde{d_G}(m_z^{x,y},   m_y^{z,x}))$;
\item[(B)]   $\{p\in X_G:\ \widetilde{d_G}(p,[y,z]\cup [y,x])\leq \delta_2 \}\supseteq [x,z].$
\end{itemize}
  Modifying  the original definition of Gromov slightly  \cite[p. 8, Definition 1.5]{Alonso} \cite[p.
408]{BH} \cite{Gromov}, we say     that a graph  $G$ is {\em
$(\delta_1,\delta_2)$-thin} provided every triple $(x,y,z)$ of its
vertices is $(\delta_1,\delta_2)$-thin.

\begin{lemma}   Let $G$ be a   graph. If $G$ is
$(\delta_1,\delta_2)$-thin, then it is
$(\delta_1+\delta_2)$-hyperbolic. \label{lemma1}
\end{lemma}
\begin{proof}  The proof is taken from
  \cite[p. 15, (2) implies (5)]{Alonso}.
     It suffices to establish \eqref{EQ} for any  $x,y,u,v\in V(G)$.
By Eq. \eqref{E3} and    Condition (A) for the
$(\delta_1,\delta_2)$-thinness of $(x,u,y)$, we have
\begin{equation}\label{4}(x\cdot y)_u+\delta_1 \geq \min (   \widetilde{d_G}(u,
m_x^{u,y})+ \widetilde{d_G}(m_x^{u,y}, m_u^{y,x}),
  \widetilde{d_G}(u,
m_y^{x,u})+ \widetilde{d_G}(m_y^{x,u}, m_u^{y,x}) ) \geq
\widetilde{d_G}(u, m_u^{y,x}).\end{equation} By Condition (B) for
the  $(\delta_1,\delta_2)$-thinness of $(x,v,y)$, we can suppose,
without loss of generality,  that there is $q\in [y,v]$ such that
\begin{equation}\label{5} \delta_2 \geq
\widetilde{d_G}(m_u^{y,x},q).\end{equation} It follows from
$\widetilde{d_G}(q,v)+  \widetilde{d_G}(q,y)=\widetilde{d_G}(y,v)$,
$\widetilde{d_G}(u,q)+\widetilde{d_G}(q,v)\geq
\widetilde{d_G}(u,v)$, and
$\widetilde{d_G}(u,q)+\widetilde{d_G}(q,y)\geq \widetilde{d_G}(u,y)$
that
\begin{equation}\label{6}\widetilde{d_G}(u,q)\geq (y\cdot v)_u.\end{equation}
We surely  have
\begin{equation}\label{7} \widetilde{d_G}(u,
m_u^{y,x})+ \widetilde{d_G}( m_u^{y,x},q)\geq
\widetilde{d_G}(u,q).
\end{equation}
To complete the proof, we just need to add together \eqref{4},
\eqref{5}, \eqref{6}, and \eqref{7}.
\end{proof}

According to Gromov    \cite{Gromov},  Rips invents the concept of
slimness: For any real number $\delta,$  we say that a graph $G$ is
{\em $\delta$-slim} if for every triple $(x,y,z)$ of vertices of $G,
$ we have
\begin{equation}\label{eq7} \{p\in X_G:\ \widetilde{d_G}(p,[y,z]\cup
[y,x])\leq \delta \}\supseteq [x,z].
\end{equation}
An easy observation is  that a $(\delta_1,\delta_2)$-thin graph is
$\delta_2$-slim.
 It is mentioned in  \cite[Proposition 1]{CDEHV}  that every $\delta$-slim graph is
$8\delta$-hyperbolic.  The next lemma gives a better bound.

\begin{lemma} If a graph is $\delta$-slim, it must be $(2\delta,
\delta)$-thin and hence $3\delta$-hyperbolic.
\end{lemma}
\begin{proof}  By Lemma \ref{lemma1}, our task is to prove that any $\delta$-slim graph $G$ is $(2\delta,
\delta)$-thin. For this purpose, it suffices to deduce $2\delta \geq
\min (\widetilde{d_G}(m_x^{y,z},m_y^{z,x}),
\widetilde{d_G}(m_z^{x,y}, m_y^{z,x}))$ for any triple $(x,y,z)$ of
vertices of $G.$ The following argument is almost word-for-word the
same as that given in \cite[p. 13, (1) implies (3)]{Alonso}. By
\eqref{eq7}, we suppose, as we may, that there is $w\in [y,x]$ such
that $\widetilde{d_G}(m_y^{z,x},w)\leq \delta.$ Observe that
$$\widetilde{d_G}(x,w)\geq
 \widetilde{d_G}(m_y^{z,x},x)-\widetilde{d_G}(m_y^{z,x},w)\geq
 \widetilde{d_G}(m_y^{z,x},x)-\delta=
 \widetilde{d_G}(m_z^{x,y},x)-\delta$$
 and that
$$\widetilde{d_G}(x,w)\leq
 \widetilde{d_G}(m_y^{z,x},x)+\widetilde{d_G}(m_y^{z,x},w)\leq
 \widetilde{d_G}(m_y^{z,x},x)+\delta=
 \widetilde{d_G}(m_z^{x,y},x)+\delta.$$
It then follows $\widetilde{d_G}(m_z^{x,y},w)\leq \delta$ and
henceforth $\widetilde{d_G}(m_z^{x,y},m_y^{z,x})\leq
  \widetilde{d_G}(m_z^{x,y},w)+ \widetilde{d_G}(w,m_y^{z,x})\leq
  2\delta,$ as desired.
\end{proof}

\begin{lemma}
 Every   $k$-chordal graph is
$ (\frac{k}{2},  \frac{k}{2})$-thin.\label{lemma2}
\end{lemma}
\begin{proof}
Consider any triple $(x,y,z)$ of vertices of  $G.$ By an abuse of
notation as usual, denote by $[x,y],[y,z]$ and $[z,x]$ three
geodesic segments joining the corresponding endpoints and put
$m_x^{y,z}\in [y,z],$ $ m_y^{z,x}\in [z,x]$ and $m_z^{x,y}\in [x,y]$
be three points of $X_G$ satisfying Eq. \eqref{E3}. For any
nonnegative number $t\leq (y\cdot z)_x$, there is a unique point $u$
lying in $ [x,y]$ such that $\widetilde{d_G}(u,x)=t$; we use the
notation $(z;x)_t$ for this point $u.$ Similarly, we define
$(y;x)_t$ for any $0\leq t\leq (y\cdot z)_x$ and so on.
 By symmetry, it
suffices to show that $\widetilde{d_G}((z;x)_t,(y;x)_t)\leq
\frac{k}{2} $  for any $0\leq t\leq (y\cdot z)_x$. Take the maximum
$t'\leq t$ such that $(z;x)_{t'}=(y;x)_{t'}$. The case of $t'=t$ is
trivial and so we assume  that $t'<t.$

\paragraph{\sc Case 1:}
There exists  $t''$ such that $(z;x)_{t''}=(y;x)_{t''}$ and
$(y\cdot z)_x\geq t''>t$. We can assume that
$(z;x)_{t'''}\not=(y;x)_{t'''}$  for any $t'<t'''<t''$.

Clearly, walking along $[x,y]$ from  $(z;x)_{t'}$ to $(z;x)_{t''}$
and then go back to $(y;x)_{t'}$ along $[x,z]$ gives rise to a cycle
$C$ in $G.$  This cycle might contain chords. But, surely $C$  has
no chord which connects one point whose distance to $x$ is less than
$t$ to another point whose distance to $x$ is larger than $t$. This
means that $(z;x)_{t}$ and $(y;x)_{t}$
 must appear in a circle in  $X_G$ corresponding to a  chordless cycle of  $G$. Since $G$ is $k$-chordal, we arrive
 at  $\widetilde{d_G}((z;x)_t,(y;x)_t)\leq
\frac{k}{2}$.

\paragraph{\sc Case 2:}  There exists no $t''$ such that $(z;x)_{t''}=(y;x)_{t''}$ and  $(y\cdot z)_x\geq t''>t$.

Let $\Lambda =  \{(z;x)_s, (y;x)_s:\   0\leq s\leq (y\cdot x)_x\}$
  and   $\Upsilon=\{(z;y)_s:\ 0\leq s<   (z\cdot x)_y\}\cup \{(y;z)_s:\ 0\leq s<
(x\cdot y)_z\}.$ For  any  $y\in \Upsilon,$
$\widetilde{d_G}(x,y)>(y\cdot z)_x$ holds and for any $y\in
\Lambda$,    $\widetilde{d_G}(x,y)\leq (y\cdot z)_x$ holds. This
says that
\begin{equation}\label{eq88}\Lambda \cap \Upsilon=\emptyset . \end{equation}
Analogously, by considering both the
distance to $y$ and the distance to  $z$, we
   have  \begin{equation}\label{eq99} \Lambda\cap [y,z]=\emptyset
.\end{equation} Combining Eqs. \eqref{eq88} and \eqref{eq99}, we get
that there is a geodesic    $P$ connecting  $(z;x)_t$ and $(y;x)_t$
whose internal points fall inside  $ [y,z]\cup \Upsilon$.   We
produce a circle in  $X_G$ as follows: Walk along   $P$ from
$(z;x)_t$ to $(y;x)_t$ and then go along $[x,z]$ from $(y;x)_t$ to
$(y;x)_{t'}$ and finally return to $(z;x)_t$ by following $[x,y]$.
This circle naturally corresponds to a cycle of  $G.$ This cycle
might have chords. But for each chord which splits the circle into
two smaller circles, our assumption guarantees that the two vertices
$(z;x)_t$ and $(y;x)_t$ will still appear in one of them
simultaneously. This means that there is a circle of  $X_G$
corresponding to a chordless cycle of  $G$  and  passing from both
$(z;x)_t$ and $(y;x)_t$. As $G$ is $k$-chordal,
$\widetilde{d_G}((z;x)_t,(y;x)_t)\leq \frac{k}{2}$ follows,
 as  expected.
\end{proof}

\section{Proofs}\label{proofs}

\subsection{Lemmas}\label{lemmas}

The proof of our main results, namely Theorems \ref{main} and
\ref{main1},  is divided into a sequence of lemmas/corollaries.

In the course of our proof, we will frequently make  use of the
triangle inequality for the shortest-path metric, namely $ab+bc\geq
ac$, without any claim.  Besides this, we will also freely    apply
the ensuing simple observation, which is so simple that we need not
bother to give any proof here.

\begin{lemma}  \label{EASY} Let  $H$ be a vertex induced subgraph of a graph  $G.$
Then $H$ is an isometric subgraph of  $G$ if and only if
$d_{H}(u,v)=d_G(u,v)$ for each pair of vertices $(u,v)\in V(G)\times
V(G)$ satisfying  $d_H(u,v)\geq 3.$  In particular,  $H$ must be
isometric   if its diameter   is at most $2$.
\end{lemma}

One small matter of convention here and in what follows. When we
refer to a graph, say a graph depicted in Fig. \ref{fig0}, we
sometimes indeed mean that graph together with the special labeling
of its vertices as indicated when it is introduced and sometimes we
mean a graph which is isomorphic to it.    We just leave it to
readers to decide from the context which usage it is. Two immediate
corollaries of Lemma \ref{EASY}  are given subsequently. We state
them with the above convention and omit their routine proofs.

\begin{corollary}\label{lem29}
Let $G$ be a  graph. Let $H\in \{H_1,H_2,H_4\}$ be an induced
subgraph of  $G$ such that $d_G(x,y)=d_G(u,v)=3$. Then $H$ is an
isometric subgraph of  $G$.
\end{corollary}

\begin{corollary}\label{lem30}
Let $G$ be a  graph and $H_3$ be an induced subgraph of  $G$. If
$d_G(x,y)=3$, then $H_3$ is an  isometric  subgraph of  $G$.
\end{corollary}

It is time to  deliver some formal proofs.

\begin{corollary}\label{lem31}
Let $G$ be a  graph and $H_5$ be an induced subgraph of  $G$. If
$d_G(x,y)=d_G(u,v)=3$ and $d_G(b,c)=4$. Then $H_5$ is an isometric
subgraph of  $G$.
\end{corollary}

\begin{proof}
Based on the fact that $d_G(b,c)=4$, we can derive from the triangle
inequality that  $d_G(u,b)=d_G(y,b)=d_G(c,x)=d_G(c,v)=3$. The result
then follows from Lemma \ref{EASY} as $$\{x,y\}, \{ u,v\},  \{
u,b\}, \{ y,b\},  \{ c,x\},   \{ c,v\},  \{ b,c\}$$  are all pairs
inside $V(H_5)\choose 2$ which are of distance at least $3$ apart in
$H_5$.
\end{proof}

\begin{lemma}\label{first}   Let  $G$     be a graph and let $x,y,u,v\in V(G)$. Then $\delta_G
(x,y,u,v)\leq \min (uv,xy,ux,yv,uy,xv)$. \end{lemma}

\begin{proof}   Suppose that
 $d_G(x,S)=d_1, d_G(y,S)=d_2$, where $S=\{ u,v\}$. We can check the
 following:
$$\left\{ \begin{array}{l }
    xy+uv\leq (d_1+d_2+uv)+uv=d_1+d_2+2uv,\\
  d_1+d_2\leq  xu+yv\leq (d_1+uv)+(d_2+uv)= d_1+d_2+2uv,\\
    d_1+d_2\leq xv+yu\leq (d_1+uv)+(d_2+uv)=d_1+d_2+2uv,
\end{array}\right.$$
  from which we get $\delta_G
(x,y,u,v)\leq uv$ and hence  our claim follows by symmetry.
\end{proof}

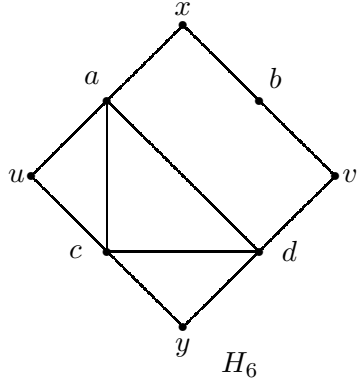
\begin{figure}
\hspace{-33mm}
\unitlength 1mm 
\linethickness{0.4pt}
\ifx\plotpoint\undefined\newsavebox{\plotpoint}\fi \hspace{30mm}
\begin{center}
\begin{picture}(50,40)

\put(30,0){\circle*{1}}\put(30,-4){\makebox(0,0)[b]{$y$}}
\put(10,20){\circle*{1}}\put(7,20){\makebox(0,0)[l]{$u$}}
\put(50,20){\circle*{1}}\put(53,20){\makebox(0,0)[r]{$v$}}
\put(30,40){\circle*{1}}\put(30,43){\makebox(0,0)[t]{$x$}}
\put(20,30){\circle*{1}}\put(17,33){\makebox(0,0)[l]{$a$}}
\put(40,30){\circle*{1}}\put(43,33){\makebox(0,0)[r]{$b$}}
\put(20,10){\circle*{1}}\put(15,10){\makebox(0,0)[l]{$c$}}
\put(40,10){\circle*{1}}\put(45,10){\makebox(0,0)[r]{$d$}}
\put(35,-5){\makebox(0,0)[l]{$H_{6}$}}
\qbezier(10,20)(20,10)(30,0)\qbezier(30,0)(40,10)(50,20)\qbezier(10,20)(20,30)(30,40)\qbezier(30,40)(40,30)(50,20)
\qbezier(20,30)(20,20)(20,10)
\qbezier(20,30)(30,20)(40,10)\qbezier(20,10)(30,10)(40,10)
 \end{picture}
\end{center}
\caption{A  $5$-chordal graph with hyperbolicity $1$.}\label{fig8}
\end{figure}

The next two simple lemmas concern the graph $H_6$ as given in Fig.
\ref{fig8}, which  is obviously a $5$-chordal graph with
hyperbolicity $1.$

\begin{lemma}\label{lead}  Let $H$ be a graph  satisfying
  $V(H)=V(H_2)=V(H_5)=\{x,y,u,v,a,b,c,d\}$ and  $E(H_5)\subseteq E(H)\subseteq E(H_2)$.
Let $t$ be the size  of   $E(H)\cap \{\{a,b\},\{b,d\}, \{d,c\},
\{c,a\}\}.$ If  $t\in \{1,2,3\}$, then  either $H$ contains an
induced $C_4$ or there is an isomorphism from $H$ to $H_6$.
\end{lemma}

\begin{proof}  For any  $v_1,v_2\in V(H)$,  $v_1v_2$  always refers to  $d_H(v_1,v_2)$ in the following.

\paragraph {\sc Case 1:}  $bc=1.$

Let us show  that  $H$  contains an induced $C_4$ in this case.
Since $t\in \{1,2,3\}$, by symmetry, we can assume that either
$ac=1, cd\not=1$ or  $cd=1,bd\not= 1.$
 In the former case,
$[acyd]$ is an induced $4$-cycle of  $H$ and in the latter case
$[cdvb]$
   is an induced $4$-cycle of  $H$.

\paragraph {\sc Case 2:}  $bc\not= 1.$

First observe that replacing the two edges $\{ a,c\}$ and  $\{d,c\}$
by the two new edges $\{ a,b\}$ and  $\{d,b\}$  will transform $H_6$
into another graph which is isomorphic to $H_6.$ Thus, by symmetry,
the condition that $t\in \{1,2,3\}$ means  it is sufficient to
consider the case that $ac=1,cd>1$ and the case that $ac=cd=1,
ab=bd=2.$ For the first case, $[acyd]$ is an induced $4$-cycle of
$H$; for the second case, $H$  itself is exactly $H_6$ after
identifying vertices of the same labels.
\end{proof}

\begin{lemma} \label{lem10}
Suppose that    $G$ is a     $5$-chordal graph  which has $H_6$ as
an induced subgraph.  If  $d_G(x,y)=d_G(u,v)=3$, then  $G$ contains
either $C_4$ or $H_2$ or $H_3$ as an isometric subgraph.
\end{lemma}

\begin{proof}   We can check that  the subgraph of   $H_6$, and hence of  $G,$ induced by $x,a,c,d,v,b$ is
isomorphic to $H_3.$  If $d_G(b,c)=3$,   Corollary  \ref{lem30}
shows that  $G$  contains $H_3$ as an isometric subgraph.  Thus, in
the remaining discussions we will assume that
\begin{equation}  \label{eq2+} d_G(b,c)=2.
\end{equation}

\paragraph {\sc Case 1:} $\min(d_G(b,u),d_G(b,y))=2$.

Assume, as we may, that  $d_G(b,u)=2.$  Take, accordingly, $w\in
V(G)$ satisfying  $d_G(b,w)=d_G(w,u)=1.$ As $d_{H_6}(b,u)=3,$ we see
that  $w\notin V(H_6)$. Observe that
$$2=3-1=d_G(u,v)-d_G(u,w)\leq d_G(v,w)\leq d_G(v,b)+d_G(b,w)=2,$$
which gives \begin{equation}d_G(v,w)=2. \label{China}
\end{equation}

\paragraph {\sc Case 1.1:}  $d_G(w,d)=1$.

In this case, it follows from  Eq.  \eqref{China} that $[wbvd]$  is
an isometric $C_4$ of $G.$

\paragraph {\sc Case 1.2:}   $d_G(w,d)\geq  2.$

Since  $G$ is $5$-chordal, we know that the $6$-cycle  $[wbvdau]$
cannot be chordless in $G.$   By  Eq. \eqref{China}  and the current
assumption that   $d_G(w,d)\geq  2,$  we can  draw the conclusion
that $d_G(w,a)=1 $ and hence find that the  subgraph of  $G$ induced
by $w,u,a,d,v,b$ is isomorphic to $H_3. $  As we already
   assumed that $d_G(u,v)=3$, this induced  $H_3$ is   even an isometric subgraph of
$G,$  taking into account Corollary \ref{lem30}.

\paragraph {\sc Case 2:}  $d_G(b,u)=d_G(b,y)=3$.

By Eq.  \eqref{eq2+},  we can choose  $w\in V(G)$ such that
$d_G(b,w)=d_G(w,c)=1.$   Since  $d_{H_6}(b,c)=3,$  we know that
$w\notin V(H).$ In addition, we have \begin{equation}d_G(u,w)\geq
d_G(u,b)-d_G(b,w)=3-1=2, \ \text{and}\ \ d_G(y,w)\geq
d_G(y,b)-d_G(b,w)=3-1=2. \label{full-moon}
\end{equation}
Clearly, the map which swaps  $u$ and  $y$, $a$ and $d$ and   $x$
and  $v$ is an automorphism of $H_6$ and the requirement to specify
our Case 2 will not be affected after applying this automorphism of
$H_6.$ Therefore, noting  that  $d_G(w,v), d_G(w,d), d_G(w,x),
d_G(w,a) \in \{1,2\},$    we may  take advantage of this  symmetry
of $H_6$ and merely consider
 the following   situations.

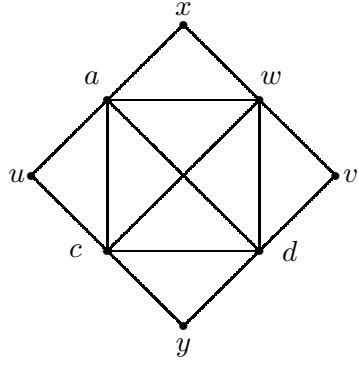
\begin{figure}
\hspace{-33mm}
\unitlength 1mm 
\linethickness{0.5pt}
\ifx\plotpoint\undefined\newsavebox{\plotpoint}\fi \hspace{30mm}
\begin{center}
\begin{picture}(220,50)

\put(100,0){\circle*{1}}\put(100,-4){\makebox(0,0)[b]{$y$}}

\put(80,20){\circle*{1}}\put(77,20){\makebox(0,0)[l]{$u$}}
\put(120,20){\circle*{1}}\put(123,20){\makebox(0,0)[r]{$v$}}

\put(100,40){\circle*{1}}\put(100,43){\makebox(0,0)[t]{$x$}}

\put(90,30){\circle*{1}}\put(87,33){\makebox(0,0)[l]{$a$}}
\put(110,30){\circle*{1}}\put(113,33){\makebox(0,0)[r]{$w$}}
\put(90,10){\circle*{1}}\put(85,10){\makebox(0,0)[l]{$c$}}
\put(110,10){\circle*{1}}\put(115,10){\makebox(0,0)[r]{$d$}}

\qbezier(80,20)(90,10)(100,0)\qbezier(100,0)(110,10)(120,20)\qbezier(80,20)(90,30)(100,40)\qbezier(100,40)(110,30)(120,20)
\qbezier(90,30)(90,20)(90,10)
\qbezier(90,30)(100,20)(110,10)\qbezier(90,10)(100,10)(110,10)\qbezier(90,30)(100,30)(110,30)\qbezier(110,10)(110,20)(110,30)\qbezier(90,10)(100,20)(110,30)
 \end{picture}
\end{center}

\caption{Case 2.1 in the proof of Lemma \ref{lem10}.}\label{fig3}
\end{figure}

\paragraph {\sc Case 2.1:}  $d_G(w,v)=d_G(w,d)=d_G(w,x)=d_G(w,a)=1.$

From Eq. \eqref{full-moon} and our assumption it follows  that
 the subgraph of  $G$ induced by
$x,a,u,c,y,d,v,w$ is isomorphic to $H_2$;  see Fig.  \ref{fig3}.
Because $d_G(x,y)=d_G(u,v)=3$, Corollary  \ref{lem29} now tells us
that $G$ contains   $H_2$ as an isometric subgraph.

\paragraph {\sc Case 2.2:}  $d_G(w,v)=d_G(w,d)=2$.

It is not difficult to check that the subgraph  of   $G$   induced
by $w,b,v,d,c,y$ is isomorphic to $H_3$.  The condition  that
$d_G(b,y)=3$ then enables us to appeal to  Corollary \ref{lem30} and
conclude that $G$ contains the graph $H_3$ as an isometric subgraph.

\paragraph {\sc Case 2.3:}
 $d_G(w,v)=2$ and $ d_G(w,d)=1$.

$G$ contains the induced $4$-cycle  $[wbvd]$.

\paragraph {\sc Case 2.4:}
 $d_G(w,v)=1$ and $ d_G(w,d)=2$.

 $[wvdc]$ is a required induced $C_4$ of  $G.$
\end{proof}

\begin{lemma}\label{lem27}
Let $G$ be a $5$-chordal graph which has   $H_5$ as an induced
subgraph.  If $d_G(x,y)=d_G(u,v)=3$, then $G$ contains at least one
of the subgraphs $C_4,H_2,H_3$ and $H_5$ as an isometric subgraph.
\end{lemma}

\begin{proof}  Because  $H_5$ is an induced subgraph
 of  $G,$  it is clear that   $d_G(b,u),d_G(b,y),d_G(c,x),d_G(c,v)\in
 \{2,3\}$.
There are thus two cases to consider.

\paragraph {\sc Case 1:} $\min  (d_G(b,y),d_G(b,u),d_G(c,x),d_G(c,v))=2$.

Without loss of generality, let us assume that  $d_G(b,y)=2.$  There
is then a vertex  $w$ of   $G$ such that  $d_G(b,w)=d_G(w,y)=1.$
Observe that
\begin{equation}d_G(w,x)\geq d_G(x,y)-d_G(w,y)=3-1=2.
\label{John}
\end{equation}

\paragraph {\sc Case 1.1:}  $d_G(w,a)=1.$

By Eq. \eqref{John}, $[wbxa]$ is an induced $4$-cycle of  $G.$

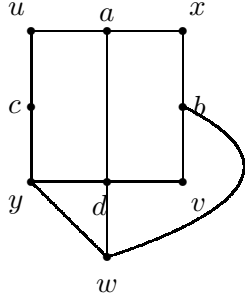
\begin{figure}
\hspace{-33mm}
\unitlength 1mm 
\linethickness{0.5pt}
\ifx\plotpoint\undefined\newsavebox{\plotpoint}\fi \hspace{30mm}
\begin{center}
\begin{picture}(70,50)

\put(10,10){\circle*{1}}\put(7,7){\makebox(0,0)[l]{$y$}}
\put(10,20){\circle*{1}}\put(7,20){\makebox(0,0)[l]{$c$}}
\put(20,30){\circle*{1}}\put(20,33){\makebox(0,0)[t]{$a$}}
\put(30,20){\circle*{1}}\put(33,20){\makebox(0,0)[r]{$b$}}
\put(20,10){\circle*{1}}\put(20,7){\makebox(0,0)[r]{$d$}}
\put(10,30){\circle*{1}}\put(7,33){\makebox(0,0)[l]{$u$}}
\put(30,10){\circle*{1}}\put(33,7){\makebox(0,0)[r]{$v$}}
\put(30,30){\circle*{1}}\put(33,33){\makebox(0,0)[r]{$x$}}
\put(20,0){\circle*{1}}\put(20,-3){\makebox(0,0)[t]{$w$}}

\qbezier(10,30)(10,20)(10,10)\qbezier(10,30)(20,30)(30,30)\qbezier(30,30)(30,20)(30,10)
\qbezier(10,10)(20,10)(30,10)\qbezier(20,10)(20,20)(20,30)\qbezier(20,0)(15,5)(10,10)\qbezier(20,0)(20,5)(20,10)\qbezier(20,0)(50,10)(30,20)
 \end{picture}
\end{center}
\caption{Case 1.2 in the proof of Lemma \ref{lem27}.}\label{fig2.2}
\end{figure}

\paragraph {\sc Case 1.2:}  $d_G(w,a)>1.$

Consider the $6$-cycle  $[bwydax].$  As $G$ is $5$-chordal, this
cycle has a chord in  $G$. According to Eq. \eqref{John} and our
assumption that   $d_G(w,a)>1$, the only possibility is that such a
chord connects $w$ and $d.$ We now examine the subgraph of $G$
induced by $w,b,x,a,d,y$ and realize that it is isomorphic with
$H_3$; see Fig. \ref{fig2.2}. Armed with Corollary \ref{lem30}, our
assumption that $d_G(x,y)=3$ shows that this $H_3$ is even an
isometric subgraph of $G,$ as wanted.

\paragraph {\sc Case 2:} $d_G(b,u)=d_G(b,y)=d_G(c,x)=d_G(c,v)=3$.

By Corollary   \ref{lem31}, $H_5$ is an isometric subgraph of   $G$
provided  $d_G(b,c)=4$. Thus, we shall
    restrict our attention to the cases that $ d_G(b,c)\in \{2, 3\}.$

\paragraph{\sc Case 2.1:} $d_G(b,c)=2$.

Pick a $w\in V(G)$ such that $d_G(b,w)=d_G(w,c)=1.$  It is  clear
that  $[bwcydv]$ is a $6$-cycle in the $5$-chordal graph $G$ and
hence must have a chord. We contend that  this chord can be nothing
but
 $\{w, d\}.$  To see this, one simply needs to notice the following:
\begin{equation*}\left\{ \begin{array}{l }
  d_G(w,y)\geq d_G(b,y)-d_G(b,w)=3-1=2;\\
d_G(w,v)\geq d_G(c,v)-  d_G(c,w)=3-1=2.
\end{array}\right.
\end{equation*}
 From the structure of the subgraph of  $G$ induced by  $b,w,c,y,d,v$  we deduce that    both
 $[cwdy]$ and $[bwdv]$ are  isometric $4$-cycles in  $G$,
 establishing our claim in this case.

\paragraph{\sc Case 2.2:} $d_G(b,c)=3$.

We choose $p,q\in V(G)$ such that $d_G(b,p)=d_G(p,q)=d_G(q,c)=1$. We
first note that
$$d_G(q,v)\geq d_G(c,v)-d_G(c,q)=3-1=2.$$
Due to the symmetry of  $H_5$, it is manifest then that
\begin{equation}  \label{Texas}  \min (d_G(q,v),d_G(p,y),d_G(q,x),d_G(p, u))\geq 2.
\frac{}{}\end{equation}
               We also  have  $ d_G(q,y)\in \{1,2\}$ as it holds
$$2=1+1=d_G(q,c)+d_G(c,y) \geq  d_G(q,y)\geq
d_G(b,y)-d_G(b,q)=3-2=1.$$ Arguing by analogy, we indeed have
\begin{equation} \label{gang}
 d_G(q,y),d_G(q,u),d_G(p,v),d_G(p,x)\in \{1,2\}.
\end{equation}
Eq.  \eqref{Texas} along with Eq. \eqref{gang}  shows that
\begin{equation}\{p,q\}\cap V(H_5)=\emptyset .
\label{xuhui}
\end{equation}

Assume, as we may, that
\begin{equation}d_G(q,y)+d_G(p,v)\geq d_G(q,u)+d_G(p,x).
\label{eqn3}
\end{equation}

\paragraph{\sc Case 2.2.1:} $d_G(q,y)=d_G(p,v)=2$.

We start with two observations:  Thanks to Eq. \eqref{Texas}, we
have $d_G(p,y)\geq 2, d_G(q,v)\geq 2$  while as $b,p,q,c$ is a
geodesic, we obtain $d_G(p,c)=d_G(q,b)=2.$ Now, let us take a look
at the $7$-cycle $[bpqcydv]$  of the $5$-chordal graph
 $G$.  The cycle must have a chord, which, according to our previous observations  and our assumption that  $d_G(q,y)=d_G(p,v)=2$,
  can only be the one  connecting $d$ to $p$ or to  $q.$ Without loss of
generality, let $d_G(q,d)=1$. Then, we can find a $4$-cycle
$[cqdy]$, as desired.

\begin{figure}
\hspace{-33mm}
\unitlength 1mm 
\linethickness{0.5pt}
\ifx\plotpoint\undefined\newsavebox{\plotpoint}\fi \hspace{30mm}
\begin{center}
\begin{picture}(40,30)

\put(10,10){\circle*{1}}\put(7,10){\makebox(0,0)[l]{$y$}}
\put(10,20){\circle*{1}}\put(7,20){\makebox(0,0)[l]{$c$}}
\put(20,30){\circle*{1}}\put(20,33){\makebox(0,0)[t]{$a$}}
\put(30,20){\circle*{1}}\put(33,20){\makebox(0,0)[r]{$b$}}
\put(20,10){\circle*{1}}\put(20,7){\makebox(0,0)[r]{$d$}}
\put(10,30){\circle*{1}}\put(7,33){\makebox(0,0)[l]{$u$}}
\put(30,10){\circle*{1}}\put(33,10){\makebox(0,0)[r]{$v$}}
\put(30,30){\circle*{1}}\put(33,33){\makebox(0,0)[r]{$x$}}
\put(10,0){\circle*{1}}\put(10,-3){\makebox(0,0)[t]{$q$}}
\put(30,0){\circle*{1}}\put(30,-3){\makebox(0,0)[t]{$p$}}
\qbezier(10,30)(10,20)(10,10)\qbezier(10,30)(20,30)(30,30)\qbezier(30,30)(30,20)(30,0)
\qbezier(10,10)(20,10)(30,10)\qbezier(20,10)(20,20)(20,30)\qbezier(10,0)(20,0)(30,0)\qbezier(10,0)(0,10)(10,20)
\qbezier(30,0)(25,5)(20,10)\qbezier(30,0)(40,10)(30,20)
 \end{picture}
\end{center}
\caption{Case 2.2.2 in the proof of Lemma
\ref{lem27}.}\label{fig1.2.2}
\end{figure}
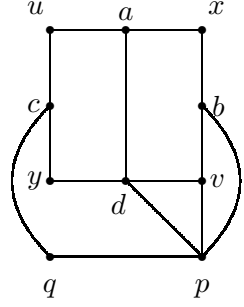

\paragraph{\sc Case 2.2.2:} $d_G(q,y)=2, d_G(p,v)=1$.

In this case, the $5$-chordal graph   $G$ possesses  the  $6$-cycle
 $[pqcydv]$, which must have a chord.  We already assume that  $d_G(q,y)=2;$ as $b,p,q,c$ is a geodesic, we get $d_G(p,c)=2;$ finally, we have
\begin{equation*}\left\{ \begin{array}{l }
  d_G(p,y)\geq d_G(b,y)-d_G(b,p)=3-1=2;\\
d_G(q,v)\geq d_G(c,v)-  d_G(c,q)=3-1=2.
\end{array}\right.
\end{equation*}
Consequently,  it happens either $d_G(q,d)=1$  or  $d_G(p,d)=1$.
 If $d_G(q,d)=1$, we will come to an isometric $4$-cycle  $[cqdy]$.  When $d_G(p,d)=1$ and $d_G(q,d)\geq 2$,   the subgraph of $G$ induced by
 $p,q,c,y,d,v$ is isomorphic to $H_3$  as shown by Fig.  \ref{fig1.2.2}, which is even an isomeric
 subgraph in view of Corollary \ref{lem30} as well as the assumption  that
 $d_G(c,v)=3.$

\paragraph{\sc Case 2.2.3} $d_G(q,y)=1, d_G(p,v)=2$.

 This case can be disposed of as Case 2.2.2.

\paragraph{\sc Case 2.2.4:}   $d_G(q,y)=d_G(p,v)=1.$

By Eqs. \eqref{gang}   and \eqref{eqn3}, we obtain
$d_G(q,u)=d_G(p,x)=1.$ Noting Eq.  \eqref{xuhui}, we further get
\begin{equation}\left\{ \begin{array}{l }
 1\leq d_G(p,d)\leq d_G(p,v)+d_G(v,d)=1+1=2;\\
1\leq d_G(q,d)\leq d_G(q,y)+d_G(y,d)=1+1=2;\\
1\leq d_G(p,a)\leq d_G(p,x)+d_G(x,a)=1+1=2;\\
1\leq d_G(q,a)\leq d_G(q,u)+d_G(u,a)=1+1=2.\\
\end{array}\right.
\label{picb}
\end{equation}
Because of Eq.  \eqref{picb}, it is only necessary to consider the
following three cases, since all others would follow by symmetry.

\paragraph{\sc Case 2.2.4.1:}  $d_G(p,d)=d_G(q,d)=2$.

The subgraph of  $G$ induced by   $c,y,q,p,v,d$ is isomorphic to
$H_3$; see Fig.   \ref{fig1.2.4.1}.  But $d_G(c,v)=3$ is among the
standing assumptions for Case 2  and hence Corollary
     \ref{lem30}  demonstrates that  $G$ has this $H_3$ as an
     isometric subgraph.

\begin{figure}
\hspace{-33mm}
\unitlength 1mm 
\linethickness{0.5pt}
\ifx\plotpoint\undefined\newsavebox{\plotpoint}\fi \hspace{30mm}
\begin{center}
\begin{picture}(40,30)

\put(10,10){\circle*{1}}\put(7,10){\makebox(0,0)[l]{$y$}}
\put(10,20){\circle*{1}}\put(7,20){\makebox(0,0)[l]{$c$}}
\put(20,30){\circle*{1}}\put(20,33){\makebox(0,0)[t]{$a$}}
\put(30,20){\circle*{1}}\put(33,20){\makebox(0,0)[r]{$b$}}
\put(20,10){\circle*{1}}\put(20,7){\makebox(0,0)[r]{$d$}}
\put(10,30){\circle*{1}}\put(7,33){\makebox(0,0)[l]{$u$}}
\put(30,10){\circle*{1}}\put(33,10){\makebox(0,0)[r]{$v$}}
\put(30,30){\circle*{1}}\put(33,33){\makebox(0,0)[r]{$x$}}
\put(10,0){\circle*{1}}\put(10,-3){\makebox(0,0)[t]{$q$}}
\put(30,0){\circle*{1}}\put(30,-3){\makebox(0,0)[t]{$p$}}
\qbezier(10,30)(10,20)(10,0)\qbezier(10,30)(20,30)(30,30)\qbezier(30,30)(30,20)(30,0)
\qbezier(10,10)(20,10)(30,10)\qbezier(20,10)(20,20)(20,30)\qbezier(10,0)(20,0)(30,0)\qbezier(10,0)(-10,15)(10,30)\qbezier(10,0)(0,10)(10,20)
\qbezier(30,0)(50,15)(30,30)\qbezier(30,0)(40,10)(30,20)
 \end{picture}

\end{center}
\caption{Case 2.2.4.1 in the proof of Lemma
\ref{lem27}.}\label{fig1.2.4.1}
\end{figure}

     \paragraph{\sc Case 2.2.4.2:}
$\{d_G(p,d),d_G(q,d)\}=\{ 1, 2\}$.

There is no loss of generality in assuming that $d_G(p,d)=1$  and
$d_G(q,d)=2.$  In such a situation,  by recalling from  Eq.
\eqref{Texas}
     that $d_G(p,y)\geq 2$, we find that   $[pdyq]$  is an induced
$4$-cycle of  $G$, as wanted.

\begin{figure}
\hspace{-33mm}
\unitlength 1mm 
\linethickness{0.5pt}
\ifx\plotpoint\undefined\newsavebox{\plotpoint}\fi \hspace{30mm}
\begin{center}
\begin{picture}(50,50)

\put(20,0){\circle*{1}}\put(20,-4){\makebox(0,0)[b]{$y$}}

\put(0,20){\circle*{1}}\put(-3,20){\makebox(0,0)[l]{$u$}}
\put(40,20){\circle*{1}}\put(43,20){\makebox(0,0)[r]{$v$}}

\put(20,40){\circle*{1}}\put(20,43){\makebox(0,0)[t]{$x$}}

\put(10,30){\circle*{1}}\put(7,33){\makebox(0,0)[l]{$a$}}
\put(30,30){\circle*{1}}\put(33,33){\makebox(0,0)[r]{$p$}}
\put(10,10){\circle*{1}}\put(7,7){\makebox(0,0)[l]{$q$}}
\put(30,10){\circle*{1}}\put(33,7){\makebox(0,0)[r]{$d$}}

\qbezier(0,20)(10,10)(20,0)\qbezier(20,0)(30,10)(40,20)\qbezier(40,20)(30,30)(20,40)\qbezier(20,40)(10,30)(0,20)
\qbezier(10,30)(20,20)(30,10)\qbezier(10,10)(20,20)(30,30)\qbezier(30,30)(30,20)(30,10)\qbezier(30,10)(20,10)(10,10)
\qbezier(10,10)(10,20)(10,30)\qbezier(10,30)(20,30)(30,30)

 \end{picture}

\end{center}
\caption{Case 2.2.4.3 in the proof of Lemma
\ref{lem27}.}\label{fig1.2.4.3}
\end{figure}

   \paragraph{\sc Case 2.2.4.3:}    $d_G(p,d)=d_G(q,d)=d_G(p,a)=d_G(q,a)=1$.

After checking all those  existing  assumptions on pairs of adjacent
vertices  as well as the fact that $\{q,v\},\{ p,y\}\notin E(G)$ as
guaranteed by Eq. \eqref{Texas}, we are led to the conclusion that
the subgraph of $G$ induced by $x,a,u,q,y,d,v,p$ is just $H_2$; see
Fig. \ref{fig1.2.4.3}. Noting further our governing assumption in
the lemma that
 $d_G(x,y)=d_G(u,v)=3$,     Corollary  \ref{lem29}  then enables us reach the conclusion that this $H_2$ is indeed an isometric subgraph
 of  $G,$  as was to be shown.
\end{proof}

The next simple result resembles \cite[Lemma 2.2]{BKM01} closely.

\begin{lemma} \label{lem2.1} Let $G$ be a $k$-chordal graph and let  $C=[x_1\cdots
x_kx_{k+1}\cdots x_{k+t}]$ be  a cycle of $G$.
 If no chord of  $C$ has an endpoint in $\{ x_2,\cdots,x_{k-1}\}$, then $x_1x_{k}=1.$
  \end{lemma}

\begin{proof}   We consider the induced subgraph
$H=G[x_1,x_{k}, x_{k+1},\ldots ,x_{k+t}]$. There must exist a
shortest path in $H$ connecting $x_1$ and  $x_{k}$, say $P$. If the
length of $P$ is greater than $1$, then we walk along $P$ from
$x_{k}$ to $x_1$ and then continue with $x_2,x_3,\ldots,$ and
finally get back to $x_{k}$, creating a chordless cycle of length at
least $k+1,$ which is absurd as $G$ is $k$-chordal.  This proves
that $x_1x_k=1,$ as desired.
\end{proof}

Let $G$ be a graph.   When studying  $\delta _G(x,y,u,v)$ for some
vertices $x,y,u,v$ of $G,$ it is natural to look at a  {\em geodesic
quadrangle}   $\mathcal {Q}(x,u,y,v)$  with {\em corners} $x,u,y$
and $v$, which is just the subgraph of  $G$   induced by the union
of all those vertices on four geodesics connecting $x$ and $u,$ $u$
and  $y$,  $y$ and  $v$, and  $v$ and  $x$, respectively. Let us fix
some notation to be used throughout the paper.

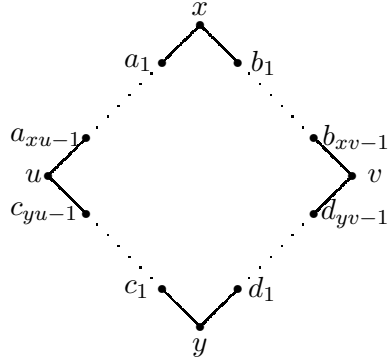
\begin{figure}
\hspace{-33mm}
\unitlength 1mm 
\linethickness{0.4pt}
\ifx\plotpoint\undefined\newsavebox{\plotpoint}\fi \hspace{30mm}
\begin{center}
\begin{picture}(60,60)

\put(20,0){\circle*{1}}\put(20,-4){\makebox(0,0)[b]{$y$}}
\put(0,20){\circle*{1}}\put(-3,20){\makebox(0,0)[l]{$u$}}
\put(40,20){\circle*{1}}\put(42,20){\makebox(0,0)[l]{$v$}}
\put(20,40){\circle*{1}}\put(20,43){\makebox(0,0)[t]{$x$}}

\put(5,25){\circle*{1}}\put(-5,25){\makebox(0,0)[l]{$a_{xu-1}$}}

\put(15,35){\circle*{1}}\put(10,35){\makebox(0,0)[l]{$a_{1}$}}
\put(5,15){\circle*{1}}\put(-5,15){\makebox(0,0)[l]{$c_{yu-1}$}}

\put(15,5){\circle*{1}}\put(10,5){\makebox(0,0)[l]{$c_{1}$}}
\put(25,5){\circle*{1}}\put(30,5){\makebox(0,0)[r]{$d_{1}$}}

\put(35,15){\circle*{1}}\put(45,15){\makebox(0,0)[r]{$d_{yv-1}$}}

\put(25,35){\circle*{1}}\put(30,35){\makebox(0,0)[r]{$b_{1}$}}

\put(35,25){\circle*{1}}\put(45,25){\makebox(0,0)[r]{$b_{xv-1}$}}

\qbezier[7](5,15)(10,10)(15,5)\qbezier(0,20)(2,22)(5,25)\qbezier(0,20)(2,18)(5,15)
\qbezier(15,5)(18,2)(20,0)\qbezier(20,0)(22,2)(25,5)\qbezier(35,15)(38,18)(40,20)
\qbezier[7](25,5)(30,10)(35,15)\qbezier(15,35)(18,38)(20,40)\qbezier(20,40)(22,38)(25,35)
\qbezier(35,25)(38,22)(40,20)\qbezier[7](25,35)(30,30)(35,25)\qbezier[7](5,25)(10,30)(15,35)

 \end{picture}\end{center}
\caption{The geodesic quadrangle   $\mathcal
{Q}(x,u,y,v)$.}\label{fig1}
\end{figure}

\paragraph{\sc \bf  Assumption I:}
 Let  us assume that  $x,u,y,v$ are four different vertices of a graph  $G$ and
  the    four geodesics corresponding to  the  geodesic quadrangle  $\mathcal
 {Q}(x,u,y,v)$  are
\begin{equation*}\left\{ \begin{array}{l }
   P_{a}: x=a_0,a_1,\ldots,a_{xu}=u;\\
 P_b: x=b_0,b_1,\ldots,b_{xv}=v;\\
   P_c: y=c_0,c_1,\ldots,c_{yu}=u;\\
P_d: y=d_0,d_1,\ldots,d_{yv}=v.
\end{array}\right.
\end{equation*}
We call $P_a$, $P_b, P_c$  and $P_d$ the four {\em sides} of
$\mathcal {Q}(x,u,y,v)$  and often just refer to them as vertex
subsets of  $V(G)$  rather than vertex sequences. We write $\mathcal
{P}(x,u,y,v)$ for the pseudo-cycle $$[x, a_1,\ldots , a_{xu-1}, u,
c_{yu-1},\ldots,c_1,y,d_1,\ldots,d_{yv-1},v,b_{xv-1},\ldots, b_1].$$
Note that $\mathcal {P}(x,u,y,v)$ is not necessarily a cycle as the
vertices appearing in the sequence may not be all different. Let us
say that $x$ is {\em opposite} to $P_c$ and $P_d$, say that $x$ and
$y$ are {\em opposite corners}, say that $x$ and $v$ are {\em
adjacent corners}, say that $x$ is the {\em common peak } of $P_a$
and $P_b$, say that $P_a$ and $P_b$ are {\em adjacent} to each
other, say that $P_a$ and $P_d$ are {\em opposite} to each other,
  and say that
    those vertices inside  $P_a\setminus \{x,u\}$  are {\em ordinary vertices} of
    $P_a,$  etc..  An edge of $\mathcal
 {Q}(x,u,y,v)$   which intersects with two adjacent sides      but do not lie in any single  side  is called an $\mathbb{A}$-edge and an edge of
    $\mathcal
 {Q}(x,u,y,v)$ which intersect with two opposite sides but does not lie in any single side is called an
 $\mathbb{H}$-edge.
 Suppose that $a_i=v_0,v_1,\ldots, v_{a_id_j}=d_j $ is a geodesic
 connecting  $a_i$  and  $d_j$ in  $G.$ We call the two
  walks $$x,a_1,\ldots ,a_i,v_1,\ldots, v_{a_id_j-1},  d_j,d_{j-1},\ldots ,d_1,y$$
  and
  $$u,a_{xu-1},\ldots,a_i,v_1,\ldots, v_{a_id_j-1}, d_j,d_{j+1},\ldots,d_{yv-1},v$$
 {\em   $\mathcal
 {Z}$-walks of   $\mathcal
 {Q}(x,u,y,v)$   through $\{a_i,d_j\}$} or just  {\em   $\mathcal
 {Z}$-walks of   $\mathcal
 {Q}(x,u,y,v)$  between $P_a$  and  $P_d.$} In an apparent way, we
 define   similar  concepts for {\em   $\mathcal
 {Z}$-walks of   $\mathcal
 {Q}(x,u,y,v)$    between $P_b$  and  $ P_c$.}


\rz

\begin{lemma} \label{lem11} Let  $G$ be a graph and let    $\mathcal
 {Q}(x,u,y,v)$ be one of its geodesic quadrangles for which Assumption I holds.
 Suppose any two adjacent sides of   $\mathcal
 {Q}(x,u,y,v)$ has only one common vertex and that vertex is  their common peak.
Then $\mathcal
 {Q}(x,u,y,v)$ contains  a   cycle on which $b_1,x, a_1$ appear in that
order consecutively.   Moreover, if $\min(d(P_a,P_d),d(P_b,
P_c))\geq t$ for some $t$, then we may even require that the length
of the cycle is no shorter than $4t$.
\end{lemma}

\begin{proof}
If $\min(d(P_a,P_d),d(P_b, P_c))\geq t$ for $t>0$, then $\mathcal
 {P}(x,u,y,v)$
itself gives rise to  a required cycle.  Otherwise, without loss of
generality, assume that $d(P_a, P_d)=0.$  Take the minimum $i$ such
that $d(a_i, P_d)=0.$  There is a unique $j>0$ such that $a_i=d_j$.
It is plain that  $i>0$  and  $j<yv.$
 This then shows that $[
d_j\cdots d_{yv-1}b_{xv}\cdots b_1xa_1\cdots a_{i-1}]$ is a required
cycle.
   \end{proof}

\begin{lemma}  \cite[p. 67, Claim 2]{BKM01}  We make Assumption I.    Further assume that  \begin{equation}\label{eq1}ub_i=1\end{equation}  for some  $i\geq 1$ and  that
\begin{equation}\label{eq2} xy+uv\geq xv+yu+2.\end{equation} Then, $b_1u<xu.$ \label{koolen2}
\end{lemma}

\begin{proof} We first check the following: \begin{equation*}\begin{array}{cll} xu+uv-2 & \geq  &
(xy-yu)+uv-2\\ & \geq & xv \ \ \ \text{(By Eq. \eqref{eq2})}\\ & = &
xb_i+b_iv\\ & = & (xb_i+b_iu)+(ub_i+b_iv)-2    \ \ \ \text{(By Eq.
\eqref{eq1})}\\
 & \geq & xu+uv-2.
\end{array}  \end{equation*}
Clearly, equalities have to  hold throughout all the above
inequalities. In particular, we have $xu=xb_i+b_iu$. This implies
that there is a geodesic between    $x$ and  $u$ passing through
$b_1$ and hence it is straightforward to see   $b_1u<xu,$ as wanted.
\end{proof}

The next lemma is some variation of  Lemma \ref{first} and will play
an   important role in our short proof of Theorem \ref{main} as to
be presented in Section \ref{Proof}.

\begin{lemma}  Let  $G$ be a graph and let    $\mathcal
 {Q}(x,u,y,v)$ be one of its geodesic quadrangles for which Assumptions I   holds.
If \begin{equation}\label{Ringel} 2\delta_G
(x,y,u,v)=(xy+uv)-\max(xu+yv,xv+yu),
\end{equation}
 then $\delta_G(x,y,u,v)\leq
\min(d(P_a,P_d),d(P_b,P_c))$.\label{lem14}
\end{lemma}

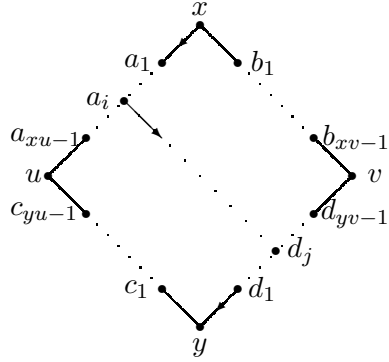
\begin{figure}
\hspace{-33mm}
\unitlength 1mm 
\linethickness{0.4pt}
\ifx\plotpoint\undefined\newsavebox{\plotpoint}\fi \hspace{30mm}
\begin{center}
\begin{picture}(60,60)

\put(20,0){\circle*{1}}\put(20,-4){\makebox(0,0)[b]{$y$}}
\put(0,20){\circle*{1}}\put(-3,20){\makebox(0,0)[l]{$u$}}
\put(40,20){\circle*{1}}\put(42,20){\makebox(0,0)[l]{$v$}}
\put(20,40){\circle*{1}}\put(20,43){\makebox(0,0)[t]{$x$}}
\put(10,30){\circle*{1}}\put(5,30){\makebox(0,0)[l]{$a_{i}$}}
\put(30,10){\circle*{1}}\put(35,10){\makebox(0,0)[r]{$d_{j}$}}
\put(5,25){\circle*{1}}\put(-5,25){\makebox(0,0)[l]{$a_{xu-1}$}}

\put(15,35){\circle*{1}}\put(10,35){\makebox(0,0)[l]{$a_{1}$}}
\put(5,15){\circle*{1}}\put(-5,15){\makebox(0,0)[l]{$c_{yu-1}$}}

\put(15,5){\circle*{1}}\put(10,5){\makebox(0,0)[l]{$c_{1}$}}
\put(25,5){\circle*{1}}\put(30,5){\makebox(0,0)[r]{$d_{1}$}}

\put(35,15){\circle*{1}}\put(45,15){\makebox(0,0)[r]{$d_{yv-1}$}}

\put(25,35){\circle*{1}}\put(30,35){\makebox(0,0)[r]{$b_{1}$}}

\put(35,25){\circle*{1}}\put(45,25){\makebox(0,0)[r]{$b_{xv-1}$}}
\put(10,30){\vector(1,-1){5}}\put(20,40){\vector(-1,-1){3}}
\put(25,5){\vector(-1,-1){3}}

\qbezier[7](5,15)(10,10)(15,5)\qbezier(0,20)(2,22)(5,25)\qbezier(0,20)(2,18)(5,15)
\qbezier(15,5)(18,2)(20,0)\qbezier(20,0)(22,2)(25,5)\qbezier(35,15)(38,18)(40,20)
\qbezier[7](25,5)(30,10)(35,15)\qbezier(15,35)(18,38)(20,40)\qbezier(20,40)(22,38)(25,35)
\qbezier(35,25)(38,22)(40,20)\qbezier[7](25,35)(30,30)(35,25)\qbezier[7](5,25)(10,30)(15,35)

\qbezier[10](10,30)(20,20)(30,10)
 \end{picture}
\end{center}
\caption{$xy\leq i+a_id_j+j$.}\label{fig2}
\end{figure}
\begin{proof} Without loss of generality, we assume that   there exist $i$ and  $j$ such that  \begin{equation}
\label{Kanazawa}a_id_j=\min(d(P_a,P_d),d(P_b,P_c)).
\end{equation}
Focusing on a $\mathcal
 {Z}$-walk  of   $\mathcal
 {Q}(x,u,y,v)$   through $\{a_i,d_j\}$ connecting $x$  and  $y$, we
 find that \begin{equation}xy\leq xa_i+a_id_j+d_jy=i+a_id_j+j;
\label{eq917}
\end{equation} see Fig. \ref{fig2}. Analogously, we have \begin{equation}   uv\leq
(xu-i)+a_id_j+(yv-j). \label{eq917+}\end{equation} Henceforth, we
arrive at the following:
\begin{equation*} \begin{array}{cll} 2\delta
(x,y,u,v)&=&(xy+uv)-\max(xu+yv,xv+yu)\ \ \text{(By Eq.
\eqref{Ringel})}
\\
&\leq& (xy+uv)-(xu+yv)
\\
&\leq& (i+a_id_j+j)+((xu-i)+a_id_j+(yv-j))-(xu+yv)\ \ \text{(By
 Eqs.  \eqref{eq917}  and \eqref{eq917+})   }\\&=& 2a_id_j.\label{eq918}\end{array} \end{equation*}
 Combining  this with Eq. \eqref{Kanazawa},    we finish  the proof of the lemma.
\end{proof}

As regards the inequality asserted in Lemma \ref{lem14}, we need to
say something more for the purpose of deriving Theorem \ref{main1}.

\begin{lemma} \label{lemma41} Let  $G$ be a graph and let    $\mathcal
 {Q}(x,u,y,v)$ be one of its geodesic quadrangles for which both Assumption I and Eq.  \eqref{Ringel} hold.

  (i)  If $\delta (x,y,u,v)=d(P_a, P_d)=a_id_j,$    then
\begin{itemize}
\item any $\mathcal
 {Z}$-walk of   $\mathcal {Q}(x,u,y,v)$   through $\{a_i,d_j\}$
 between $P_a$  and  $P_d$ must be a geodesic, and
 \item  under the additional assumption that
$a_id_j=1$ and $G$ is $5$-chordal,  either   $G$ has an isometric
$4$-cycle or $\{a_i,d_j\}$ is the only edge intersecting both $P_a$
and  $P_d$.
\end{itemize}

(ii)
   If   $\delta (x,y,u,v)=d(P_b, P_c)=b_pc_q ,$   then
\begin{itemize}
\item any $\mathcal
 {Z}$-walk of   $\mathcal {Q}(x,u,y,v)$   through $\{b_p,c_q\}$
 between $P_b$  and  $ P_c$ must be a geodesic, and
 \item
under the additional assumption that $b_pc_q=1$ and $G$ is
$5$-chordal,  either   $G$ has an isometric $4$-cycle or
$\{b_p,c_q\}$ is the only edge intersecting both $P_b$ and  $P_c$.
\end{itemize}
\end{lemma}

\begin{proof} (i) Let us  continue   our discussion  launched in the proof of Lemma
\ref{lem14}.

In the event that    $\delta (x,y,u,v)= d(P_a, P_d)=a_id_j,$ the
equalities in both Eq. \eqref{eq917} and Eq. \eqref{eq917+} must
occur, which clearly shows that any $\mathcal
 {Z}$-walk of   $\mathcal
 {Q}(x,u,y,v)$   through $\{a_i,d_j\}$ between $P_a$  and  $
 P_d$ must be a geodesic.

We now  further assume that $a_id_j=d(P_a, P_d)=1$ and $G$ is
$5$-chordal. Set
$$\mathcal {I}=\{ (k,\ell):\ a_{k}d_{\ell}= 1\}.$$ For any  $(k,\ell)\in  \mathcal {I} $, by considering each
$\mathcal{Z}$-walk through $\{a_k,d_\ell\}$ connecting $x$ and  $y$,
which is a geodesic as we already know, we come to
 \begin{equation}  \mathcal {I}= \{  (k,\ell):\
a_{k}d_{\ell}= 1\} \subseteq \{ (k,\ell):\ k+\ell =xy-1 \} .
\label{eqn8}
\end{equation} Eq. \eqref{eqn8}  means that there exists $(i',j')$
and  $0=t_0<t_1<\cdots <t_{m}$ such that
$$\mathcal {I}=\{  (i'-t_\alpha,j'+t_\alpha):\  \alpha =0,1,\ldots, m  \}.$$
Suppose  that $\{a_i,d_j\}$ is not the only edge intersecting both
$P_a$ and  $P_d$. This means that
 $|\mathcal {I}|-1=m \geq 1$ and then we see that  $$[a_{i'}a_{i'-1}\cdots
a_{i'-t_1}b_{j'+t_1}b_{j'+t_1-1}\cdots b_{j'}]$$ is a chordless
cycle of length $2t_1+2$.  Since  $G$ is $5$-chordal and $t_1$ is a
positive integer, this cycle can only be an isometric $4$-cycle of
$G,$ finishing the proof  of  (i).

(ii)  The proof can be carried out in the same way as that of (i).
\end{proof}

\begin{lemma}  Let $G$ be a graph and we will adopt Assumption  I.
 We choose $j$ to be the
maximum number such that $a_jb_j\leq 1$, $i$ the minimum number such
that $b_id_{yv-xv+i}\leq 1,$   $\ell$ the maximum number such that
$c_\ell d_\ell \leq 1$, and $m$ the minimum number such that
$a_mc_{yu-xu+m}\leq 1.$ Put \begin{equation}\left\{ \begin{array}{l
}
   \pi (b)=i-j+\frac{a_jb_j+b_id_{yv-xv+i}}{2},\\
 \pi(d)=(yv-xv+i)-\ell+\frac{b_id_{yv-xv+i}+c_\ell d_\ell}{2},\\
   \pi (c)=(yu-xu+m)-\ell+\frac{a_mc_{yu-xu+m}+c_\ell d_\ell}{2},\\
   \pi (a)=m-j+\frac{a_jb_j+a_m c_{yu-xu+m}}{2}.
\end{array}\right.
\label{PI}
\end{equation}
If Eq. \eqref{Ringel}   is valid, then $\delta_G(x,y,u,v)\leq \min
(\pi (a), \pi (b), \pi (c), \pi (d))\leq 1+\min ( i-j,
(yv-xv+i)-\ell,
   (yu-xu+m)-\ell, m-j).$
     \label{lem49}
\end{lemma}

\begin{figure}
\hspace{-33mm}
\unitlength 1mm 
\linethickness{0.4pt}
\ifx\plotpoint\undefined\newsavebox{\plotpoint}\fi \hspace{30mm}
\begin{center}
\begin{picture}(60,60)

\put(20,0){\circle*{1}}\put(20,-4){\makebox(0,0)[b]{$y$}}
\put(0,20){\circle*{1}}\put(-3,20){\makebox(0,0)[l]{$u$}}
\put(40,20){\circle*{1}}\put(41,20){\makebox(0,0)[l]{$v$}}
\put(20,40){\circle*{1}}\put(20,43){\makebox(0,0)[t]{$x$}}
\put(12,32){\circle*{1}}\put(7,32){\makebox(0,0)[l]{$a_{j}$}}
\put(32,28){\circle*{1}}\put(36,28){\makebox(0,0)[r]{$b_{i}$}}
\put(5,25){\circle*{1}}\put(-5,25){\makebox(0,0)[l]{$a_{xu-1}$}}
\put(28,32){\circle*{1}}\put(32,32){\makebox(0,0)[r]{$b_{j}$}}
\put(32,12){\circle*{1}}\put(48,12){\makebox(0,0)[r]{$d_{yv-xv+i}$}}
\put(15,35){\circle*{1}}\put(11,35){\makebox(0,0)[l]{$a_{1}$}}
\put(5,15){\circle*{1}}\put(-5,15){\makebox(0,0)[l]{$c_{yu-1}$}}

\put(15,5){\circle*{1}}\put(10,5){\makebox(0,0)[l]{$c_{1}$}}
\put(25,5){\circle*{1}}\put(30,5){\makebox(0,0)[r]{$d_{1}$}}
\put(12,32){\vector(1,0){15}} \put(32,28){\vector(0,-1){15}}
\put(0,20){\vector(1,1){3}}\put(20,40){\vector(1,-1){3}}
\put(25,5){\vector(-1,-1){3}}\put(35,25){\vector(1,-1){3}}
\put(35,15){\circle*{1}}\put(45,15){\makebox(0,0)[r]{$d_{yv-1}$}}

\put(25,35){\circle*{1}}\put(29,35){\makebox(0,0)[br]{$b_{1}$}}

\put(35,25){\circle*{1}}\put(45,25){\makebox(0,0)[r]{$b_{xv-1}$}}

\qbezier[7](5,15)(10,10)(15,5)\qbezier(0,20)(2,22)(5,25)\qbezier(0,20)(2,18)(5,15)
\qbezier(15,5)(18,2)(20,0)\qbezier(20,0)(22,2)(25,5)\qbezier(35,15)(38,18)(40,20)
\qbezier[7](25,5)(30,10)(35,15)\qbezier(15,35)(18,38)(20,40)\qbezier(20,40)(22,38)(25,35)
\qbezier(35,25)(38,22)(40,20)\qbezier[7](25,35)(30,30)(35,25)\qbezier[7](5,25)(10,30)(15,35)

\qbezier(32,12)(32,20)(32,28)\qbezier(12,32)(20,32)(28,32)

 \end{picture}\end{center}
\caption{$xy\leq i+b_id_{yv-xv+i}+(yv-(xv-i))$, $uv\leq
(xu-j)+a_jb_j+(xv-j)$.} \label{fig5}
\end{figure}
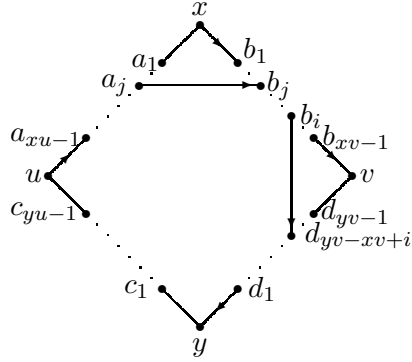
\begin{proof} By symmetry, we only need to show that $\delta_G(x,y,u,v)\leq \pi (b).$

 Taking into account the fact that we can walk from $x$ to  $y$
in $G$ by following $x,b_1,\ldots,b_i$ and then moving in at most
one step from $b_i$ to $d_{yv-xv+i}$ and finally traversing from
$d_{yv-xv+i}$ to $y$ along $P_d$, we get that \begin{equation}xy\leq
i+b_id_{yv-xv+i}+(yv-(xv-i)). \label{jian}\end{equation} Similarly,
starting from $u,$ we can first walk along $P_a$ and then jump from
$a_j$ to $b_j$ in at most one step and then walk along $P_b$ to
arrive at $v$. This gives us \begin{equation}uv\leq
(xu-j)+a_jb_j+(xv-j). \label{guo}\end{equation} See Fig. \ref{fig5}.
Accordingly, we have
\begin{equation}\begin{array}{cll} 2\delta _G(x,y,u,v) & = &
(xy+uv)- \max(xu+yv,xv+yu)   \ \ \ \ \text{(By Eq.  \eqref{Ringel})}\\ & \leq & (xy+uv)-(xu+yv)\\
& \leq & (i+b_id_{yv-xv+i} +(yv-(xv-i)))+((xu-j)+a_jb_j+(xv-j))-(xu+yv)\\
& = & 2\pi (b),
\end{array} \label{eq53}  \end{equation}
which is exactly what we want.
\end{proof}

Brinkmann,  Koolen and Moulton \cite{BKM01} introduced an
extremality argument to deduce upper bounds of hyperbolicity of
graphs. We follow their approach and make the following  standing
assumption in the main step of proving Theorems \ref{main} and
\ref{main1} and thus in several subsequent lemmas.

\paragraph{\sc \bf Assumption II:}  We assume  $x,y,u,v$ are four different vertices of $G$ such that
 the sum $xy+uv$ is minimal
subject to the condition
\begin{equation}\label{key}  xy+uv=\max ( xu+yv,xv+yu )+2\delta^*(G).
\end{equation}

  \rz

\begin{lemma}  \cite[p. 67, Claim 1]{BKM01}   \cite[p. 690, Claim 1]{KM02}  Let $G$ be any graph
and $u,v,x,y\in V(G)$.    Under the Assumptions I and II, we have
            $a_{1}v\geq xv$, $a_{xu-1}y\geq uy$, $b_{1}u\geq xu$, $b_{xv-1}y\geq
vy$, $c_{1}v\geq yv$, $c_{yu-1}x\geq ux$,  $d_{1}u\geq yu$,
$d_{yv-1}x\geq vx$. \label{koolen1}
\end{lemma}

\begin{proof} By symmetry, we only need to show that  $a_{1}v\geq
xv$.  If $a_{1}v< xv$, then,  as a result of  $a_1v\geq
xv-xa_1=xv-1$, we have
\begin{equation}\label{eq4}a_{1}v=xv-1.
\end{equation} Notice the obvious fact that
\begin{equation}\label{eq4+}a_{1}u=xu-1.
\end{equation} We then  come to  the following:
\begin{equation}\begin{array}{cll} a_{1}y+uv & \geq  &
(xy-xa_1)+uv\\ & = & (xy-1)+uv\\& = & (xy+uv)-1\\ & = & \max
(xu+yv-1, xv+yu-1
  )+2\delta^{*}(G)\ \ \ \ \ \ \   \text{(By Eq. \eqref{key})}\\ & = &
\max  ( a_1u+yv, a_1v+yu )+2\delta^{*}(G).\ \ \ \ \ \ \ \ \text{(By
Eqs. \eqref{eq4}  and \eqref{eq4+})}
\end{array}  \label{eq3}\end{equation}

According to the definition of $\delta^{*}(G)$,  we read from  Eq.
\eqref{eq3} that    $a_{1}y+uv= \max  (a_1u+yv, a_1v+yu
  )+2\delta^{*}(G)$   and hence that      $a_{1}y+uv = xy+uv-1$. This
 contrasts with
  the minimality of the sum $xy+uv$ (Assumption II), completing the
proof.
\end{proof}

\begin{corollary} \cite[p. 67, Claim 2]{BKM01}       Under Assumptions I and II and stipulating that
 $\delta^* (G)\geq 1,$   we have that each corner  of  $\mathcal
 {Q}(x,u,y,v)$  is not adjacent to its opposite corner and any ordinary vertex
 of  its  opposite sides and hence has degree $2$ in  $\mathcal
 {Q}(x,u,y,v)$.  \label{cor2.1}  \end{corollary}

\begin{proof}  By symmetry, we only need to prove the claim for the
corner $u$, which
  directly  follows from Lemmas \ref{koolen2}  and \ref{koolen1}.\end{proof}

\begin{lemma} Suppose that Assumptions I and II are met.  (i)  Any two adjacent sides of $\mathcal
 {Q}(x,u,y,v)$ only intersect at their common peak.  In particular,
  no corner of         $\mathcal
 {Q}(x,u,y,v)$  can be an ordinary vertex of some side of  $\mathcal
 {Q}(x,u,y,v)$.   (ii)  For any geodesic  connecting two adjacent corners of $\mathcal
 {Q}(x,u,y,v)$, say  $P=w_0,w_1,\ldots ,w_m$,  no corner of  $\mathcal
 {Q}(x,u,y,v)$ can be found among  $\{w_1,\ldots ,w_{m-1}\}.$    (iii)  Let
$w$ be the common peak of two adjacent sides $P$ and $P'$ of
$\mathcal
 {Q}(x,u,y,v)$.
Let $\alpha \in P\setminus \{w\}$ and $\alpha '\in P'\setminus
\{w\}$ be two vertices of $\mathcal
 {Q}(x,u,y,v)$  such that  $\alpha \alpha '=1$, then
  $\alpha w=\alpha 'w$.
 \label{lem2.4}
\end{lemma}

\begin{proof}
(i)  By symmetry, it suffices to prove that  $a_p\not= b_q$ for any
$p\geq q >0.$ Suppose otherwise, it   then follows   that
$b_1,b_2,\ldots, b_q=a_p, a_{p+1},\ldots,a_{xu}=u$ is a path
connecting $b_1$ and $u$ and so $b_1u<xu ,$ violating  Lemma
\ref{koolen1}.

(ii)  Assume the contrary, we can replace the side of $\mathcal
 {Q}(x,u,y,v)$ connecting the two asserted adjacent corners by the
 geodesic $P$ and get to a new geodesic quadrangle for which Assumptions I and II still
 hold  but for which
a corner   appears as an ordinary vertex in the side $P,$ yielding a
contradiction to (i).

(iii)
 It is no loss to merely   prove that if $i,j>0$ and $a_ib_j=1$
then
 $i=j$.    In the case of   $i>j,$  $b_1,b_2,\ldots, b_j,
a_{i},a_{i+1},\ldots,a_{xu}=u$ is a path connecting $b_1$ and $u$ of
length smaller than $xu ,$ contrary to Lemma \ref{koolen1}.
Similarly, $i<j$
 is impossible as well.
\end{proof}

\begin{corollary} \label{cor45} Let $G$ be a   graph with  $\delta^*(G)>0$  and let $\mathcal
 {Q}(x,u,y,v)$ be a geodesic quadrangle for which  Assumptions  I and II
 hold. Then $\mathcal
 {P}(x,u,y,v)$ is a cycle.  Moreover, if $\delta^*(G)>\frac{1}{2}$,
 then all chords of $\mathcal
 {P}(x,u,y,v)$ must be either $\mathbb{A}$-edges   or $\mathbb{H}$-edges.
\end{corollary}
\begin{proof} This follows from Lemma   \ref{lem14}, Lemma \ref{lem2.4} (i)
and    Corollary \ref{cor2.1}   in a straightforward fashion.
\end{proof}

The next result is very essential to our proof of Theorem
\ref{main1} and both its statement and its  proof have their origin
in  \cite[p. 65, Prop. 3.1]{BKM01} and \cite[p. 691, Claim 2]{KM02}.


\begin{lemma}
Suppose that $G$ is a graph for which Assumptions I and II are met
and $\mathcal
 {Q}(x,u,y,v)$   has  at least one
$\mathbb{A}$-edge. Then we have $xu+yv=xv+yu$.\label{cor12}
\end{lemma}

\begin{proof}
 If the claim were false, without loss of generality, we suppose that
\begin{equation}xu+yv>xv+yu.
\label{eq5}
\end{equation}
By symmetry    and because of  Lemma   \ref{lem2.4}  (iii), let us
work under the assumption  that $a_ib_i=1$. It  clearly holds
\begin{equation}a_i\not= x. \label{clear}
\end{equation}
Before moving on, let us prove that
\begin{equation} a_i\not= u.\label{SJTU}
\end{equation}
Suppose for a contradiction that $a_i=u$, we find that $$
\begin{array}{cll} yv & \geq  &
yu+xv-xu+1\ \ \ \ \ \ \   \text{(By Eq. \eqref{eq5})}\\ & = &
yu+xv-xa_i+1\\
& = &  yu+xv-i+1\\
& = & yu+b_iv+1.
\end{array}
$$
 But we surely have
$yv\leq yu+ub_i+b_iv=yu+b_iv+1$ and so we conclude that we can  get
a geodesic $P$ connecting  $y$ to  $v$ in  $G$  by first walking
along $P_c$ to go from $y$ to $u$, then moving from $u$ to $b_i$ in
one step and finally traversing  along $P_b$ from $b_i$ to $v$.
Since this geodesic passes through $u$ in the middle, we obtain a
contradiction to Lemma \ref{lem2.4} (ii) and hence establish Eq.
\eqref{SJTU}.

To go one step further, let us check the following:
\begin{equation}\begin{array}{cll} xv+1  & =  &  xb_i+b_iv+1= xb_i+b_iv+a_ib_i\\ & \geq  &
xb_i+a_iv=xa_i+a_iv \\
& =  &xa_1+a_1a_i+a_iv \ \ \ \ \ \ \ \text{(By Eq.  \eqref{clear})}\\
& =  & 1+a_1a_i+a_iv\geq
 1+a_1v\\ & \geq  & 1+xv.  \ \ \ \ \ \ \ \text{(By Lemma
\ref{koolen1})}
\end{array}  \label{autumn}\end{equation}
Clearly,    equalities hold throughout Eq.  \eqref{autumn}.
 In particular, we have
\begin{equation}b_iv+1=a_iv.
\label{eqn12}
\end{equation}
 From  Eq. \eqref{eqn12}  and  $xv=b_iv+i$  we   deduce that
\begin{equation}\label{eq8} a_iv=xv-(i-1).\end{equation}
  Here comes the punch line of the proof:
\begin{equation}\begin{array}{cll} a_iy+uv & \geq  &
(xy-xa_i)+uv\\& =  & (xy-i)+uv\\ & = & \max ( xu+yv,xv+yu )
+2\delta^*(G)-i \ \ \ \ \ \ \text{(By Eq. \eqref{key})}
\\ & = &
\max ( (xu-i)+yv, xv+yu-(i-1) ) +2\delta^*(G)   \ \ \ \ \ \
\text{(By Eq. \eqref{eq5})}
\\ & = &
\max
 ((xu-i)+yv,a_iv+yu) +2\delta^*(G)   \ \ \ \ \ \
\text{(By Eq. \eqref{eq8})}
\\ & = & \max
 (a_iu+yv,a_iv+yu) +2\delta^*(G).
\end{array}  \label{eq6}\end{equation}
According to   Eqs. \eqref{clear} and \eqref{SJTU}, we can apply
Lemma \ref{lem2.4} (i) to find that $a_i,y,u,v$ are four different
vertices. We further conclude from the definition of $\delta^*(G)$
that Eq. \eqref{eq6} should hold equalities throughout, hence that
$xy+uv\leq a_iy+uv $ as a result of the minimality  of $xy+uv$
 as indicated in  our Assumption II, and finally that the first inequality in Eq.
\eqref{eq6}   must be strict in light of Eq. \eqref{clear}, getting
a contradiction with the assertion that  all equalities in Eq.
\eqref{eq6}  hold. This is the end of the proof.
\end{proof}

\begin{lemma}  \label{lem15}  Let  $G$ be a graph     for which Assumptions I and II
are required.   Suppose that  $\mathcal
 {Q}(x,u,y,v)$ has  an     $\mathbb{A}$-edge. If there is $1\leq i\leq xu-1$ and
  $0\leq j\leq yv$  such that  $a_id_j\leq 1$, then  $a_id_j= 1$,
 $a_iu+d_jy=yu$ and  $a_ix+d_jv=xv.$
\end{lemma}

\begin{proof} We start with an easy observation:
\begin{equation}\label{EWN} \left\{ \begin{array}{l }
   a_iu+d_jy=a_{xu-1}a_i+1+d_jy\geq a_{xu-1}a_i+a_id_j+d_jy\geq a_{xu-1}y,\\
  a_ix+d_jv=a_{1}a_i+1+d_jv\geq a_{1}a_i+a_id_j+d_jv\geq a_{1}v.
\end{array}\right.
\end{equation}
 In view of    Lemma \ref{koolen1}, this says that
\begin{equation}a_iu+d_jy\geq uy,  \ \    a_ix+d_jv\geq xv.\label{eq13}
\end{equation}
 Adding together the two inequalities in Eq. \eqref{eq13}, we obtain
\begin{equation}\label{eqn14}xu+yv\geq xv+yu.
\end{equation} But, it follows from
 Lemma   \ref{cor12}  and  the
existence of an $\mathbb{A}$-edge of   $\mathcal
 {Q}(x,u,y,v)$ that the equality in Eq. \eqref{eqn14} must occur.
 Consequently, none of the inequalities  in Eqs. \eqref{EWN}  and \eqref{eq13} can be strict, which
 is
 exactly what we want to prove.
\end{proof}

    With a little bit of luck, the forthcoming lemma contributes the number  $ \frac{ \lfloor \frac{k}{2}\rfloor}{2}$, which
    is just the mysterious one we find
in Theorem
 \ref{main}. Note that $ \frac{ \lfloor
 \frac{k}{2}\rfloor}{2}$ is  the smallest half integer that is
 greater than $\frac{k-2}{4}$.

\begin{lemma} \label{lemma54} Let $G$ be a $k$-chordal graph for some  $k\geq 4$ and let $\mathcal
 {Q}(x,u,y,v)$ be a geodesic quadrangle for which  Assumptions  I and II hold.
 Then we have  $\delta^*(G)\leq \frac{ \lfloor \frac{k}{2}\rfloor}{2}$ provided
\begin{equation} \min(d(P_a,P_d),d(P_b,P_c))> 1.\label{eq00}
\end{equation}
\end{lemma}

\begin{proof} Take $i,j,m,\ell$ as specified in Lemma \ref{lem49} and follow all the convention made
in the statement and the proof of  Lemma \ref{lem49}.
 Surely, the result is a direct consequence of   Lemma \ref{lem49}  when
     \begin{equation}\min (\pi (a), \pi (b), \pi (c), \pi (d))\leq     \frac{ \lfloor \frac{k}{2}\rfloor}{2}.
   \label{xiang}\end{equation}

\begin{figure}
\hspace{-33mm}
\unitlength 1mm 
\linethickness{0.4pt}
\ifx\plotpoint\undefined\newsavebox{\plotpoint}\fi \hspace{30mm}
\begin{center}
\begin{picture}(60,60)

\put(20,0){\circle*{1}}\put(20,-4){\makebox(0,0)[b]{$y$}}
\put(0,20){\circle*{1}}\put(-3,20){\makebox(0,0)[l]{$u$}}
\put(40,20){\circle*{1}}\put(42,20){\makebox(0,0)[l]{$v$}}
\put(20,40){\circle*{1}}\put(20,43){\makebox(0,0)[t]{$x$}}

\put(7,27){\circle*{1}}\put(2,27){\makebox(0,0)[l]{$a_{m}$}}

\put(15,35){\circle*{1}}\put(10,35){\makebox(0,0)[l]{$a_{j}$}}

\put(7,13){\circle*{1}}\put(-10,13){\makebox(0,0)[l]{$c_{yu-xu+m}$}}

\put(15,5){\circle*{1}}\put(10,5){\makebox(0,0)[l]{$c_{\ell}$}}
\put(25,5){\circle*{1}}\put(30,5){\makebox(0,0)[r]{$d_{\ell}$}}

\put(33,13){\circle*{1}}\put(35,13){\makebox(0,0)[l]{$d_{yv-xv+i}$}}

\put(25,35){\circle*{1}}\put(30,35){\makebox(0,0)[r]{$b_{j}$}}

\put(33,27){\circle*{1}}\put(38,27){\makebox(0,0)[r]{$b_{i}$}}

\qbezier(15,35)(15,35)(25,35)\qbezier(7,27)(7,13)(7,13)\qbezier(7,27)(7,27)(9,29)
\qbezier[9](9,29)(9,29)(13,33)\qbezier(13,33)(13,33)(15,35)\qbezier(7,13)(7,13)(9,11)
\qbezier(25,5)(25,5)(27,7)\qbezier[9](27,7)(27,7)(31,11)\qbezier[7](0,20)(5,15)(7,13)\qbezier(33,13)(33,13)(31,11)
\qbezier[7](0,20)(2,22)(7,27)\qbezier(15,5)(15,5)(13,7)\qbezier[9](13,7)(9,11)(9,11)\qbezier(15,5)(15,5)(25,5)
\qbezier[7](15,5)(18,2)(20,0)\qbezier[7](20,0)(22,2)(25,5)\qbezier[7](33,13)(38,18)(40,20)
\qbezier(25,5)(25,5)(27,7)\qbezier[9](27,7)(27,7)(31,11)\qbezier[7](15,35)(18,38)(20,40)\qbezier[7](20,40)(22,38)(25,35)
\qbezier[7](33,27)(38,22)(40,20)\qbezier[9](27,33)(27,33)(31,29)\qbezier(25,35)(25,35)(27,33)\qbezier(33,27)(33,27)(31,29)

\qbezier(33,27)(33,20)(33,13)

 \end{picture}\end{center}
\caption{A chordless cycle in   $\mathcal
{Q}(x,u,y,v)$.}\label{figlem48}
\end{figure}
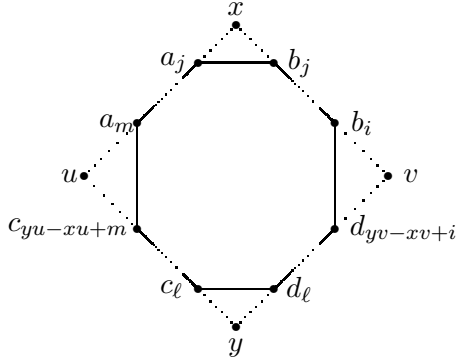

Suppose, for a contradiction,   that the inequality \eqref{xiang}
does not hold. In this event, as $\frac{ \lfloor
\frac{k}{2}\rfloor}{2}\geq 1$, we know that $\min ( i-j,
(yv-xv+i)-\ell,
   (yu-xu+m)-\ell, m-j) \geq \min (\pi(a),\pi(b),\pi(c),\pi(d))-1> 0$.    By virtue of  Lemma   \ref{lem2.4} (i)  and Eq.  \eqref{eq00},
   this implies that
\begin{equation}C=[a_jb_jb_{j+1}\cdots b_id_{yv-xv+i}\cdots
d_{\ell-1}d_{\ell}c_{\ell}c_{\ell+1}\cdots
c_{yu-xu+m}a_ma_{m-1}\cdots a_{j+1}] \label{cycle}\end{equation} is
a cycle, where the redundant $a_j$ should be deleted from the above
notation when $a_j=b_j=x$,  the redundant $b_i$ should be deleted
from the above notation when $b_i=d_{yv-xv+i}=v$, etc.;    see  Fig.
  \ref{figlem48}. Moreover, by Lemma \ref{lem2.4}  (iii), Eq. \eqref{eq00}
and the choice of $i,j,\ell, m,$ we know that $C$ is even a
chordless cycle. But the length of $C$ is just $\pi (a)+\pi (b)+\pi
(c)+\pi (d)$, which, as the assumption is that \eqref{xiang} is
violated, is no smaller than $4(\frac{1}{2}+ \frac{ \lfloor
\frac{k}{2}\rfloor}{2})$ and hence is at least $k+1.$ This
contradicts the assumption that $G$ is $k$-chordal, finishing the
proof.
\end{proof}

\begin{lemma}   Let    $G$ be a  $5$-chordal  graph and
we demand that
 Assumptions I and II hold.
Take  $i,j,\ell,m$  to be the  numbers as specified in Lemma
\ref{lem49}. Suppose that $\mathcal
 {Q}(x,u,y,v)$ has no  $\mathbb{H}$-edges and  $\min(xu,xv,$ $yu,yv$, $2\delta^*(G))$ $\geq
 2$.
Then we have
\begin{equation}a_jb_j+b_id_{yv-xv+i}+c_{\ell}d_{\ell}+a_mc_{yu-xu+m}\geq
2. \label{Nippon}
\end{equation} Furthermore, we have the following conclusions:    if $a_jb_j=1$, then   $(i,m)\in  \{  (j,j), (j, xu), (xv,j)\};$
if $b_id_{yv-xv+i}=1$, then  $(j,\ell)\in \{  (i, yv-xv+i), (i,0),
(0, yv-xv+i) \}$; if  $c_{\ell}d_{\ell}=1,  $    then $(yu-xu+m,
yv-xv+i)\in \{(\ell , \ell), (\ell ,yv), (yu, \ell )$; if  $a_m
c_{yu-xu+m}=1,$ then  $(j, \ell)\in \{  (m, yu-xu+m), (m, 0), (0,
yu-xu+m)\}$. \label{lem92}
\end{lemma}

\begin{proof}
 Since  $\delta^*(G)\geq 1$, it
follows from Lemma \ref{lem49} that
\begin{equation}  \min (i-j,
(yv-xv+i)-\ell, (yu-xu+m)-\ell, m-j)\geq 0. \label{mgs}
\end{equation}
Using Lemma  \ref{lem14} instead, we obtain from $\delta^*(G)\geq 1$
 that $\min (d(P_a,P_d),d(P_b,P_c))\geq 1$.
Considering Lemma   \ref{lem2.4} (i) additionally, this says  that
$G$ has a cycle $C$  as displayed in Eq. \eqref{cycle} whose length
is $\pi (a)+\pi (b)+\pi (c)+\pi (d)$,   where  $\pi (a), \pi (b),
\pi (c), \pi (d)$ stand for  the numbers   introduced in Eq.
\eqref{PI}. As $\mathcal
 {Q}(x,u,y,v)$ has no  $\mathbb{H}$-edges, by the choice of  $i,j,\ell,m$ and by Lemma  \ref{lem2.4} (iii), we see that     $C$ is chordless and are thus led
  to $\pi (a)+\pi (b)+\pi (c)+\pi (d)\leq 5.$

We first observe that
$a_jb_j+b_id_{yv-xv+i}+c_{\ell}d_{\ell}+a_mc_{yu-xu+m}>0$; as
otherwise $C$ will be a chordless cycle of length $xv+xu+yv+yu\geq
8,$ contradicting $\mathbbm{l}\mathbbm{c}(G)\leq 5.$ Let us proceed
to consider the case that
$(a_jb_j,b_id_{yv-xv+i},c_{\ell}d_{\ell},a_mc_{yu-xu+m})=(1,0,0,0)$.
Note that Corollary \ref{cor2.1} guarantees that $0< j
<\min(xu,xv)$. Evoking our assumption $\min (yv,yu)\geq 2$, it is
then obvious that the cycle $C$ contains at least $6$ different
vertices $a_j,b_j,v,d_1,y,c_1,u$, which is absurd as $G$ is
$5$-chordal.  By symmetry, Eq. \eqref{Nippon} is thus established.

Among the four conclusions, let us now only deal with the one
accompanied  with the  assumption   that $a_jb_j=1.$ If
$a_{m}c_{yu-xu+m}=b_id_{yv-xv+i}=0,$ Eq. \eqref{Nippon}
 implies $c_{\ell}d_{\ell}=1$ and so  $C$ will have at least $6$ different vertices
$u,a_j,b_j,v,d_{\ell},c_{\ell}$, contrary to
$\mathbbm{l}\mathbbm{c}(G)\leq 5$. To this point, we can conclude
that $\max (b_id_{yv-xv+i}, a_mc_{yu-xu+m})=1$ and so it suffices to
prove that  $i\in \{ j, xv\}$  and  $m\in \{ j, xu\}.$ We only prove
the first claim and the second one will follow by symmetry. Since
we already have $i\geq j$ as guaranteed by Eq. \eqref{mgs}, our task
is now to get  from $i>j$ to $i=xv.$ If this
 is not true,  the chordless cycle  $C$ will already have four different
 vertices $a_j,b_j,b_i, d_{yv-xv+i}$, which are all outside of $P_c$
 according to Corollary  \ref{cor2.1}. Consequently, due to Corollary  \ref{cor2.1} and $\mathbbm{l}\mathbbm{c}(G)\leq
 5,$ we find that $C$ must have the fifth vertex $c\in P_c\setminus \{ u,y\}$ such
 that   $ca_j=cd_{yv-xv+i}=1.$  In view of  Lemma \ref{lem2.4}
 (iii), we then see that  $$cu=a_ju,cy=
 d_{yv-xv+i}y,a_jx=b_jx,b_iv=d_{yv-xv+1}v.$$ We sum them up and yield
  $xv+yu=xu+yv+b_ib_j>xu+yv$, which is a contradiction with Lemma
  \ref{cor12}.
 This completes the
proof of the lemma.
\end{proof}

\begin{lemma}\label{cor10}  Let $G$ be a  graph for which we will make  Assumptions I and II.
Let $P$ and $P'$    be two adjacent sides of  $\mathcal
 {Q}(x,u,y,v)$  whose common peak is   $w$.
Let $\alpha,  \beta\in P\setminus \{w\}$ and $\alpha', \beta'\in
P'\setminus \{w\}$ be four vertices of $\mathcal
 {Q}(x,u,y,v)$  such that  $\alpha\alpha'=1$ and
  $\beta w=\beta 'w<\alpha w=\alpha 'w$. Then it holds   $\beta \beta
  '=1$ in the case that  $G$ is $4$-chordal as well as in the case that   $G$ is $5$-chordal and  $\beta w=\beta
  'w>1$.
\end{lemma}

\begin{proof}  By symmetry, we only need to show that  for any $i\geq
3$ ($i\geq 2$) we can obtain from $a_ib_i=1$ that $a_{i-1}b_{i-1}=1$
provided $G$ is $5$-chordal ($4$-chordal). But   Lemma \ref{lem2.4}
(i) states that $C=[a_{i-1}a_ib_ib_{i-1}\cdots b_1a_0a_1$ $\cdots
a_{i-2}]$ is a
  cycle of length at
 least $7$ ($5$). Thus,  Lemma
 \ref{lem2.1}  in conjunction with Lemma   \ref{lem2.4} (iii) applies to give
 $a_{i-1}b_{i-1}=1$, as wanted.
\end{proof}

\begin{corollary} \label{cor15}   Let $G$ be a $5$-chordal graph without isometric $C_4$ for which we will make  Assumptions I and
II. If there is an   $\mathbb{A}$-edge connecting $\alpha$ and
$\alpha '$ lying in two adjacent sides $P$ and $P'$ with common peak
$w$, respectively, then   this  is the only $\mathbb{A}$-edge
between  $P$  and  $P'$  and  $\alpha w=\alpha 'w \leq 2.$
\end{corollary}

\begin{proof}
This follows directly from Lemmas \ref{lem2.4} (iii) and
\ref{cor10}.
\end{proof}

\begin{lemma} \label{lem23}        Let  $G$ be a $5$-chordal  graph with $\delta^*(G)\geq 1$ and let  Assumptions I and II hold.
Assume that  $\mathcal {Q}(x,u,y,v)$
 has no $\mathbb{A}$-edges. (i)
  If  there  is an $\mathbb{H}$-edge
 between $P_b$ and  $P_c$, then $\max (xu, yv)\leq 2.$  (ii)  If  there  is an $\mathbb{H}$-edge
 between $P_a$ and  $P_d$, then $\max (xv, yu)\leq 2.$
\end{lemma}

\begin{proof}  We only prove  $xu\leq 2$   under the assumption that   there  is an $\mathbb{H}$-edge
 between $P_b$ and  $P_c$ and all the other claims follow similarly.
Take the minimum $i$ such that $b_i$ is incident with an
$\mathbb{H}$-edge and then pick the maximum  $j$ such that $b_ic_j$
is an $\mathbb{H}$-edge.  By Corollary  \ref{cor2.1}, we have $\min
(i,  yu-j)\geq 1$.
      Since $\mathcal {Q}(x,u,y,v)$
 has no $\mathbb{A}$-edges, we find that  $$[b_1\cdots b_ic_jc_{j+1}\cdots c_{yu-1}ua_{xu-1}\cdots
 a_1x]$$
 is a chordless cycle of length $xu+1+i+(yu-j)\geq xu+3$. Finally, because  $G$
 is
 $5$-chordal, we conclude that $xu\leq 2,$ as desired.
\end{proof}

\begin{lemma}\label{lem19}
Let  $G$ be a  $5$-chordal  graph with   $\delta^*(G)\geq 1$.  We
keep Assumptions I and II. In addition, assume that  $\mathcal
 {Q}(x,u,y,v)$ has a side of length one. Then,  $G$ contains at least one
  graph among $C_4,H_3$ and  $H_5$  as an isometric subgraph.
\end{lemma}

\begin{proof}
It involves  no restriction  of generality in assuming that $xv=1.$
Owning to Lemma  \ref{lem2.4} (i), the walk along $P_a, P_c$ and
$P_d$ will connect $x$ and  $v$ without passing through $x$ or $v$
in the middle and hence there is a shortest path connecting $x$ and
$v$ in the graph obtained from $\mathcal
 {Q}(x,u,y,v)$ by deleting the edge $\{x,v\}$. This
 says that $\mathcal
 {Q}(x,u,y,v)$  has an induced cycle passing through $x$ and  $v$
 contiguously, say  $C=[w_1w_2\cdots w_n]$, where  $w_1=x$ and $w_2=v.$
From Corollary \ref{cor2.1} we know that $w_3=d_{yv-1}\not= a_1=w_n$
and hence  $n>3$. Since $G$ is $5$-chordal, our task is to derive
that if  $n=5$ then  $G$ contains an isometric $C_4,$ $H_3$ or
$H_5$.

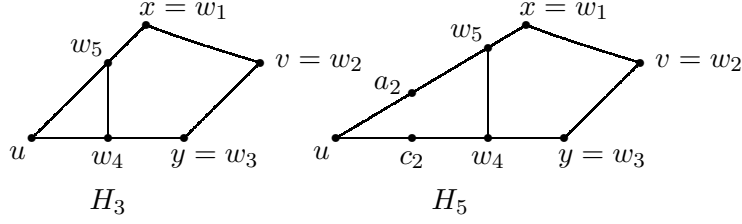
\begin{figure}
\hspace{-33mm}
\unitlength 1mm 
\linethickness{0.4pt}
\ifx\plotpoint\undefined\newsavebox{\plotpoint}\fi \hspace{30mm}
\begin{center}
\begin{picture}(80,25)

\put(20,5){\circle*{1}}\put(24,1){\makebox(0,0)[b]{$y=w_3$}}
\put(10,5){\circle*{1}}\put(10,1){\makebox(0,0)[b]{$w_4$}}
\put(10,15){\circle*{1}}\put(5,17){\makebox(0,0)[l]{$w_5$}}
\put(0,5){\circle*{1}}\put(-3,3){\makebox(0,0)[l]{$u$}}
\put(30,15){\circle*{1}}\put(32,15){\makebox(0,0)[l]{$v=w_2$}}
\put(15,20){\circle*{1}}\put(20,23){\makebox(0,0)[t]{$x=w_1$}}
\put(10,-5){\makebox(0,0)[b]{$H_3$}}
\qbezier(0,5)(10,15)(15,20)\qbezier(0,5)(0,5)(20,5)\qbezier(10,15)(10,10)(10,5)
\qbezier(20,5)(25,10)(30,15)\qbezier(30,15)(20,18)(15,20)

\put(70,5){\circle*{1}}\put(75,1){\makebox(0,0)[b]{$y=w_3$}}
\put(60,5){\circle*{1}}\put(60,1){\makebox(0,0)[b]{$w_4$}}
\put(60,17){\circle*{1}}\put(55,19){\makebox(0,0)[l]{$w_5$}}
\put(40,5){\circle*{1}}\put(37,3){\makebox(0,0)[l]{$u$}}
\put(80,15){\circle*{1}}\put(82,15){\makebox(0,0)[l]{$v=w_2$}}
\put(65,20){\circle*{1}}\put(70,23){\makebox(0,0)[t]{$x=w_1$}}

\put(50,5){\circle*{1}}\put(50,1){\makebox(0,0)[b]{$c_2$}}
\put(50,11){\circle*{1}}\put(45,12){\makebox(0,0)[l]{$a_2$}}

\qbezier(40,5)(60,17)(65,20)\qbezier(40,5)(60,5)(70,5)\qbezier(70,5)(75,10)(80,15)
\qbezier(80,15)(70,18)(65,20)\qbezier(60,5)(60,10)(60,17)
\put(55,-5){\makebox(0,0)[b]{$H_5$}}
 \end{picture}
\end{center}
 \caption{Case 1 in the proof  of Lemma \ref{lem19}.}\label{case1}
\end{figure}

\paragraph{\sc Case 1:}  $w_3$  is a
corner of  $\mathcal
 {Q}(x,u,y,v)$, namely  $yv=1$.

In light of Corollary \ref{cor2.1}, we have  $w_4=c_1$. If $c_1=u$
or  $w_5=u$ occurs,  then   $\mathcal
 {Q}(x,u,y,v) $ turns out to be a $5$-cycle and hence has hyperbolicity $\frac{1}{2}$.  This is impossible as Assumption II means that this hyperbolicity
 can be no smaller than  $\delta^*(G)\geq 1$. Accordingly, by
   Lemma \ref{lem2.4} (iii)   we know that  $w_4u$  and   $w_5u$ have a common value, say  $m.$

   If $m>3 $   or  there are two  $\mathbb{A}$-edges between  $P_a$ and  $P_c$,     Corollary \ref{cor15} says that
 $G$  contains an isometric $C_4.$

When $m=2$ and there are no two $\mathbb{A}$-edges between  $P_a$
and  $P_c$,  the graph $H_5$ as depicted on the right of Fig.
\ref{case1} is an induced graph of $G$. Utilizing Eq.  \eqref{key}
and the assumption that  $\delta^*(G)\geq 1$, we find that
$$4=3+1=ux+xv\geq uv\geq xv+uy+2\delta^*(G)-xy=2+2\delta^*(G)\geq 4.$$ This
illustrates that $uv=4$.  It follows from Lemma \ref{koolen1} that
$a_2y\geq uy=3$. In addition, we have $a_2y\leq
a_2a_1+a_1c_1+c_1y=3$ and so we see that $a_2y=3.$ Similarly, we
have $c_2x=3.$ Getting that $a_2y=c_2x=3$ and $uv=4,$ we apply
  Corollary  \ref{lem31} and conclude that the above-mentioned  $H_5$ must be an isometric
 subgraph of  $G.$

When $m=1,$ the graph $H_3$ as depicted on the left of Fig.
\ref{case1} is an induced graph of $G$.  As in the case of  $m=2$,
we make use of  Eq.  \eqref{key} and  $\delta^*(G)\geq 1$ to get an
important  information: $$3=uw_5+w_5x+xv \geq uv\geq xv+uy
+2\delta^*(G)-xy=1+2\delta^*(G)\geq 3.$$ This implies $uv=3$ and
hence we deduce from    Corollary \ref{lem30} that this $H_3$ is
even an isometric subgraph of $G$.

 \paragraph{\sc Case 2:}  $w_5$  is a
corner of  $\mathcal
 {Q}(x,u,y,v)$, namely $xu=1$.

The analysis is symmetric to that of Case 1.

 \paragraph{\sc Case 3:} Neither $w_3$ nor $w_5$  is a
corner.  In this case, Corollary \ref{cor2.1} ensures that $w_4$ is
not a corner as well.  We proceed to show that this case indeed
cannot happen.

 \paragraph{\sc Case 3.1:}  $w_4
 \in P_a$.

 Since  $P_a$ is a geodesic, we get that $w_4=a_2.$
    It is easy to see that
 $xy\leq vy+vx=vy+1$ and that
$uv=uw_2\leq w_2w_3+w_3w_4+w_4u=2+a_2u=xu.$
 Adding together, we obtain      by Assumption II that
\begin{equation*}   2\delta^*(G)=(xy+uv)-\max(xv+yu,vy+xu)\leq (xy+uv)-(vy+xu)\leq
 1,
\end{equation*}   violating the assumption that $\delta^*(G)\geq 1.$

 \paragraph{\sc Case 3.2:}  $w_4
 \in P_d$.

Reasoning as in Case 3.1   rules out the possibility that this case
may happen.

 \paragraph{\sc Case 3.3:}  $w_4
 \in P_c$.

In this case,  $\mathcal
 {Q}(x,u,y,v)$   contains  $\mathbb{A}$-edges  and hence  Lemma
 \ref{cor12} tells us
 \begin{equation}\label{e33}
  xu+yv=xv+yu.
 \end{equation} But, by
Lemma  \ref{lem2.4} (iii) we have  $xu-1=uw_5=uw_4$ and  $yv-1=w_3y=
w_4y$. We therefore get
  that $xu+yv=2+uw_4+w_4y=2+yu=1+xv+yu,$ which contradicts    Eq. \eqref{e33} and finishes the proof.
\end{proof}

\begin{lemma}     Let  $G$ be a $5$-chordal  graph  and let  Assumptions I and II hold.
If
 $\min (ux,xv,uy,yv,2\delta^*(G))\geq 2$, and $\mathcal {Q}(x,u,y,v)$
 has no $\mathbb{A}$-edges, then $\delta^*(G)= 1$  and   either   $ux=xv=uy=yv=2$  or  $G$ has an isometric $4$-cycle.
 \label{lem2.6}\end{lemma}

\begin{proof} By
Corollary \ref{cor45}, $\mathcal {P}(x,u,y,v)$ is a cycle of length
at least $8$. As $G$ is  $5$-chordal, this cycle must have chords,
which, by Corollary \ref{cor45} again and by the fact that
$\mathcal {Q}(x,u,y,v)$
 has no $\mathbb{A}$-edges, must
be
  $\mathbb{H}$-edges. So, without loss of generality, suppose that  $\mathcal {Q}(x,u,y,v)$
 has an  $\mathbb{H}$-edge between $P_a$ and  $ P_d$.
On the one hand, we can thus go to  Lemma   \ref{lem23} and get
\begin{equation}\label{eq24}   xv=yu=2.
\end{equation}
On the other hand, this allows us to apply
 Lemmas \ref{lem14}  and  \ref{lemma41}  to deduce that  $\delta^*(G)= 1$  and that
 either $G$ has an isometric  $C_4$ or  has exactly  one  $\mathbb{H}$-edge
 between  $P_a$ and  $P_d.$  If the latter case happens, say  we
 have  an  $\mathbb{H}$-edge connecting  $a_i$ and  $d_j,$ we will
 get two chordless cycles of  $G$, $[a_ia_{i-1}\cdots xb_1\cdots b_{xv-1}vd_{yv-1}\cdots d_j]$
and   $[a_ia_{i+1}\cdots uc_{yu-1}\cdots c_1yd_1\cdots d_j]$. Since
neither of these two chordless cycles can be longer than $5$, it
follows from   Eq. \eqref{eq24} that $a_ix+d_jv\leq 2$ and
$ua_i+yd_j\leq 2.$ Taking into account additionally     that $2\leq
ux=ua_i+a_ix$ and $2\leq yv=yd_j+d_jv$, we thus have  $ xu=yv=2.$
This is the end of the proof.
\end{proof}

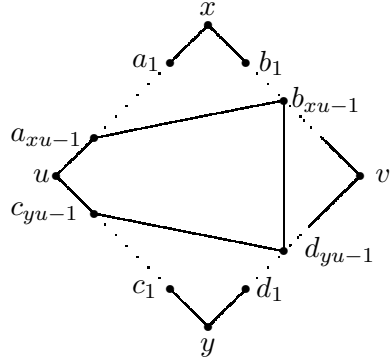
\begin{figure}
\hspace{-33mm}
\unitlength 1mm 
\linethickness{0.4pt}
\ifx\plotpoint\undefined\newsavebox{\plotpoint}\fi \hspace{30mm}
\begin{center}
\begin{picture}(60,60)

\put(20,0){\circle*{1}}\put(20,-4){\makebox(0,0)[b]{$y$}}
\put(0,20){\circle*{1}}\put(-3,20){\makebox(0,0)[l]{$u$}}
\put(40,20){\circle*{1}}\put(42,20){\makebox(0,0)[l]{$v$}}
\put(20,40){\circle*{1}}\put(20,43){\makebox(0,0)[t]{$x$}}

\put(5,25){\circle*{1}}\put(-6,25){\makebox(0,0)[l]{$a_{xu-1}$}}

\put(15,35){\circle*{1}}\put(10,35){\makebox(0,0)[l]{$a_{1}$}}
\put(5,15){\circle*{1}}\put(-6,15){\makebox(0,0)[l]{$c_{yu-1}$}}

\put(15,5){\circle*{1}}\put(10,5){\makebox(0,0)[l]{$c_{1}$}}
\put(25,5){\circle*{1}}\put(30,5){\makebox(0,0)[r]{$d_{1}$}}

\put(30,10){\circle*{1}}\put(42,10){\makebox(0,0)[r]{$d_{yu-1}$}}

\put(25,35){\circle*{1}}\put(30,35){\makebox(0,0)[r]{$b_{1}$}}

\put(30,30){\circle*{1}}\put(40,30){\makebox(0,0)[r]{$b_{xu-1}$}}

\qbezier[7](5,15)(10,10)(15,5)
\qbezier(0,20)(5,25)(5,25)\qbezier(0,20)(5,15)(5,15)
\qbezier(15,5)(18,2)(20,0)\qbezier(20,0)(22,2)(25,5)\qbezier[7](30,30)(40,20)(40,20)\qbezier(35,25)(40,20)(40,20)
\qbezier[7](25,5)(30,10)(35,15)\qbezier(15,35)(18,38)(20,40)\qbezier(20,40)(22,38)(25,35)
\qbezier[7](30,10)(30,10)(40,20)\qbezier(35,15)(30,10)(40,20)\qbezier[7](25,35)(30,30)(35,25)\qbezier[7](5,25)(10,30)(15,35)
\qbezier(5,25)(17,27.5)(30,30)\qbezier(5,15)(15,13)(30,10)\qbezier(30,30)(30,20)(30,10)
 \end{picture}\end{center}
\caption{$[ua_{xu-1}b_{xu-1}d_{yu-1}c_{yu-1}]$ is an induced
$5$-cycle in $\mathcal
 {Q}(x,u,y,v)$.}\label{fig924}
\end{figure}

\begin{lemma}  We take a $5$-chordal graph  $G$ satisfying $\delta^*(G)=1$ and   require  Assumptions I and II. Suppose that  $\mathcal
 {Q}(x,u,y,v)$ has no  $\mathbb{H}$-edge and   $[ua_{xu-1}b_{xu-1}d_{yu-1}c_{yu-1}]$  is an induced $5$-cycle of  $G$; see Fig. \ref{fig924}. Then $G$ has at
 least one graph among $C_4$, $H_3$   and  $H_5$ as an isometric subgraph.
 \label{lem54}
\end{lemma}

\begin{proof}
By  Corollary  \ref{cor2.1} and Lemma  \ref{lem2.4} (iii),
 it will be enough to consider the following  cases,
$b_{xu-1}v=d_{yu-1}v>3 $ or $b_{xu-1}v=d_{yu-1}v\in \{ 1,2\}.$

\paragraph{\sc Case 1:} $b_{xu-1}v=d_{yu-1}v>3.$

Corollary \ref{cor15} implies that $G$ contains an isometric $C_4$.

\paragraph{\sc Case 2:} $b_{xu-1}v=d_{yu-1}v\in \{1,2\}.$

Before jumping into the analysis of two separate subcases, we make
some general observations. Note that
\begin{equation}\begin{array}{cll} xu+yv+2 & = &
(xa_{xu-1}+ua_{xu-1})+(yd_{yu-1}+ d_{yu-1}v)+2
\\ &=&xb_{xu-1}+ua_{xu-1}+yd_{yu-1}+b_{xu-1}v+(a_{xu-1}b_{xu-1}+b_{xu-1}d_{yu-1})
\\ &=& (xb_{xu-1}+b_{xu-1}d_{yu-1}+yd_{yu-1})+(ua_{xu-1}+a_{xu-1}b_{xu-1}+b_{xu-1}v)\\
 &\geq& xy+(ua_{xu-1}+a_{xu-1}v) \\   &\geq& xy+uv \\ &=&\max(xu+yv,xv+yu)+2\delta^*(G)   \ \ \ \ \text{(By Assumption II)}\\
& \geq & xu+yv+2.    \ \ \ \ \text{(By $\delta^*(G)=1$)}
\end{array} \label{eq924}  \end{equation}
It follows that all  inequalities in Eq. \eqref{eq924}  are best
possible     and hence  we have
\begin{equation}uv=ua_{xu-1}+a_{xu-1}b_{xu-1}+b_{xu-1}v=2+b_{xu-1}v\label{925}
\end{equation}
and
\begin{equation}a_{xu-1}v=a_{xu-1}b_{xu-1}+b_{xu-1}v=1+b_{xu-1}v.\label{108}
\end{equation}

\paragraph{\sc Case 2.1:} $b_{xu-1}v=d_{yu-1}v=1$.

We derive from Corollary  \ref{cor2.1} that the subgraph  of  $G$
induced by $u,a_{xu-1},b_{xu-1},v,d_{yu-1},c_{yu-1}$ is isomorphic
to $H_3$ in an obvious way. Thanks to Corollary \ref{lem30}, in
order to check that this  $H_3$ is isometric, our task is to show
that $uv=3.$ But $uv=3$  is an immediate result of
 Eq.  \eqref{925},  proving the claim in this case.

\paragraph{\sc Case 2.2:} $b_{xu-1}v=d_{yu-1}v=2$.

To start things off  we look at  the following:
\begin{equation}\begin{array}{cll} b_{xu-1}v+1 & = &  d_{yu-1}v+1   =d_{yu-1}d_{yv-1}+2  \\ &=&
d_{yu-1}d_{yv-1}+   a_{xu-1}b_{xu-1}+b_{xu-1}d_{yu-1}
\\ &\geq &a_{xu-1}d_{yv-1}   \ \ \ \ \text{(By the triangle inequality)}
\\ &\geq & xd_{yv-1}-xa_{xu-1}  \ \ \ \ \text{(By the triangle inequality)}  \\
&\geq& xv -xa_{xu-1}   \ \ \ \ \text{(By Lemma   \ref{koolen1})}          \\ &=&  (xb_{xu-1}+b_{xu-1}v)-xa_{xu-1} \\
& = & b_{xu-1}v.
\end{array} \label{KIM}  \end{equation}
A consequence of  Eq. \eqref{KIM} is  that
\begin{equation}b_{xu-1}v+1\geq a_{xu-1}d_{yv-1}\geq
b_{xu-1}v.\label{44}
\end{equation}
By symmetry, we also have
\begin{equation}b_{xu-1}v+1= d_{yu-1}v+1  \geq c_{yu-1}b_{xv-1}\geq d_{yu-1}v=b_{xu-1}v.
\label{45}
\end{equation}
As a result of  Eqs.   \eqref{44} and \eqref{45} we
  obtain
 \begin{equation}3\geq \max (a_{xu-1}d_{yv-1},   c_{yu-1}b_{xv-1})\geq \min   (a_{xu-1}d_{yv-1},   c_{yu-1}b_{xv-1}) \geq
  2.\label{eq47}
  \end{equation}
Finally,  let us remark    that  $b_{xu-1}v=2$ implies $xu-1=xv-2$
and hence
  $b_{xu-1}b_{xv-1}=d_{xu-1}d_{xv-1}=1.$

According to Eq. \eqref{eq47}, the following two subcases are
exhaustive.

\paragraph{\sc Case 2.2.1:}
$\min(a_{xu-1}d_{yv-1},c_{yu-1}b_{xv-1})=2$.

 Without loss of
generality, we suppose that there is a vertex $w\in V(G)$ such that
$a_{xu-1}w=wd_{yv-1}=1$. Note that Lemma  \ref{lem2.4} (iii) says
that $w\notin \{a_{xu-1},b_{xu-1},b_{xv-1},v,d_{yv-1}\} $  and hence
    $C= [a_{xu-1}b_{xu-1}b_{xv-1}vd_{yv-1}w]$ is a $6$-cycle in  $G$.
Because $G$ is $5$-chordal, $C$  has at least one chord. Observe
that   Eq.   \eqref{108} says that
\begin{equation} a_{xu-1}v=1+2=3\label{eq57}
\end{equation}
and so
$$wv\geq   a_{xu-1}v-a_{xu-1}w =3-1=2.$$
In consequence, by  virtue of   Lemma  \ref{lem2.4} (iii), we have
  $\min ( wb_{xu-1}, wb_{xv-1},  b_{xv-1}d_{yv-1} )=1.$  There are
  three cases to dwell on.
\paragraph{\sc Case 2.2.1.1:} If  $b_{xv-1}d_{yv-1}=1,$ then
$[b_{xu-1}b_{xv-1}d_{yv-1}d_{xu-1}]$ is a required isometric
$4$-cycle.

\paragraph{\sc Case 2.2.1.2:}  If $wb_{xv-1}=1$ and  $b_{xv-1}d_{yv-1}>1$, we find that
$[b_{xv-1}vd_{yv-1}w]$ is an isometric $4$-cycle, as desired.

\paragraph{\sc Case 2.2.1.3:}
If $\min (wb_{xv-1}, b_{xv-1}d_{yv-1})    >wb_{xu-1}=1$, as a result
of Eq.  \eqref{eq57},
 we can   make use of Corollary \ref{lem30}
to yield that the subgraph induced by
$a_{xu-1},b_{xu-1},b_{xv-1},v,d_{yv-1},w$ is an isometric $H_3$ in
$G$; see Fig. \ref{fig9.24}.

\begin{figure}
\hspace{-33mm}
\unitlength 1mm 
\linethickness{0.4pt}
\ifx\plotpoint\undefined\newsavebox{\plotpoint}\fi \hspace{30mm}
\begin{center}
\begin{picture}(30,20)
\put(0,15){\circle*{1}}\put(-5,18){\makebox(0,0)[l]{$a_{xu-1}$}}
\put(0,0){\circle*{1}}\put(0,-4){\makebox(0,0)[b]{$w$}}
\put(15,5){\circle*{1}}\put(19,5){\makebox(0,0)[r]{$v$}}
\put(10,0){\circle*{1}}\put(10,-4){\makebox(0,0)[b]{$d_{yv-1}$}}

\put(10,10){\circle*{1}}\put(20,10){\makebox(0,0)[r]{$b_{xv-1}$}}
\put(5,15){\circle*{1}}\put(15,15){\makebox(0,0)[r]{$b_{xu-1}$}}

\qbezier(0,15)(0,15)(0,0)\qbezier(0,0)(10,0)(10,0)\qbezier(0,15)(0,15)(5,15)
\qbezier(5,15)(5,15)(15,5)\qbezier(15,5)(15,5)(10,0)\qbezier(0,0)(5,15)(5,15)
 \end{picture}\end{center}
\caption{Case 2.2.1.3 in the proof of Lemma
\ref{lem54}.}\label{fig9.24}
\end{figure}
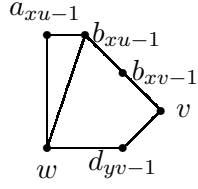

\paragraph{\sc Case 2.2.2:}           $a_{xu-1}d_{yv-1}=c_{yu-1}b_{xv-1}=3$.

From Eq.  \eqref{925} we obtain $uv=4.$ This, together with the
standing assumption of Case 2.2, enables us to deduce from Corollary
\ref{lem31}   that the subgraph induced by
$$u,a_{xu-1},b_{xu-1},b_{xv-1},v,d_{yv-1},d_{yu-1},c_{yu-1}$$
is an isometric $H_5$.
\end{proof}

\subsection{Proofs of Theorems \ref{main} and
\ref{main1}}\label{Proof}

We now have all necessary tools to prove our main results.

\rz

\par \noindent \textbf{Proof of  Theorem \ref{main}: }
Using typical compactness argument, it suffices to prove that every
 connected finite  induced  subgraph of a $k$-chordal graph $G$ is $\frac{\lfloor
\frac{k}{2}\rfloor}{2}$-hyperbolic. If $G$ has less than $4$
vertices, the result is trivial.  Thus, we can simply assume that
$4\leq |V(G)| <\infty $   and henceforth there surely exists
 a
geodesic quadrangle  $\mathcal
 {Q}(x,u,y,v)$ in  $G$  fulfilling
Assumptions I and II. When
 $\min(d(P_a,P_d),d(P_b,P_c))\leq 1$, the result is direct from  Lemma
\ref{lem14} and the fact that  $1\leq \frac{\lfloor
\frac{k}{2}\rfloor}{2}$         while when
 $\min(d(P_a,P_d),d(P_b,P_c))>
1$ we are done by  Lemma \ref{lemma54}. {\QED\par \bigskip \par}

\par \noindent \textbf{Proof of  Theorem \ref{main1}: }
Consider a $5$-chordal graph $G$ with  $\delta^*(G)=1$.   We surely
can get a  geodesic quadrangle  $\mathcal
 {Q}(x,u,y,v)$ in  $G$  for which   Assumption I and Assumption II
 hold. Passing to the proof that $G$ contains one graph from  Fig. \ref{fig0} as an isometric
 subgraph, we have to distinguish four main  cases.

\paragraph {\sc Case 1:}
 $\min (xu,xv,yu,yv)=1. $

  Lemma  \ref{lem19} tells
us that $G$ has either  an isometric $C_4$ or an isometric $H_3$ or
an isometric  $H_5$.

\paragraph {\sc Case 2:}
 $\min (xu,xv,yu,yv)\geq 2$ and there exist   no
$\mathbb{A}$-edges.

\paragraph {\sc Case 2.1:}   $\max (xu,xv,yu,yv)>2$.

By Lemma \ref{lem2.6},  $G$ must have an isometric
 $C_4$.

\paragraph {\sc Case 2.2:}  $xu=xv=yu=yv=2$.

By Corollary \ref{cor45},  $\mathcal
 {Q}(x,u,y,v)$  must have an $\mathbb{H}$-edge. By Corollary
 \ref{cor2.1}, we may assume, without loss of generality, that
 $a_1d_1=1.$  It then follows from Lemma \ref{lemma41} (i) that $xy=uv=3.$

\paragraph {\sc Case 2.2.1:}
 $\mathcal
 {Q}(x,u,y,v)$  has only one  $\mathbb{H}$-edge  and hence the subgraph of  $G$ induced by  its vertices is isomorphic to  $H_5$.

By     Lemma  \ref{lem27}, $G$ has one of $C_4,H_2,H_3$   and  $H_5$
as an isometric subgraph.

\paragraph {\sc Case 2.2.2:}
 $\mathcal
 {Q}(x,u,y,v)$  has  two  $\mathbb{H}$-edges    and hence the subgraph of  $G$ induced by  its vertices is isomorphic to  $H_4$.

By Corollary \ref{lem29},  $G$  contains $H_4$ as an isometric
subgraph.

\paragraph {\sc Case 3:}   $\min (xu,xv,yu,yv)\geq 2$ and there exist no
$\mathbb{H}$-edges.

Take  $i,j,\ell,m$  to be the  numbers as specified in Lemma
\ref{lem49}. By  Lemma \ref{lem92}, Eq. \eqref{Nippon} holds. So,
without loss of generality, we can assume that $i=j,$ $a_jb_j=1$ and
$b_jd_{yv-xv+j}=1.$

\paragraph {\sc Case 3.1:}  $d_{\ell}c_{\ell}=a_mc_{yu-xu+m}=1$.

By Lemma  \ref{lem92}, the chordless cycle displayed  in Eq.
\eqref{cycle} is an isometric $C_4.$

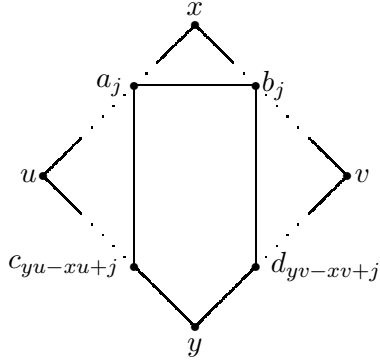
\begin{figure}
\hspace{-33mm}
\unitlength 1mm 
\linethickness{0.4pt}
\ifx\plotpoint\undefined\newsavebox{\plotpoint}\fi \hspace{30mm}
\begin{center}
\begin{picture}(60,60)

\put(20,0){\circle*{1}}\put(20,-4){\makebox(0,0)[b]{$y$}}
\put(0,20){\circle*{1}}\put(-3,20){\makebox(0,0)[l]{$u$}}
\put(40,20){\circle*{1}}\put(41,20){\makebox(0,0)[l]{$v$}}
\put(20,40){\circle*{1}}\put(20,43){\makebox(0,0)[t]{$x$}}
\put(12,32){\circle*{1}}\put(7,32){\makebox(0,0)[l]{$a_{j}$}}

\put(28,32){\circle*{1}}\put(32,32){\makebox(0,0)[r]{$b_{j}$}}
\put(28,8){\circle*{1}}\put(30,8){\makebox(0,0)[l]{$d_{yv-xv+j}$}}

\put(12,8){\circle*{1}}\put(10,8){\makebox(0,0)[r]{$c_{yu-xu+j}$}}

\qbezier(12,32)(12,20)(12,8)
\qbezier[7](5,15)(10,10)(15,5)\qbezier(0,20)(2,22)(5,25)\qbezier(0,20)(2,18)(5,15)
\qbezier(12,8)(18,2)(20,0)\qbezier(20,0)(22,2)(28,8)\qbezier(35,15)(38,18)(40,20)
\qbezier[7](25,5)(30,10)(35,15)\qbezier(15,35)(18,38)(20,40)\qbezier(20,40)(22,38)(25,35)
\qbezier(35,25)(38,22)(40,20)\qbezier[7](25,35)(30,30)(35,25)\qbezier[7](5,25)(10,30)(15,35)

\qbezier(28,32)(28,20)(28,8)\qbezier(12,32)(20,32)(28,32)

 \end{picture}\end{center}
\caption{Case 3.2 in the proof of Theorem \ref{main1}.}
\label{case3.2}
\end{figure}

\paragraph {\sc Case 3.2:}  $(d_{\ell}c_{\ell}, a_mc_{yu-xu+m})=
(0,1)$ or  $(1,0)$.

We only consider the case that  $(d_{\ell}c_{\ell}, a_mc_{yu-xu+m})=
(0,1)$.  For now, the chordless cycle   shown  in Eq. \eqref{cycle}
is just  the  $5$-cycle $[a_jb_jd_{yv-xv+j}yc_{yu-xu+j}];$ see
  Fig.  \ref{case3.2}. Lemma \ref{lem54}  demonstrates that  $G$   contains one
graph among $C_4,H_3$ and  $H_5$ as an isometric subgraph.

\paragraph {\sc Case 3.3:}  $d_{\ell}c_{\ell}=a_mc_{yu-xu+m}=0$.

This case is impossible as the chordless cycle  demonstrated  in Eq.
\eqref{cycle} will contain $6$ different vertices  $a_j,b_j,
d_{yv-xv+j}, y,c_1, u$.

\paragraph {\sc Case 4:}   $\min (xu,xv,yu,yv)\geq 2$ and there exist both $\mathbb{H}$-edges
and $\mathbb{A}$-edges.

Before delving into the case by case analysis, here are some general
observations.
 First note that
 Lemma  \ref{cor12} can be applied to give
 \begin{equation} \label{eq30}  xu+yv=xv+yu.
 \end{equation}
Secondly, according to Corollary \ref{cor2.1},
 we can suppose that
there are \begin{equation}1\leq i\leq xu-1\  \ \text{ and} \ \ 1\leq
j\leq yv-1 \label{eq31}\end{equation} such that $a_id_j=1$ and,
 by Lemma \ref{lem15},  hence  that
\begin{equation}a_iu+d_jy=yu\ \ \text{and} \ \  a_ix+d_jv=xv.
\label{eq23}
\end{equation}
 Thirdly, as    $\delta^*(G)=a_id_j=1$,  Lemma \ref{lem14} gives
 \begin{equation}d(P_a,P_d)=1.
 \label{JAIST}
 \end{equation}
Finally,
  Lemma \ref{lemma41} (i) says that the  $\mathcal
 {Z}$-walks of   $\mathcal
 {Q}(x,u,y,v)$   through the  $\mathbb{H}$-edge  $\{a_i,d_j\}$ must
 be
 geodesics.  Since any subpath of a geodesic is still a geodesic, we
 come to
 \begin{equation}
ud_j=ua_i+a_id_j=ua_i+1 \ \text{and } \ \  a_iy=a_id_j+d_jy=1+d_jy.
\label{eqn33}
\end{equation}

\paragraph {\sc Case 4.1:} $yu=xv=2$.

In this case,  Eq.  \eqref{eq30} forces $xu=xv=yu=yv=2$ and so Eq.
\eqref{eq31} tells us that    $i=j=1.$      It follows that  $\max
(xy, uv)\leq 3$ due to the existence of the path $x,a_1,d_1,y$ and
the path $u,a_1,d_1,v.$ For the moment, in view of Eq. \eqref{key},
we can get  \begin{equation}xy=uv=3.\label{eq33}
\end{equation}

Identifying  $a_1,b_1,c_1,d_1$ with $a,b,c,d,$ respectively,
Corollary \ref{cor2.1} says that  $\mathcal
 {Q}(x,u,y,v)$ is obtained from the graph  $H_5$ as depicted in Fig.   \ref{fig0}  by adding
 $t$ additional edges among
$\{a,b\},\{b,d\},  \{d,c\}, \{c,a\}$, where $t\in \{1,2,3,4\}$, and
adding possibly the edge $\{b,c\}$.

If $t=4$ and $bc=1$,   we easily infer from Eq. \eqref{eq33} and
Corollary \ref{lem29}  that $\mathcal
 {Q}(x,u,y,v)$  is an isometric
subgraph of $G$ which is isomorphic to $H_2$.

If   $t=4$  and $bc>1$, we can check that   $\mathcal
 {Q}(x,u,y,v)$  is an induced
subgraph of $G$ isomorphic to   $H_1$  and then, again    by
  Eq.  \eqref{eq33}
and   Corollary
 \ref{lem29},  $G$  contains an  isometric   $H_1$.

If $t<4$, as a consequence of Lemma \ref{lead},  either $C_4$
 is   an induced subgraph of  $G$ or   $\mathcal
 {Q}(x,u,y,v)$ is isomorphic with $H_6$.          Accordingly,  Eq.  \eqref{eq33} together with    Lemma  \ref{lem10}   implies  that  $G$
 has an isometric  subgraph which is isomorphic to either  $C_4$  or
  $H_2$  or  $H_3.$

   \paragraph {\sc Case 4.2:} $\max (yu, xv) >2$.

We will show that $G$ contains an isometric subgraph which is
isomorphic to    $H_3$, under the assumption that $G$ has no
isometric $C_4$.  Note that the nonexistence of an isometric $C_4$
in $G$  together with Eq.  \eqref{JAIST} yields that there exists
exactly one $\mathbb{H}$-edge between $P_a$ and $P_d$, namely
$\{a_i,d_j\}$, as a result of Lemma \ref{lemma41} (i).

It is no loss of generality in setting
\begin{equation}      yu>2.  \label{eqn35}\end{equation}
 By Lemma \ref{lem2.4} (i)  and Eq.  \eqref{JAIST}, the following is a set   of pairwise different
vertices:
 $$y,c_1,\ldots, c_{yu-1},u,
a_{xu-1},\ldots,a_i,d_j,d_{j-1},\ldots ,d_1.$$ In the subgraph $F$
induced by these  vertices  in  $G$, $a_i$ and $d_j$ are connected
by a path disjoint from the edge $\{a_i,d_j\}.$ This means that
there is a chordless cycle $[w_1w_2\cdots w_n]$ in $F$ where $n\geq
3$ and $w_1 =a_i,w_2=d_j$.   Recall that it is already stipulated
that the $5$-chordal graph $G$ has no isometric $4$-cycle  and hence
$n$ can only take value either $3$ or  $5$.

\begin{figure}
\hspace{-33mm}
\unitlength 1mm 
\linethickness{0.4pt}
\ifx\plotpoint\undefined\newsavebox{\plotpoint}\fi \hspace{30mm}
\begin{center}
\begin{picture}(50,50)

\put(30,0){\circle*{1}}\put(30,-4){\makebox(0,0)[b]{$y$}}

\put(20,0){\circle*{1}}\put(20,-4){\makebox(0,0)[b]{$c_1$}}
\put(10,0){\circle*{1}}\put(10,-4){\makebox(0,0)[b]{$c_2$}}

\put(0,20){\circle*{1}}\put(-4,23){\makebox(0,0)[lt]{$a_i$}}

\put(30,10){\circle*{1}}\put(35,10){\makebox(0,0)[r]{$d_1$}}
\put(30,20){\circle*{1}}\put(35,23){\makebox(0,0)[rt]{$d_2$}}

\qbezier(10,0)(5,10)(0,20)\qbezier(0,20)(10,20)(30,20)\qbezier(30,20)(30,10)(30,0)\qbezier(30,0)(20,0)(10,0)
\qbezier(10,0)(20,10)(30,20)

 \end{picture}
\end{center}

\caption{Case 4.2.1 in the proof of Theorem
\ref{main1}.}\label{fig4.2.1}
\end{figure}
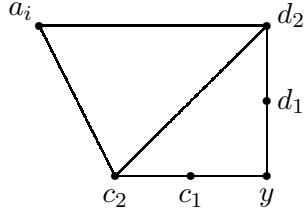

 \paragraph {\sc Case 4.2.1:}  $n=3$.

Since there is exactly one $\mathbb{H}$-edge between $P_a$ and
$P_d$, $w_3$ is neither on
  $P_a$ nor on   $P_d$. Hence,  there is   $0< q < yu$ such that   $w_3=c_q$. It follows from
   Lemma \ref{lem2.4} (iii) that $ua_i=uc_q$ and  $yc_q=yd_j$.
From  Eq. \eqref{eqn35} we obtain $\max (yc_q, c_qu) \geq 2.$
Without loss of generality, assume that $yc_q=\max (yc_q, c_qu) \geq
2$. Since $G$ contains no isometric $C_4$, we infer from Corollary
\ref{cor15} that $q=j=2$ and $c_1d_1=2$. This then demonstrates that
the subgraph induced by the vertices $a_i,d_2,d_1,y,c_1,c_2$ is
 isomorphic    to $H_3$;  see Fig.  \ref{fig4.2.1}.
Granting that $a_iy=3$,
 Corollary  \ref{lem30} will bring to us that $G$ contains $H_3$ as an
 isometric subgraph.
 But $a_iy=3$  follows from  Eq. \eqref{eqn33} and  $d_jy=j=2.$

 \paragraph {\sc Case 4.2.2:}  $n=5$.

  We aim to prove that this case will never happen by deducing contradictions in all the following subcases.

 \paragraph {\sc Case 4.2.2.1:}  Both $w_3$  and  $w_5$ belong to $P_c.$

First consider the case that both $w_3$  and  $w_5$ are ordinary
 vertices of $P_c.$
 From Lemma   \ref{lem2.4} (iii) we obtain $a_iu=uw_5$ and  $d_jy=w_3y.$
 It then follows  $uy=uw_5+w_3y$   by means of Eq. \eqref{eq23}. Since $w_3$ and  $w_5$ are on the same geodesic connecting
 $u$ and  $y$, this is  possible only when
 $w_3=w_5$, yielding a contradiction.

 Next the case that at least one of  $w_3$  and  $w_5$ is a  corner.  We could assume that
  $w_3$ is a corner, and then, in view of Corollary   \ref{cor2.1},
it holds  $w_3=y$. This implies that $w_5\not= u$, as otherwise we
obtain $yu=w_3w_5= 2,$ contradicting Eq.  \eqref{eqn35}.
Accordingly, it follows from  Lemma   \ref{lem2.4} (iii) that
$a_iu=uw_5$. But, we surely have
 $ 2 = w_5w_3=w_5y$  and  $  yd_j=w_3w_2=1.$ Putting together, we get
$a_iu+yd_j=uw_5+1<uw_5+2=uw_5+w_5y=uy$, contradicting Eq.
\eqref{eq23}.

 \paragraph {\sc Case 4.2.2.2:} Neither $w_3$  nor  $w_5$ belongs to $P_c.$

Because there is just one $\mathbb{H}$-edge between $P_a$ and $P_d$,
we see that $w_3,w_4,w_5$ are ordinary vertices of  $P_d, P_c$ and
$P_a$, respectively.  By Lemma   \ref{lem2.4} (iii), it occurs
 that  $a_iu=1+w_5u=1+w_4u$  and  $d_jy=1+w_3y=1+w_4y$. Consequently, we arrive at
 $a_iu+d_jy=2+uy,$ which is contrary to Eq.    \eqref{eq23}.

  \paragraph {\sc Case 4.2.2.3:}  One of  $w_3$  and  $w_5$ is outside of
  $P_c$ and the other lies inside $P_c.$

Incurring no loss of generality,  we make the assumption that
$w_5\notin P_c$ and $w_3\in P_c$.  As there exists only one
$\mathbb{H}$-edge between $P_a$ and $P_d$, we know that
$w_5=a_{i+1}$ and either $w_4=a_{i+2},w_3\not=y$ or $w_4\in
P_c\setminus \{u\}$. By Lemma \ref{lem2.4} (iii), the former implies
that $uy=uw_3+w_3y= uw_4+w_3y= uw_4+yd_j=ua_i+yd_j-2$, violating Eq.
\eqref{eq23}.
 In consequence, we must have $w_4\in P_c\setminus \{u\}$ and hence it
holds either $w_3=c_t,w_4=c_{t+1}$
 or $w_3=c_{t+1},w_4=c_{t}$  for some $t<yu-1$.
If it happens the latter case, we deduce from  Lemma \ref{lem2.4}
(iii) that $uy=uw_4+w_3y-1=uw_5+d_jy-1=ua_i+d_jy-2, $ yielding a
contradiction with the first part of   Eq.  \eqref{eq23}. At this
point, our object is to
 exclude the first possibility as well.
By way of contradiction, let us assume that this case happens and
turn to  the quartet $(w_1,w_2,w_3,u)$. The following calculation
can be trivially verified:
$$uw_1+w_2w_3=uw_1+1;$$
$$\begin{array}{cll} uw_2+w_1w_3 & = &  (uw_1+1) +2  \ \ \ \ \text{(By the first part of Eq.  \eqref{eqn33})}
\\
& = & uw_1+3;
\end{array} $$
$$\begin{array}{cll} uw_3+w_1w_2 & = &  (uw_4+1)+1
\\
& = & (uw_5+1)+1     \ \ \ \ \text{(By Lemma  \ref{lem2.4} (iii))}\\
& = &  uw_1+1.
\end{array} $$
  This gives $\delta
(u, w_1,w_2,w_3)=1=\delta^*(G)$ and $\max (uw_2+w_1w_3, uw_1+w_2w_3,
uw_3+w_1w_2 )=uw_1+3=(uw_1+1)+2\leq ux+yv\leq
xy+uv-2\delta^*(G)=xy+uv-2$, which is the desired contradiction to
Assumption II  on $\mathcal
 {Q}(x,u,y,v)$.
 {\QED\par \bigskip \par}

\end{document}